\documentclass[review,3p,square]{elsarticle}

\usepackage{amssymb}
\usepackage{amsmath}
\usepackage{amsthm}
\usepackage{stmaryrd}
\usepackage{bbm}
\usepackage{enumitem}

\usepackage{fancyhdr}
\usepackage{hyperref}
\hypersetup{%
  pdftitle={BSDEs},%
  pdfsubject={BSDEs},%
  pdfauthor={ShengJun FAN, Ying Hu},%
  pdfkeywords={BSDEs},%
  pdfstartview=FitH,%
  CJKbookmarks=true,%
  bookmarksnumbered=true,%
  bookmarksopen=true,%
  colorlinks=true, linkcolor=blue, urlcolor=blue, citecolor=blue, %
}

\usepackage{cleveref}
\crefname{thm}{Theorem}{Theorems}
\crefname{pro}{Proposition}{Propositions}
\crefname{lem}{Lemma}{Lemmas}
\crefname{rmk}{Remark}{Remarks}
\crefname{cor}{Corollary}{Corollaries}
\crefname{dfn}{Definition}{Definitions}
\crefname{ex}{Example}{Examples}
\crefname{section}{Section}{Sections}
\crefname{subsection}{Subsection}{Subsections}

\newcommand{\eps}{\varepsilon}

\newcommand{\To}{\rightarrow}
\newcommand{\as}{{\rm d}\mathbb{P}\times{\rm d} t\hbox{\rm -}a.e.}
\newcommand{\ass}{{\rm d}\mathbb{P}\times{\rm d} s-a.e.}
\newcommand{\ps}{\mathbb{P}\hbox{\rm -}a.s.}

\newcommand{\F}{\mathcal{F}}
\newcommand{\E}{\mathbb{E}}

\newcommand{\s}{\mathcal{S}}
\newcommand{\lcal}{\mathcal{L}}
\newcommand{\mcal}{\mathcal{M}}

\newcommand{\T}{[0,T]}

\newcommand{\R}{{\mathbb R}}

\newcommand{\RE}{\forall}

\newcommand {\Dis}{\displaystyle}

\newtheorem{thm}{Theorem}[section]
\newtheorem{lem}[thm]{Lemma}
\newtheorem{pro}[thm]{Proposition}
\newtheorem{rmk}[thm]{Remark}

\newtheorem{ex}[thm]{Example}

\journal{ArXiv}
\begin{document}
\begin{frontmatter}

\title{{Existence, uniqueness and comparison theorem on unbounded solutions of scalar super-linear BSDEs}\tnoteref{found}}
\tnotetext[found]{Shengjun Fan is supported by the State Scholarship Fund from the China Scholarship Council (No. 201806425013). Ying Hu is partially supported by Lebesgue Center of Mathematics ``Investissements d'avenir" program-ANR-11-LABX-0020-01, by CAESARS-ANR-15-CE05-0024 and by MFG-ANR-16-CE40-0015-01. Shanjian Tang is supported by National Science Foundation of China (No. 11631004).
\vspace{0.2cm}}


\author{Shengjun Fan\vspace{-0.7cm}\corref{cor1}}
\author{\ \ \ \ Ying Hu\corref{cor2}}
\author{\ \ \ \ Shanjian Tang$^\dag$\corref{cor3}}

\cortext[cor1]{\ School of Mathematics, China University of Mining and Technology, Xuzhou 221116, China. E-mail: f\_s\_j@126.com \vspace{0.2cm}}

\cortext[cor2]{\ Univ. Rennes, CNRS, IRMAR-UMR6625, F-35000, Rennes, France; School of Mathematical Sciences, Fudan University, Shanghai 200433,
China. E-mail: ying.hu@univ-rennes1.fr\vspace{0.2cm}}

\cortext[cor3]{${}^\dag$Department of Finance and Control Sciences, School of Mathematical Sciences, Fudan University, Shanghai 200433, China. E-mail: sjtang@fudan.edu.cn}

\begin{abstract}
This paper is devoted to the existence, uniqueness and comparison theorem on unbounded solutions of  a scalar backward stochastic differential equation (BSDE) whose generator grows (with respect to both unknown variables $y$ and $z$) in a super-linear way like $|y||\ln |y||^{(\lambda+1/2)\wedge 1}+|z||\ln |z||^{\lambda}$ for some $\lambda\geq 0$. For the following four different ranges of the growth power parameter $\lambda$: $\lambda=0$, $\lambda\in (0,1/2)$, $\lambda=1/2$ and $\lambda>1/2$, we give reasonably weakest possible different  integrability conditions of the terminal value for the existence of an unbounded solution to the  BSDE. In the first two cases, they are stronger than the $L\ln L$-integrability and weaker than any $L^p$-integrability with $p>1$; in the third case,  the integrability condition is just some $L^p$-integrability for $p>1$; and in the last case,  the integrability condition is stronger than any $L^p$-integrability with $p>1$ and weaker than any $\exp(L^\eps)$-integrability with $\eps\in (0,1)$. We also establish the comparison theorem, which yields naturally the uniqueness, when either generator of both BSDEs is convex (concave) in both unknown variables $(y,z)$, or satisfies a one-sided Osgood condition in the first unknown variable $y$ and a uniform continuity condition in the second unknown variable $z$. \vspace{0.2cm}
\end{abstract}

\begin{keyword}
Backward stochastic differential equation \sep Existence and uniqueness\sep Super-linear growth\sep Comparison theorem \sep One-sided Osgood condition\sep Uniform continuity condition. \vspace{0.2cm}

\MSC[2010] 60H10\vspace{0.2cm}
\end{keyword}

\end{frontmatter}
\vspace{-0.4cm}

\section{Introduction}
\label{sec:1-Introduction}
\setcounter{equation}{0}

We fix a positive real number $T>0$ and a positive integer $d$, and let $(\Omega, \F, \mathbb{P})$ be a complete probability space, $(B_t)_{t\in\T}$ an $\R^d$-valued standard Brownian motion defined on this space, and $(\F_t)_{t\in\T}$ the natural filtration generated by $B_\cdot$ and augmented by all $\mathbb{P}$-null sets of $\F$. Any progressively measurability with respect to processes will refer to this filtration in this paper.

Let $\R_+$ be the set of all nonnegative real numbers, $x\cdot y$ the scalar inner product of two vectors $x,y\in \R^d$, ${\bf 1}_A$ the indicator function of set $A$, and ${\rm sgn}(x):={\bf 1}_{x>0}-{\bf 1}_{x\leq 0}$. We recall that a real-valued and progressively measurable process $(X_t)_{t\in\T}$ belongs to class (D) if the family of random variables $\{X_\tau: \tau \hbox{ \rm is any $(\F_t)$ -stopping time taking values in } \T\}$ is uniformly integrable.

For any real $p\geq 1$, let $L^p$ denote the set of all real-valued and $\F_T$-measurable random variables $\xi$ satisfying $\E[|\xi|^p]<+\infty$, $\lcal^p$ the set of all real-valued and progressively measurable processes $(X_t)_{t\in\T}$ such that
$$
\|X\|_{\lcal^p}:=\left\{\E\left[\left(\int_0^T |X_t|{\rm d}t\right)^p\right]\right\}^{1/p}<+\infty,\vspace{0.1cm}
$$
$\s^p$ the set of all real-valued, progressively measurable and continuous processes $(Y_t)_{t\in\T}$ satisfying
$$\|Y\|_{{\s}^p}:=\left(\E[\sup_{t\in\T} |Y_t|^p]\right)^{1/p}<+\infty,\vspace{0.1cm}$$
and $\mcal^p$ the set of all $\R^d$-valued and progressively measurable processes $(Z_t)_{t\in\T}$ such that
$$
\|Z\|_{\mcal^p}:=\left\{\E\left[\left(\int_0^T |Z_t|^2{\rm d}t\right)^{p/2}\right] \right\}^{1/p}<+\infty.\vspace{0.1cm}
$$

We consider the following backward stochastic differential equation (BSDE for short):
\begin{equation}\label{eq:1.1}
  Y_t=\xi+\int_t^T g(s,Y_s,Z_s){\rm d}s-\int_t^T Z_s \cdot {\rm d}B_s, \ \ t\in\T,
\end{equation}
where the terminal condition $\xi$ is a real-valued $\F_T$-measurable random variable, and the generator
$$g(\omega, t, y, z):\Omega\times\T\times\R\times\R^d \to \R$$
is jointly continuous in $(y,z)$ and progressively measurable for each $(y,z)$. A solution of \eqref{eq:1.1} is a pair of processes $(Y_t,Z_t)_{t\in\T}$ with values in $\R\times\R^d$, which is progressively measurable such that $\ps$, the functions $t\mapsto Y_t$ is continuous, $t\mapsto Z_t$ is square-integrable, $t\mapsto g(t,Y_t,Z_t)$ is integrable, and Equation~\eqref{eq:1.1} is satisfied. We denote by BSDE $(\xi,g)$  the BSDE with the terminal value $\xi$ and the generator $g$.

In this paper, let us always assume that $\beta,\gamma\geq 0$ are two nonnegative real number and $(\alpha_t)_{t\in\T}$ is an $\R_+$-valued  progressively measurable process. We are interested in the unbounded solution of BSDE \eqref{eq:1.1} with a super-linearly growing generator $g$, which satisfies, roughly speaking, that $\as$,
\begin{equation}\label{eq:1.2}
\RE\ (y,z)\in \R\times\R^d,\quad  |g(\omega,t,y,z)|\leq \alpha_t(\omega)+\beta|y|(\ln |y|)^{\delta}{\bf 1}_{|y|>1}+\gamma |z||\ln |z||^{\lambda}
\end{equation}
for some $(\delta, \lambda)\in (0,1]\times [0,+\infty)$.

The generator $g$ is said to be of a linear growth when $\delta=\lambda=0$.  In this case, if the terminal condition $(\xi,\alpha_\cdot)$ belongs to $L^p\times\lcal^p$ for some $p>1$, then BSDE$(\xi,g)$ admits a solution in the space $\s^p\times\mcal^p$, and the solution is unique in this space if $g$ further satisfies the uniformly Lipschitz condition in $(y,z)$. The reader is referred to \cite{PardouxPeng1990SCL,ElKarouiPengQuenez1997MF,LepeltierSanMartin1997SPL,
BriandDelyonHu2003SPA,FanJiang2012JAMC,Fan2016SPA} for details.

Recently, the authors in \cite{HuTang2018ECP,BuckdahnHuTang2018ECP,FanHu2019ECP} found step by step a weakest integrability condition for $(\xi,\alpha_\cdot)$ to guarantee the existence of a solution to a linearly growing BSDE $(\xi,g)$, where the terminal condition $(\xi,\alpha_\cdot)$ only needs to satisfy an $L\exp\left(\mu\sqrt{2\ln(1+L)}\right)$-integrability condition for a positive parameter $\mu=\gamma \sqrt{T}$. They showed that this integrability condition is weaker than the usual $L^p$-integrability ($p>1$) and stronger than the usual $L\ln L$-integrability, and that the preceding integrability for a positive parameter $\mu<\gamma \sqrt{T}$ is not sufficient to ensure the existence of a solution. They also established the uniqueness of the unbounded solution under the conditions that the generator $g$ satisfies the monotone condition in the state variable $y$ and the uniformly Lipschitz condition in the state variable $z$, and then extended some results mentioned in the above paragraph.

Furthermore, for quadratic BSDEs, namely, the generator $g$ has a quadratic growth in the second unknown variable $z$, much progress has been made in the last two decades. Roughly speaking, to have the existence of a solution to BSDE $(\xi,g)$, the boundedness or at least some exponential integrability of the terminal condition $(\xi,\alpha_\cdot)$ are required. The reader is referred to \cite{Kobylanski2000AP,BriandHu2006PTRF,BriandHu2008PTRF,DelbaenHuRichou2011AIHPPS,
BriandElie2013SPA,BarrieuElKaroui2013AoP,HuTang2015SPA,Fan2016SPA,
FanHuTang2020CRM,FanHu2021SPA} among others. In addition, we mention that a class of quadratic BSDEs with $L^p$-integrable  ($p>1$) terminal conditions are investigated in several recent works, see for example \cite{BahlaliEddahbiOuknine2017AoP,Yanghanlin2017Arxiv,
Bahlali2019Arxiv,BahlaliTangpi2019Arxiv}.

Finally, we would like to especially mention that Bahlali et al.~\cite{BahlaliElAsri2012BSM,BahlaliKebiri2017Stochastics} proved the existence of a solution to BSDE $(\xi,g)$ in the space $\s^p\times\mcal^2$ for some sufficiently large $p>2$, when the terminal condition $(\xi,\alpha_\cdot)$ belongs to $L^p\times\lcal^p$ and the generator $g$ satisfies \eqref{eq:1.2} with $\delta=1$ and $\lambda=1/2$.  They also established the uniqueness of the solution when $g$ further satisfies a locally monotonicity condition in $(y,z)$. Some related works on BSDEs of super-linearly growing generators are available in \cite{BahlaliEssakyHassani2010CRM,BahlaliHakassouOuknine2015Stochastics,
BahlaliEssakyHassani2015SIAM}.

The paper is devoted to the existence, uniqueness and comparison theorem for unbounded solutions of BSDE$ (\xi,g)$ with generator $g$ having a super-linear  growth~\eqref{eq:1.2}. We distinguish four different situations: (i) $\lambda=0$ and $\delta=1/2$, (ii) $\lambda\in (0,1/2)$ and $\delta=\lambda+1/2$, (iii) $\lambda=1/2$ and $\delta=1$, (iv) $\lambda>1/2$ and $\delta=1$, and put forward four reasonably weakest possible integrability conditions of the terminal value $(\xi,\alpha_\cdot)$ to ensure the existence of an unbounded solution  to the  BSDE. In the first two cases, they are stronger than the $L\ln L$-integrability and weaker than any $L^p$-integrability with $p>1$; in the third case,  the integrability condition is just some $L^p$-integrability for $p>1$; and in the last case,  the integrability condition is stronger than any $L^p$-integrability with $p>1$ and weaker than any $\exp(L^\eps)$-integrability with $\eps\in (0,1)$. We also establish the comparison theorem, which yields naturally the uniqueness, when either generators of both BSDEs is convex (concave) in both unknown variables $(y,z)$, or satisfies a one-sided Osgood condition in the first unknown variable $y$ and a uniform continuity condition in the second unknown variable $z$.  We would like to remark that all these results obtained in this paper are not covered by any existing results.

The paper is organized as follows. In next section, we introduce the whole strategy of the paper and establish three useful propositions, which will play an important role later, and in subsequent four sections we establish, respectively for four kinds of different cases mentioned above, the existence, uniqueness and comparison theorem  for unbounded solutions of BSDEs in four kinds of different spaces. In the last section, we conclude the whole paper.

\section{The whole strategy and three useful propositions}
\label{sec:2-whole-idea}
\setcounter{equation}{0}

\subsection{Existence}

For the existence of an unbounded solution to BSDE $(\xi,g)$ with generator $g$ satisfying \eqref{eq:1.2}, our whole idea is to establish some uniform a priori estimate and then apply the localization procedure put forward initially in \cite{BriandHu2006PTRF}. To this end, we first establish a general stability result for the solutions of BSDEs.

We assume that $\xi$ is a terminal condition and  the generator $g$  is continuous in $(y,z)$ and satisfies the following general growth assumption:
\begin{enumerate}
\renewcommand{\theenumi}{(EX1)}
\renewcommand{\labelenumi}{\theenumi}
\item \label{A:EX1} The generator $g$ has a general growth in $y$ and a quadratic growth in $z$, i.e., there exists a real-valued progressively measurable process $(f_t)_{t\in\T}$ with $\ps$, $\int_0^T f_t{\rm d}t<+\infty$, a nonnegative real function $H(x)$ defined on $\R_+$ with $H(0)=0$, and a positive real $c>0$ such that $\as$, for each $(y,z)\in \R\times\R^d$,
    $$
    |g(\omega,t,y,z)|\leq f_t(\omega)+H(|y|)+c|z|^2.\vspace{-0.1cm}
    $$
\end{enumerate}
For each pair of positive integers $n,p\geq 1$, set\vspace{-0.1cm}
\begin{equation}\label{eq:2.1}
\xi^{n,p}:=\xi^+\wedge n-\xi^-\wedge p\ \ \ \ {\rm and}\ \ \ \ g^{n,p}(\omega,t,y,z):=g^+(\omega,t,y,z)\wedge n-g^-(\omega,t,y,z)\wedge p,\vspace{-0.1cm}
\end{equation}
where and hereafter, $a^+$ stands for the maximum of $a$ and $0$, $a^-:=(-a)^+$ and $a\wedge b$ represents the minimum of $a$ and $b$. It follows from \cite{Kobylanski2000AP} that for each $n,p\geq 1$, the following BSDE$(\xi^{n,p},g^{n,p})$ admits a minimal (maximal) bounded solution $(Y^{n,p}_t,Z^{n,p}_t)_{t\in\T}$ such that $Y^{n,p}_\cdot$ is bounded and $Z^{n,p}_\cdot\in \mcal^2$:
\begin{equation}\label{eq:2.2}
  Y^{n,p}_t=\xi^{n,p}+\int_t^T g^{n,p}(s,Y^{n,p}_s,Z^{n,p}_s){\rm d}s-\int_t^T Z^{n,p}_s \cdot {\rm d}B_s, \ \ t\in\T.
\end{equation}
By the comparison theorem, $Y^{n,p}_\cdot$ is nondecreasing in $n$ and non-increasing in $p$. Furthermore, we have the following monotone stability theorem, which slightly generalizes those of ~\cite{Kobylanski2000AP} and \cite{BriandHu2008PTRF}.

\begin{pro}\label{Pro:2.1--stability}
Assume that $\xi$ is a terminal condition and the generator $g$ is continuous in the state variables $(y,z)$ and satisfies \ref{A:EX1}. For each pair of positive integers $n,p\geq 1$, let $\xi^{n,p}$ and $g^{n,p}$ be defined in \eqref{eq:2.1}, and $(Y^{n,p}_\cdot,Z^{n,p}_\cdot)$ be the minimal (maximal) bounded solution of \eqref{eq:2.2}. If there exists an $\R_+$-valued, progressively measurable and continuous process $(X_t)_{t\in\T}$ such that
$$\as,\ \ \RE\ n,p\geq 1,\ \ \ |Y^{n,p}_\cdot|\leq X_\cdot,$$
then there exists an $\R^d$-valued  progressively measurable process $(Z_t)_{t\in\T} $ such that $(Y_\cdot:=\inf_p\sup_n Y^{n,p}_\cdot, \ Z_\cdot)$ is a solution to BSDE $(\xi,g)$.
\end{pro}

\begin{proof}
We use the localization procedure to construct the desired solution. For each integer $m\geq 1$, define the following stopping time:\vspace{-0.1cm}
\[
 \sigma_m:=\inf\bigg\{t\in \T: X_t+\int_0^tf_s {\rm d}s\geq m\bigg\}\wedge T\vspace{0.1cm}
\]
with the convention that $\inf\emptyset=+\infty$. Then $(Y^{n,p}_m(t), Z^{n,p}_m(t)):=(Y^{n,p}_{t\wedge\sigma_m}, Z^{n,p}_t{\bf 1}_{t\leq\sigma_m})$ solves the BSDE:\vspace{0.1cm}
\[
Y^{n,p}_m(t)=Y^{n,p}_{\sigma_m}+\int_t^T{{\bf 1}_{s\leq \sigma_m}g^{n,p}(s, Y^{n,p}_m(s), Z^{n,p}_m(s)){\rm d}s}-\int_t^T Z^{n,p}_m(s) \cdot {\rm d}B_s,\ \ t\in\T.\vspace{0.1cm}
\]
Note that for fixed $m\geq 1$, $Y^{n,p}_m(\cdot)$ is nondecreasing in $n$, non-increasing in $p$ and bounded by $m$, and that $\as$, the sequence $(g^{n,p})_{n,p}$ converges locally uniformly in $(y,z)$ to the generator $g$ as $n,p\To \infty$. Furthermore, in view of the facts that $|g^{n,p}|\leq |g|$ and $g$ satisfies \ref{A:EX1}, it follows that $\ass$,
\[
\RE\ (y, z)\in [-m, m]\times\R^d,\ \ \ \sup\limits_{n,p\geq 1}\left({\bf 1}_{s\leq \sigma_m}\left|g^{n,p}(s, y, z)\right|\right)\leq {\bf 1}_{s\leq \sigma_m}f_s+ H(m)+c|z|^2,\ \ \vspace{-0.1cm}
\]
with $\ps$, $\int_0^T {\bf 1}_{s\leq \sigma_m} f_s{\rm d}s\leq m$. Thus, we can apply the stability property for bounded solutions of BSDEs (see e.g. Proposition 3.1 in \cite{LuoFan2018SD}). Setting $Y_m(t):=\inf_p\sup_n Y^{n,p}_{t\wedge\sigma_m}$, then $Y_m(\cdot)$ is continuous and there exists a process $Z_m(\cdot)$ such that $$\lim_{n,p\To\infty}Z^{n,p}_{t}{\bf 1}_{t\le \sigma_m}=Z_m(t)\ \ {\rm in}\ \mcal^2$$
and\vspace{-0.1cm}
$$
Y_m(t)=\inf_p\sup_n Y^{n,p}_{\sigma_m}+\int_t^T{{\bf 1}_{s\le \sigma_m}g(s,Y_m(s),Z_m(s)){\rm d}s}-\int_t^T{Z_m(s)\cdot {\rm d}B_s},\ \ \ t\in \T.\vspace{0.1cm}
$$
Finally, in view of the assumptions of \cref{Pro:2.1--stability}, the stability of stopping time $\sigma_m$ and the fact that $\as$, for each $m\geq 1$,
\[
Y_{m+1}(t\wedge \sigma_m)=Y_m(t\wedge \sigma_m)=\inf_{p\geq1}\sup_{n\geq 1} Y^{n,p}_{t\wedge\sigma_m}
\]
and\vspace{-0.1cm}
\[
Z_{m+1}{\bf 1}_{t\leq \sigma_m}=Z_m{\bf 1}_{t\leq \sigma_m}=\lim_{n,p\To\infty}Z^{n,p}_{t}{\bf 1}_{t\le \sigma_m},\vspace{0.1cm}
\]
the conclusion of \cref{Pro:2.1--stability} follows immediately by picking
$$
Z_t:=\sum_{m=1}^{+\infty}Z_m(t){\bf 1}_{t\in [\sigma_{m-1},\sigma_m)},\ \ \ \  t\in\T
$$
with $\sigma_0:=0$. The proof is complete.
\end{proof}

Now, to apply the preceding stability theorem, we need the uniform a priori bound $X_\cdot$ for the processes $Y^{n,p}_\cdot$, which is crucial. Our strategy is to first search for an appropriate function $\phi(s, x)$ and then apply It\^{o}-Tanaka's formula to $\phi(s,|Y^{n,p}_s|)$ on the time interval $s\in [t,\tau_m]$ with $t\in (0,T]$ and $\tau_m$ being an $(\F_t)$-stopping time valued in $[t,T]$. More specifically, we need to find a constant $\bar\delta\geq 0$ and a nonnegative  smooth function $\phi(\cdot, \cdot): (0,T]\times \R_+\to \R_+$ such that $\phi_x(s,x)>0$, $\phi_{xx}(s,x)>\bar\delta$ and
\begin{equation}\label{eq:2.3}
-\phi_x(s,x)\left(\beta x(\ln x)^{\delta}{\bf 1}_{x>1}+\gamma |z||\ln |z||^{\lambda}\right)+{1\over 2}\left(\phi_{xx}(s,x)-\bar\delta\right)|z|^2+\phi_s(s,x)\geq 0
\end{equation}
for all $(s,x,z)\in (0,T]\times \R_+\times\R^d$.
Here and hereafter, $\phi_s$ is the first-order partial derivative of $\phi$ with respect to the first variable, and $\phi_x$ and $\phi_{xx}$ are the first- and second-order partial derivatives of $\phi$ with respect to the second variable.

To find the function $\phi$ satisfying inequality \eqref{eq:2.3}, our next idea is to find, for the term in the left hand side of \eqref{eq:2.3}, a perfect lower bound $\Phi(s,x)$ without containing the variable $z$, and transfer \eqref{eq:2.3} to the inequality $\Phi(s,x)\geq 0$. For the case $\lambda=0$, we use the inequality $2ab\leq a^2+b^2$ to get that for each $(s,x,z)\in (0,T]\times \R_+\times\R^d$,
$$
\begin{array}{lll}
\Dis -\gamma \phi_x(s,x)|z|+{1\over 2}\left(\phi_{xx}(s,x)-\bar\delta\right)|z|^2&=& \Dis {1\over 2}\left(\phi_{xx}(s,x)-\bar\delta\right)\left( -2{\gamma\phi_x(s,x)\over \phi_{xx}(s,x)-\bar\delta}|z|+|z|^2\right)\vspace{0.2cm}\\
&\geq & \Dis -{\gamma^2\over 2}{(\phi_x(s,x))^2\over \phi_{xx}(s,x)-\bar\delta}.
\end{array}
$$
Hence, \eqref{eq:2.3} holds if the function $\phi(\cdot,\cdot)$ satisfies the following condition:\vspace{0.1cm}
\begin{equation}\label{eq:2.4}
\RE\ (s,x)\in (0,T]\times \R_+, \ \ \ -\beta\phi_x(s,x) x(\ln x)^{\delta}{\bf 1}_{x>1}-{\gamma^2\over 2}{(\phi_x(s,x))^2\over \phi_{xx}(s,x)-\bar\delta}+\phi_s(s,x)\geq 0.\vspace{0.1cm}
\end{equation}
However, for the other case of $\lambda\in (0,+\infty)$, the situation seems to be more complicated.

We first establish the following inequality.

\begin{pro}\label{Pro:2.2Inequality}
Given $\lambda>0$. For each $k>1$, there is a positive constant $C_{k,\lambda}>0$ depending only on $(k,\lambda)$ such that\vspace{-0.1cm}
\begin{equation}\label{eq:2.5}
\RE\ x,y>0,\ \  \ 2xy|\ln y|^\lambda\leq x^2\left(k 4^{(\lambda-1)^+}
\left|\ln x\right|^{2\lambda}+C_{k,\lambda}\right)+y^2.
\end{equation}

Moreover, for each $k\leq 1/4^{(\lambda-1)^+},$  there is no constant $C_{k,\lambda}$ to satisfy the last inequality.
\end{pro}

\begin{proof} Note first that \eqref{eq:2.5} can be rewritten as
$$
2\frac{y}{x}\left|\ln \frac{y}{x}+ \ln x\right|^\lambda\leq k4^{(\lambda-1)^+}\left|\ln x\right|^{2\lambda}+C_{k,\lambda}
+\left(\frac{y}{x}\right)^2.
$$
Set $a:=y/x>0$. Then, to prove \eqref{eq:2.5}, it is equivalent to prove that
\begin{equation}\label{eq:2.6}
\RE\ a,x>0,\ \ \ 2 a|\ln a+\ln x|^\lambda\leq k4^{(\lambda-1)^+}\left|\ln x\right|^{2\lambda}+a^2+C_{k,\lambda}.
\end{equation}

Now, we let $k>1$ and prove \eqref{eq:2.6}. By virtue of Young's inequality, observe that for each $a,x>0$,
\begin{equation}\label{eq:2.7}
|\ln a+\ln x|^\lambda\leq 2^{(\lambda-1)^+}|\ln a|^\lambda+2^{(\lambda-1)^+}|\ln x|^\lambda
\end{equation}
and
\begin{equation}\label{eq:2.8}
2 a\left(2^{(\lambda-1)^+}|\ln x|^\lambda\right)\leq \frac{1}{k}a^2+k4^{(\lambda-1)^+}\left|\ln x\right|^{2\lambda}.
\end{equation}
And, note that the function
$$
f(a;k,\lambda):=2 a \left(2^{(\lambda-1)^+}|\ln a|^\lambda\right)-\left(1-\frac{1}{k}\right)a^2\vspace{0.1cm}
$$
is continuous on $a\in (0,+\infty)$, and respectively tends to $0$ and $-\infty$ as $a$ tends to $0^+$ and $+\infty$, and then there must exist a constant $C_{k,\lambda}>0$ depending only on $(k,\lambda)$ such that
\begin{equation}\label{eq:2.9}
\RE\ a>0,\ \ \ f(a;k,\lambda)\leq C_{k,\lambda}.
\end{equation}
Combining \eqref{eq:2.7}, \eqref{eq:2.8} and \eqref{eq:2.9} yields that
$$
\begin{array}{ll}
&\Dis 2 a|\ln a+\ln x|^\lambda- k4^{(\lambda-1)^+}\left|\ln x\right|^{2\lambda}-a^2\vspace{0.2cm}\\
\leq & \Dis 2 a\left(2^{(\lambda-1)^+}|\ln a|^\lambda\right)-\left(1-\frac{1}{k}\right)a^2+2 a\left(2^{(\lambda-1)^+}|\ln x|^\lambda\right)- k4^{(\lambda-1)^+}\left|\ln x\right|^{2\lambda}-\frac{1}{k}a^2
\vspace{0.2cm}\\
\leq & \Dis f(a;k,\lambda)\leq C_{k,\lambda}, \ \ \ \RE\ a,x>0.
\end{array}
$$
This is just \eqref{eq:2.6}. Thus, we have proved  \eqref{eq:2.5}  for each $k>1$.

Next, we show the second part of this proposition. In fact, it suffices to prove that \eqref{eq:2.6} does not hold for $k=1/4^{(\lambda-1)^+}$, i.e., there does not exist a positive constant $C_{\lambda}>0$ such that
\begin{equation}\label{eq:2.10}
\RE\ a,x>0,\ \ \ 2 a|\ln a+\ln x|^\lambda\leq \left|\ln x\right|^{2\lambda}+a^2+C_{\lambda}.
\end{equation}
Let $\bar a(x):=|\ln x|^\lambda$. Then $2 \bar a(x) |\ln x|^\lambda
=\left|\ln x\right|^{2\lambda}+\bar a^2(x)$. Thus, if \eqref{eq:2.10} holds, then we have
$$
g(x;\lambda):=2 \bar a(x) |\ln \bar a(x)+\ln x|^\lambda-2 \bar a(x)|\ln x|^\lambda\leq C_{\lambda},\ \ \  x>0,
$$
which is impossible since $g(x;\lambda)$ tends to positive infinity as $x\To +\infty$. The proof is then complete.
\end{proof}

With \cref{Pro:2.2Inequality} in the hand, coming back to \eqref{eq:2.3}, observe that for each $k>1$, $\lambda>0$, $\gamma>0$ and $\bar\delta\geq 0$, there exists a constant $C_{k,\lambda}>0$ such that for each $(s,x,z)\in \T\times \R_+\times\R^d$,
$$
\begin{array}{ll}
&\Dis -\gamma\phi_x(s, x)|z||\ln|z||^{\lambda}+{1\over 2}\left(\phi_{xx}(s, x)-\bar\delta\right)|z|^2\vspace{0.3cm}\\
= & \Dis {1\over 2}\left(\phi_{xx}(s, x)-\bar\delta\right)\left(-2\frac{\gamma\phi_x(s, x)}
{\phi_{xx}(s, x)-\bar\delta}|z||\ln|z||^{\lambda}
+|z|^2\right)\vspace{0.3cm}\\
\geq  & \Dis -{\gamma^2\over 2}
\frac{\left(\phi_x(s, x)\right)^2}
{\phi_{xx}(s, x)-\bar\delta}\left(k4^{(\lambda-1)^+}\left|\ln \frac{\gamma\phi_x(s, x)}{\phi_{xx}(s, x)-\bar\delta}
\right|^{2\lambda}+C_{k,\lambda}\right).
\end{array}
$$
Then, \eqref{eq:2.3} holds if for some $k>1$ and $C_{k,\lambda}>0$ associated with $k$ by \cref{Pro:2.2Inequality}, the function $\phi(s,x)$ satisfies that for each $(s,x)\in (0,T]\times \R_+$, with $\gamma>0$,\vspace{0.1cm}
\begin{equation}\label{eq:2.11}
-\beta\phi_x(s, x) x(\ln x)^{\delta} {\bf 1}_{x>1} \Dis -{\gamma^2\over 2}
\frac{\left(\phi_x(s, x)\right)^2}
{\phi_{xx}(s, x)-\bar\delta}\left(k4^{(\lambda-1)^+}\left|\ln \frac{\gamma\phi_x(s, x)}{\phi_{xx}(s, x)-\bar\delta}
\right|^{2\lambda}+C_{k,\lambda}\right)+\phi_s(s, x)\geq 0.\vspace{0.1cm}
\end{equation}
In the following four sections, for four different cases,  we shall find the constant $\bar\delta\geq 0$ and the function $\phi$ satisfying \eqref{eq:2.4} and \eqref{eq:2.11}, and then the uniform a priori bound $X_\cdot$ for $Y^{n,p}_\cdot$ in \cref{Pro:2.1--stability}.

\subsection{Uniqueness}

For the uniqueness and comparison theorem of the unbounded solutions to BSDE $(\xi,g)$ with generator $g$ satisfying \eqref{eq:1.2}, more assumptions on the generator $g$ are required  as usual than those required for the existence.

First, let us introduce the following three  assumptions.
\begin{enumerate}
\renewcommand{\theenumi}{(UN\arabic{enumi})}
\renewcommand{\labelenumi}{\theenumi}
\item \label{A:UN1} The generator $g$ satisfies the one-sided Osgood condition in $y$, i.e., there exists a nonnegative, increasing, continuous and concave function $\rho(\cdot)$ defined on $\R_+$ satisfying $\rho(0)=0$, $\rho(u)>0$ for $u>0$ and
    $$\int_{0^+}\frac{{\rm d}u}{\rho(u)}:=\lim\limits_{\eps\To 0^+}\int_0^\eps\frac{{\rm d}u}{\rho(u)}=+\infty\vspace{0.1cm}$$
    such that $\as$, for each $(y_1,y_2,z)\in \R\times\R\times\R^d$,
    $${\rm sgn}(y_1-y_2) \left(g(\omega,t,y_1,z)-g(\omega,t,y_2,z)\right)\leq \rho(|y_1-y_2|).$$

\item \label{A:UN2} The generator $g$ is uniformly continuous in $z$, i.e., there exists a nonnegative, nondecreasing and continuous function $\kappa(\cdot)$ defined on $\R_+$ with $\kappa(0)=0$ such that $\as$, for each $(y,z_1,z_2)\in \R\times\R^d\times\R^d$,
    $$|g(\omega,t,y,z_1)-g(\omega,t,y,z_2)|\leq \kappa (|z_1-z_2|).$$

\item\label{A:UN3} $\as$, the generator $g$ is convex or concave in $(y,z)$.
\end{enumerate}

\begin{rmk}\label{rmk:2.3} Assumptions \ref{A:UN1} and \ref{A:UN2} are respectively strictly weaker than the usual monotonicity condition of the generator $g$ in $y$ and the uniformly Lipschitz continuity condition of $g$ in $z$. They are usually used when the generator $g$ has a general growth in $y$ and a linear growth in $z$. See \cite{FanJiang2013AMSE,Fan2015JMAA,Fan2016SPA} among others for more details. And, to study the uniqueness, Assumption \ref{A:UN3} seems to be very natural for a non-linear growth function, see for example \cite{BriandHu2008PTRF} and \cite{DelbaenHuRichou2011AIHPPS}, where the convexity (concavity) condition of the generator $g$ in $z$ are required to ensure the uniqueness of the solution.
\end{rmk}

\begin{rmk}\label{rmk:2.4}
Without loss of generality, we always assume that the functions $\rho(\cdot)$ and $\kappa(\cdot)$ defined respectively in \ref{A:UN1} and \ref{A:UN2} are of linear growth, i.e., there exists a constant $A>0$ such that
$$
\RE\ u\in\R_+, \ \ \ \rho(u)\leq A(u+1)\ \ {\rm and}\ \ \kappa(u)\leq A(u+1).
$$
\end{rmk}

Now, we establish a general comparison theorem for solutions of BSDEs under Assumptions \ref{A:UN1} and \ref{A:UN2}, which naturally yields the uniqueness of the solution. In particular, it is clear that this comparison theorem strengthens the uniqueness part of~\citet[Theorem 3.1]{FanHu2019ECP}.

\begin{pro}\label{pro:2.5}
Let $\xi$ and $\xi'$ be two terminal conditions, $g$ and $g'$ be two generators which are both continuous in the state variables $(y,z)$, and $(Y_t, Z_t)_{t\in\T}$ and $(Y'_t, Z'_t)_{t\in\T}$ be respectively a solution to BSDE $(\xi, g)$ and BSDE $(\xi', g')$ such that both $\psi(|Y_\cdot|,\mu_\cdot)$ and $\psi(|Y'_\cdot|,\mu_\cdot)$ belongs to class (D), where\vspace{0.1cm}
\begin{equation}\label{eq:2.12}
\psi(x,\mu):=x\exp\left(\mu\sqrt{2\ln (1+x)}\right),\ \ (x,\mu)\in \R_+\times \R_+.
\end{equation}
and $\mu_\cdot$ is any nonnegative, strictly increasing and continuous function defined on $\T$ with $\mu_0=0$.

Assume that $\ps$, $\xi\leq \xi'$. If $g$ (resp. $g'$) verifies Assumptions \ref{A:UN1} and \ref{A:UN2}, and $\as$,
\begin{equation}\label{eq:2.13}
{\bf 1}_{Y_t>Y'_t}\left(g(t,Y'_t,Z'_t)-g'(t,Y'_t,Z'_t)\right)\leq 0\ \ \ ({\rm resp.}\  \ {\bf 1}_{Y_t>Y'_t}\left(g(t,Y_t,Z_t)-g'(t,Y_t,Z_t)\right)\leq 0\ ),
\end{equation}
then $\ps$, for each $t\in\T$, $Y_t\leq Y'_t$.
\end{pro}

To prove the above proposition, we need the following two lemmas of~\citet{HuTang2018ECP}.\vspace{-0.1cm}

\begin{lem}\label{lem:2.6}
Let the function $\psi$ be defined in \eqref{eq:2.12}. Then we have
\begin{itemize}
\item [(i)] For each $x\in\R_+$, $\psi(x,\cdot)$ is increasing on $\R_+$.
\item [(ii)] For $\mu\in\R_+$, $\psi(\cdot,\mu)$ is a positive, increasing and strictly convex function on $\R_+$.
\item [(iii)] For all $c>1$ and $x,\mu\in\R_+$, $\psi(cx,\mu)\leq \psi(c,\mu)\psi(x,\mu)$.
\item [(iv)] For all $x_1,x_2,\mu\in\R_+$, $\psi(x_1+x_2,\mu)\leq {1\over 2}\psi(2,\mu)\left[\psi(x_1,\mu)+\psi(x_2,\mu)\right].$
\item [(v)] For each $x\in \R$, $y\in\R_+$ and $\mu>0$,
$e^xy\leq e^{x^2\over 2\mu^2}+e^{2\mu^2}\psi(y,\mu)$.
\end{itemize}
\end{lem}

\begin{lem}\label{lem:2.7}
Let $(q_t)_{t\in\T}$ be an $\R^d$-valued progressively measurable process with $|q_\cdot|\leq \epsilon$ almost surely. For each $t\in\T$, if $0\leq \lambda<{1 \over 2\epsilon^2 (T-t)}$, then
$$
\E\left[\left.e^{\lambda \left|\int_t^Tq_s\cdot {\rm d}B_s\right|^2}\right|\F_t\right]\leq {1\over \sqrt{1-2\lambda\epsilon^2(T-t)}}.\vspace{0.2cm}
$$
\end{lem}

Now, we prove~\cref{pro:2.5}.

\begin{proof}[Proof of \cref{pro:2.5}]
We only prove the case that $g$ verifies Assumptions \ref{A:UN1} and \ref{A:UN2}, and $\as$, ${\bf 1}_{Y_t>Y'_t}\left(g(t,Y'_t,Z'_t)-g'(t,Y'_t,Z'_t)\right)\leq 0$. In the same way, another case can be proved.

Define $\hat Y_\cdot:=Y_\cdot-Y'_\cdot$ and $\hat Z_\cdot:=Z_\cdot-Z'_\cdot$.
Then the pair $(\hat Y_\cdot,\hat Z_\cdot)$ verifies \vspace{0.1cm}
\begin{equation}\label{eq:2.14}
  \hat Y_t=\xi-\xi'+\int_t^T \left(g(s,Y_s,Z_s)-g'(s, Y'_s, Z'_s)\right){\rm d}s-\int_t^T \hat Z_s \cdot {\rm d}B_s, \ \ \ \ t\in\T.\vspace{0.1cm}
\end{equation}
From the assumption of ${\bf 1}_{Y_t>Y'_t}\left(g(t,Y'_t,Z'_t)-g'(t,Y'_t,Z'_t)\right)\leq 0$ and Assumptions \ref{A:UN1} and \ref{A:UN2} of the generator $g$ together with \cref{rmk:2.4}, it is not difficult to deduce  that $\ass$,
\begin{equation}\label{eq:2.15}
{\bf 1}_{\hat Y_s>0}\left(g(s,Y_s,Z_s)-g'(s, Y'_s, Z'_s)\right)\leq \rho(\hat Y_s^+)+{\bf 1}_{\hat Y_s>0}\kappa(|\hat Z_s|)\leq A \hat Y_s^+ +A{\bf 1}_{\hat Y_s>0}|\hat Z_s|+2A.
\end{equation}

For each $t\in [0,T]$ and each integer $n\geq 1$, define the following stopping times:
$$
\sigma_n:=\inf\left\{s\in [t,T]:\ |\hat Y_s|+\int_t^s |\hat Z_r|^2 {\rm d}r\geq n\right\}\wedge T.
$$
It\^{o}-Tanaka's formula applied \eqref{eq:2.14} gives\vspace{0.1cm}
$$
  \hat Y_t^+= \hat Y_{\sigma_n}^+ +\int_t^{\sigma_n} {\bf 1}_{\hat Y_s>0}\left(g(s,Y_s,Z_s)-g'(s, Y'_s, Z'_s)\right){\rm d}s-\int_t^{\sigma_n} {\bf 1}_{\hat Y_s>0}\hat Z_s \cdot {\rm d}B_s-\int_t^{\sigma_n}{\rm d}L_s, \ \ \ \ t\in\T,\vspace{0.1cm}
$$
where $L_\cdot$ is the local time of $\hat Y_\cdot$ at the origin. Then, in view of \eqref{eq:2.15},\vspace{0.1cm}
\begin{equation}\label{eq:2.16}
  e^{At}\hat Y_t^+\leq e^{AT}(\hat Y_{\sigma_n}^+ +2AT)+\int_t^{\sigma_n} v_s\cdot \tilde Z_s {\rm d}s-\int_t^{\sigma_n} \tilde Z_s \cdot {\rm d}B_s, \ \ \ \ t\in\T,
\end{equation}
with $v_s:=A\frac {\hat Z_s}{|\hat Z_s|^2}{\bf 1}_{|\hat Z_s|>0}$ and $\tilde Z_s:={\bf 1}_{\hat Y_s>0}e^{As}\hat Z_s$. It is clear that $|v_\cdot|\leq A$ and it follows from \eqref{eq:2.16} that \vspace{0.1cm}
\begin{equation}\label{eq:2.17}
\RE\ n\geq 1,\ \ \ \hat Y_t^+\leq e^{A T}\left(\E\left[\left.e^{\int_t^{\sigma_n}v_s\cdot {\rm d}B_s}\hat Y_{\sigma_n}^+\right|\F_t\right]+2AT\right),\ \ \ t\in\T.\vspace{0.1cm}
\end{equation}

Furthermore, in view of Assertion (v) of \cref{lem:2.6}, we know that for each $n\geq 1$,
\begin{equation}\label{eq:2.18}
e^{\int_t^{\sigma_n}v_s\cdot {\rm d}B_s}\hat Y_{\sigma_n}^+\leq e^{{1\over 2\mu^2_t}\left(\int_t^{\sigma_n}v_s\cdot {\rm d}B_s\right)^2}+e^{2\mu^2_t}\psi\left(\hat Y_{\sigma_n}^+,\mu_t\right), \ \ t\in (0,T].
\end{equation}
And, by virtue of the assumptions of $\mu_\cdot$, we can pick a (unique) $T_1\in (0,T)$ such that $\mu_{T_1}^2=4A^2(T-T_1)$ and $\mu_t^2\geq 4A^2(T-t)$ for each $t\in [T_1,T]$. It then follows from \cref{lem:2.7} that for all $n\geq 1$,\vspace{0.1cm}
$$
\E\left[\left| e^{{1\over 2\mu^2_t}\left(\int_t^{\sigma_n}v_s\cdot {\rm d}B_s\right)^2} \right|^2\right]=\E\left[e^{{1\over \mu^2_t}\left(\int_t^{\sigma_n}v_s\cdot {\rm d}B_s\right)^2} \right]\leq {1\over \sqrt{1-{2A^2(T-t)\over \mu^2_t}}}\leq \sqrt{2},\ \  t\in [T_1,T]
$$
and then, the family of random variables $e^{{1\over 2\mu^2_t}\left(\int_t^{\sigma_n}v_s\cdot {\rm d}B_s\right)^2}$ is uniformly integrable on the interval $[T_1,T]$. On the other hand, in view of \cref{lem:2.6} and the monotonicity of $\mu_\cdot$, observe that for all $n\geq 1$,\vspace{0.1cm}
$$
\begin{array}{lll}
\Dis e^{2\mu^2_t}\psi\left(\hat Y_{\sigma_n}^+,\mu_t\right)&\leq & \Dis  e^{2\mu^2_T}\psi\left(\left|Y_{\sigma_n}\right|+\left|Y'_{\sigma_n}\right|,
\mu_{\sigma_n}\right)\vspace{0.2cm}\\
&\leq & \Dis {e^{2\mu^2_T}\psi(2,\mu_T) \over 2}\left[\psi\left(\left| Y_{\sigma_n}\right|,\mu_{\sigma_n}\right)+\psi\left(\left| Y'_{\sigma_n}\right|,\mu_{\sigma_n}\right)\right],\ \ t\in\T.
\end{array}\vspace{0.1cm}
$$
It then follows from \eqref{eq:2.18} that for $t\in [T_1,T]$, the family of random variables $e^{\int_t^{\sigma_n}v_s\cdot {\rm d}B_s}\hat Y_{\sigma_n}^+$ is uniformly integrable. Consequently, in view of the fact that $\hat Y_T^+=(\xi-\xi')^+=0$, by letting $n\To\infty$ in \eqref{eq:2.17} we get that $\hat Y_\cdot\leq 2ATe^{AT}$ on the time interval $[T_1,T]$, which means that $\hat Y_\cdot$ is a bounded process on $[T_1,T]$. Thus, by virtue of the assumptions of \cref{pro:2.5} again, we can apply Theorem 2.1 in \cite{Fan2016SPA} to obtain $\hat Y_\cdot^+=0$ on the time interval $[T_1,T]$.

Now, in view of the fact $Y_{T_1}^+=0$ obtained above, we can respectively replace $T$ and $T_1$ in the above analysis with $T_1$ and $T_2\in (0,T_1)$ satisfying $\mu_{T_2}^2=4A^2(T_1-T_2)$ and $\mu_t^2\geq 4A^2(T_1-t)$ for each $t\in [T_2,T_1]$, and use the same argument to get that $\hat Y_\cdot^+=0$ on the time interval $[T_2,T_1]$.

Finally, repeating the above procedure, we successively obtain that $\hat Y_\cdot^+=0$ on the time intervals $[T_3, T_2]$, $\cdots$, $[T_{m+1}, T_m]$, $\cdots$, where $0<T_{m+1}<T_m<T$ for each $m\geq 1$, and the limit of sequence $(T_m)_{m\geq 1}$ must exist and is denoted by $\bar T$. Noticing that $\mu_{T_m}^2=4A^2(T_{m-1}-T_m)$ by the above construction and in view of the assumptions of $\mu_\cdot$, we can conclude that $\bar T=\lim\limits_{m\To\infty}T_m=0$. Thus,
in view of the continuity of process $\hat Y_t^+$ with respect to the time variable $t$, we have $\hat Y_\cdot^+=0$ on the whole time interval $\T$ by sending $m\to \infty$, which is the desired result. The proof is complete.
\end{proof}

At the end of this subsection, we would like to mention that in subsequent four sections, we shall also establish, under the additional convexity Assumption \ref{A:UN3}, the uniqueness theorem and the comparison theorem for unbounded solutions of BSDEs with super-linearly growing generators and four different integrability terminal conditions. To take advantage of the convexity (concavity) condition of the generator $g$, we use the $\theta$-technique developed in \cite{BriandHu2008PTRF} to establish the comparison theorem for the unbounded solutions, which yields directly the uniqueness of the solution. More precisely, instead of estimating the difference between the processes $Y_\cdot$ and $Y'_\cdot$, we estimate $Y_\cdot-\theta Y'_\cdot$ for each $\theta\in (0,1)$. We found that the uniform a priori estimate is still the key point.

\section{The case of $L\exp\left(\mu \sqrt{2\ln (1+L)}\right)$-integrable terminal conditions}
\label{sec:3-mainresult--3}
\setcounter{equation}{0}

In this section, we always assume that the generator $g$ satisfies the following assumption, which is strictly weaker than \eqref{eq:1.2} with $\delta=1/2$ and $\lambda=0$.
\begin{enumerate}
\renewcommand{\theenumi}{(A1)}
\renewcommand{\labelenumi}{\theenumi}
\item \label{A:A1} There exist two nonnegative constants $\beta\geq 0$ and $\gamma\geq 0$ with $\beta+\gamma>0$, and an $\R_+$-valued progressively measurable process $(\alpha_t)_{t\in\T}$  such that $\as$,
    $$
    \RE\ (y,z)\in \R\times\R^d,\quad
    {\rm sgn}(y)g(\omega,t,y,z)\leq \alpha_t(\omega)+\beta|y|\sqrt{\ln|y|}{\bf 1}_{|y|>1}+\gamma |z|.
    $$
  \end{enumerate}

Based on the analysis in Section 2, we need to find a smooth function $\phi(\cdot,\cdot):[0,T]\times \R_+\To \R_+$ such that $\phi_x(s,x)>0$, $\phi_{xx}(s,x)>0$ and Equation~\eqref{eq:2.4} is satisfied for $\delta=1/2$ and $\bar\delta=0$, i.e., \vspace{0.1cm}
\begin{equation}\label{eq:3.1}
\RE\ (s,x)\in (0,T]\times \R_+, \quad -\beta\phi_x(s,x) x\sqrt{\ln x}{\bf 1}_{x>1}-{\gamma^2\over 2}{\phi^2_x(s,x)\over \phi_{xx}(s,x)}+\phi_s(s,x)\geq 0.\vspace{0.1cm}
\end{equation}

We assume that $\beta,\gamma\geq 0$ with $\beta+\gamma>0$, and choose the following function
$$
\phi(s,x):=(x+e)\exp\left(\mu_s\sqrt{2\ln (x+e)}+\nu_s\right)>0, \ \ \ (s,x)\in [0,T]\times \R_+
$$
to explicitly solve the inequality \eqref{eq:3.1}, where $\mu_s,\nu_s:\T\To\R_+$ are two strictly increasing and continuous functions with zero only at the origin, and are continuously differential on the time interval $(0,T]$. For each $(s,x)\in (0,T]\times \R_+$, a simple computation gives
\begin{equation}\label{eq:3.2}
\phi_x(s,x)=\phi(s,x){\mu_s+\sqrt{2\ln (x+e)} \over (x+e)\sqrt{2\ln (x+e)}}>0,
\end{equation}
\begin{equation}\label{eq:3.3}
\phi_{xx}(s,x)=\phi(s,x){\mu_s\left(
2\ln (x+e)+\mu_s \sqrt{2\ln (x+e)}-1\right)\over (x+e)^2\left(\sqrt{2\ln (x+e)}\right)^3}>0
\end{equation}
and\vspace{-0.1cm}
\begin{equation}\label{eq:3.4}
\phi_s(s,x)=\phi(s,x)\left(\mu'_s\sqrt{2\ln (x+e)}+\nu'_s\right)>0.\vspace{0.1cm}
\end{equation}
Note that\vspace{-0.1cm}
\begin{equation}\label{eq:3.5}
\RE\ x\in\R_+,\ \ x\sqrt{\ln x}{\bf 1}_{x>1}\leq \frac{1}{\sqrt{2}}(x+e)\sqrt{2\ln (x+e)}.
\end{equation}
Denote $v=\sqrt{2\ln (x+e)}$. Substituting \eqref{eq:3.2}, \eqref{eq:3.3},  \eqref{eq:3.4} and \eqref{eq:3.5} into the left side of \eqref{eq:3.1} yields that
$$
\begin{array}{ll}
&\Dis -\beta\phi_x(s,x) x\sqrt{\ln x}{\bf 1}_{x>1}-{\gamma^2\over 2}{\phi^2_x(s,x)\over \phi_{xx}(s,x)}+\phi_s(s,x)\vspace{0.2cm}\\
\geq & \Dis -\beta\phi(s,x)\frac{\mu_s+v}{\sqrt{2}}-\frac{\gamma^2}{2}
\phi(s,x)\frac{v(\mu_s+v)^2}{\mu_s(v^2+\mu_sv-1)}
 +\phi(s,x)(\mu'_s v+\nu'_s), \ \ (s,x)\in (0,T]\times\R_+.\vspace{0.1cm}
\end{array}
$$
Furthermore, in view of the fact of $v\geq \sqrt{2}$, we know that
$$
\begin{array}{lll}
\Dis {v\left(\mu_s+v\right)^2\over v^2+\mu_s v-1}&=&\Dis {(v+\mu_s) (v^2+\mu_s v-1)+(v+\mu_s) \over  v^2+\mu_s v-1}\vspace{0.2cm}\\
&\leq& \Dis v+\mu_s +{v+\mu_s \over {1\over 2}(v^2+\mu_s v)}= v+\mu_s +{2\over v}\leq v+\mu_s +\sqrt{2},
\end{array}
$$
and then
$$
\begin{array}{ll}
&\Dis -\beta\phi_x(s,x) x\sqrt{\ln x}{\bf 1}_{x>1}-{\gamma^2\over 2}{\phi^2_x(s,x)\over \phi_{xx}(s,x)}+\phi_s(s,x)\vspace{0.2cm}\\
\geq & \Dis  \phi(s,x)\left\{\left(\mu'_s-\frac{\gamma^2}{2\mu_s}-\frac{\sqrt{2}}{2} \beta \right)v+\left[\nu'_s-\frac{\sqrt{2}}{2}\beta\mu_s-\frac{\gamma^2}{2}
\left(1+\frac{\sqrt{2}}{\mu_s}\right)
\right]\right\}, \ \ \ \  (s,x)\in (0,T]\times\R_+.\vspace{0.1cm}
\end{array}
$$
Thus, \eqref{eq:3.1} holds if the functions $\mu_s, \nu_s\in \T$ satisfies
\begin{equation}\label{eq:3.6}
\mu'_s=\frac{\gamma^2}{2\mu_s}+\frac{\sqrt{2}}{2} \beta,\quad s\in (0,T]; \quad \mu_0=0
\end{equation}
and
\begin{equation}\label{eq:3.7}
\nu'_s=\frac{\sqrt{2}}{2}\beta\mu_s+\frac{\gamma^2}{2}
\left(1+\frac{\sqrt{2}}{\mu_s}\right),\ \ \ s\in (0,T].\vspace{0.1cm}
\end{equation}

It is not very hard to verify that for each $\beta,\gamma\geq 0$ with $\beta+\gamma>0$, ODE~\eqref{eq:3.6} has a unique strictly increasing and continuous solution, denoted  by $\mu_{\beta,\gamma}(\cdot)$;   i.e.,
\begin{equation}\label{eq:3.8}
\mu_{\beta,\gamma}(0)=0\ \ \ {\rm and}\ \ \ \mu'_{\beta,\gamma}(s)=\frac{\gamma^2}{2\mu_{\beta,\gamma}(s)}+\frac{\sqrt{2}}{2} \beta,\ \ \ s\in (0,T].
\end{equation}
Furthermore, by integrating on both sides of the last equation, we have
$$
\mu_{\beta,\gamma}(T)-\mu_{\beta,\gamma}(t)=\frac{\gamma^2}{2}\int_t^T \frac{1}{\mu_{\beta,\gamma}(s)}{\rm d}s+\frac{\sqrt{2}}{2} \beta (T-t), \quad t\in (0,T].
$$
Letting $t\To 0^+$ in the last equation yields that the integral of $1/\mu_{\beta,\gamma}(s)$ on $\T$ is well defined, and then \cref{eq:3.7} admits a unique solution valued zero at the origin, denoted by $\nu_{\beta,\gamma}(\cdot)$, that is
\begin{equation}\label{eq:3.9}
\nu_{\beta,\gamma}(s)=\int_0^s \left[\frac{\sqrt{2}}{2}\beta\mu_{\beta,\gamma}(r)+\frac{\gamma^2}{2}
\left(1+\frac{\sqrt{2}}{\mu_{\beta,\gamma}(r)}\right)\right]{\rm d}r,\ \ \ s\in [0,T].\vspace{0.2cm}
\end{equation}

\begin{rmk}\label{rmk:3.1}
In general, $\mu_{\beta,\gamma}(\cdot)$ does not have an explicit expression. However, it is easy to check that $\mu_{\beta,\gamma}(s)=\gamma \sqrt{s}$ in the case of $\beta=0$, and $\mu_{\beta,\gamma}(s)=\frac{\sqrt{2}}{2}\beta s$ in the case of $\gamma=0$.\vspace{0.2cm}
\end{rmk}

In summary, we have

\begin{pro}\label{pro:3.2}
For each $\beta,\gamma\geq 0$ with $\beta+\gamma>0$, denote the function
\begin{equation}\label{eq:3.10}
\varphi(s,x):=(x+e)\exp\left(\mu_{\beta,\gamma}(s)\sqrt{2\ln (x+e)}
+\nu_{\beta,\gamma}(s)\right), \ \ \ (s,x)\in [0,T]\times \R_+
\end{equation}
with $\mu_{\beta,\gamma}(\cdot)$ and $\nu_{\beta,\gamma}(\cdot)$ being respectively defined in \eqref{eq:3.8} and \eqref{eq:3.9}. Then we have
\begin{itemize}
\item [(i)] $\varphi(\cdot,\cdot)$ is continuous on $[0,T]\times \R_+$; And, $\varphi(\cdot,\cdot)\in \mathcal {C}^{1,2}\left((0,T]\times \R_+\right)$;
\item [(ii)] $\varphi(\cdot,\cdot)$ satisfies the inequality \eqref{eq:2.3} with $\delta=1/2$, $\bar\delta=0$ and $\lambda=0$, i.e., \vspace{0.1cm}
$$
-\varphi_x(s,x)\left(\beta x\sqrt{\ln x}{\bf 1}_{x>1}+\gamma |z|\right)+{1\over 2}\varphi_{xx}(s,x)|z|^2+\varphi_s(s,x)\geq 0,\ \ (s,x,z)\in (0,T]\times \R_+\times\R^d.
$$
\end{itemize}
\end{pro}

The following existence and uniqueness theorem is the main result of this section.

\begin{thm}\label{thm:3.3}
Let the functions $\psi(\cdot,\cdot)$ and $\mu_{\beta,\gamma}(\cdot)$ be respectively defined in \eqref{eq:2.12} and \eqref{eq:3.8}, $\xi$ be a terminal condition and $g$ be a generator which is continuous in $(y,z)$. If $g$ satisfies assumptions \ref{A:EX1} and \ref{A:A1} with parameters $\alpha_\cdot$, $\beta$ and $\gamma$, and the terminal condition satisfies
\begin{equation}\label{eq:3.11}
\E\left[\psi\left(|\xi|+\int_0^T \alpha_t {\rm d}t,\ \mu_{\beta,\gamma}(T)\right)\right]<+\infty,
\end{equation}
then BSDE$(\xi,g)$ admits a solution $(Y_t,Z_t)_{t\in\T}$ such that $\left(\psi\left(|Y_t|,\ \mu_{\beta,\gamma}(t)\right)\right)_{t\in\T}$ belongs to class (D), and $\ps$, for each $t\in \T$,
\begin{equation}\label{eq:3.12}
|Y_t|\leq \psi(|Y_t|,\ \mu_{\beta,\gamma}(t))\leq C\E\left[\left.\psi\left(|\xi|+\int_0^T\alpha_t {\rm d}t,\ \mu_{\beta,\gamma}(T)\right)\right|\F_t\right]+C,
\end{equation}
where $C$ is a positive constant depending only on $(\beta,\gamma,T)$.

Moreover, if $g$ satisfies either assumptions \ref{A:UN1} and \ref{A:UN2} or assumption \ref{A:UN3}, then BSDE $(\xi,g)$ admits a unique solution $(Y_t,Z_t)_{t\in\T}$ such that $\left(\psi\left(|Y_t|,\ \mu_{\beta,\gamma}(t)\right)\right)_{t\in\T}$ belongs to class (D).
\end{thm}

In order to prove this theorem, we need the following two propositions.

\begin{pro}\label{pro:3.4}
Let the functions $\psi(x,\mu)$, $\mu_{\beta,\gamma}(s)$ and $\varphi(s,x)$ be respectively defined on \eqref{eq:2.12}, \eqref{eq:3.8} and \eqref{eq:3.10}. Then, there exists a universal constant $K>0$ depending only on $(\beta,\gamma,T)$ such that\vspace{-0.1cm}
\begin{equation}\label{eq:3.13}
\RE\ (s,x)\in [0,T]\times\R_+,\ \ \ \psi(x,\mu_{\beta,\gamma}(s))\leq \varphi(s,x)\leq K\psi(x,\mu_{\beta,\gamma}(s))+K.\vspace{-0.1cm}
\end{equation}
\end{pro}

\begin{proof}
The first inequality in \eqref{eq:3.13} is obvious. We now prove the second inequality. In fact,
$$
\begin{array}{lll}
\Dis{\varphi(s,x) \over\psi(x,\mu_{\beta,\gamma}(s))+1}&=&\Dis {(x+e)\exp\left(\mu_{\beta,\gamma}(s)\sqrt{2\ln (x+e)}
+\nu_{\beta,\gamma}(s)\right)\over x\exp\left(\mu_{\beta,\gamma}(s)\sqrt{2\ln (1+x)}\right)+1}\vspace{0.2cm}\\
&\leq & \Dis {x+e\over x}\exp\left(\mu_{\beta,\gamma}(T)\left(\sqrt{2\ln (x+e)}-\sqrt{2\ln (x+1)}\right)
+\nu_{\beta,\gamma}(T)\right)\vspace{0.2cm}\\
&=: & \Dis H_1(x;\beta,\gamma,T),\ \ \ \ \ (s,x)\in [0,T]\times [1,+\infty).
\end{array}
$$
And, in the case of $x\in [0,1]$,
$$
{\varphi(s,x) \over\psi(x,\mu_{\beta,\gamma}(s))+1}\leq (1 +e)\exp\left(\mu_{\beta,\gamma}(T)\sqrt{2\ln (1+e)}+\nu_{\beta,\gamma}(T)\right)=: H_2(\beta,\gamma, T),\ \ \ s\in\T.
$$
Hence, for all $x\in \R_+$, we have
\begin{equation}\label{eq:3.14}
{\varphi(s,x) \over\psi(x,\mu_{\beta,\gamma}(s))+1}\leq H_1(x;\beta,\gamma,T){\bf 1}_{x\geq 1}+H_2(\beta, \gamma, T){\bf 1}_{0\leq x<1},\ \ \ s\in\T.\vspace{0.1cm}
\end{equation}
Thus, in view of \eqref{eq:3.14} and the fact that the function $H_1(x;\beta, \gamma, T)$ is continuous on $[1,+\infty)$ and tends to $\exp\left(\nu_{\beta,\gamma}(T)\right)$ as $x\To +\infty$, the second inequality in \eqref{eq:3.13} follows immediately. The proof is complete.
\end{proof}

\begin{pro}\label{pro:3.5}
Let the functions $\psi(\cdot,\cdot)$ and $\mu_{\beta,\gamma}(\cdot)$ be respectively defined in \eqref{eq:2.12} and \eqref{eq:3.8}, $\xi$ be a terminal condition and $g$ be a generator which is continuous in $(y,z)$. If $g$ satisfies assumption \ref{A:A1} with parameters $\alpha_\cdot$, $\beta$ and $\gamma$, $|\xi|+\int_0^T \alpha_t {\rm d}t$ is a bounded random variable, and $(Y_t,Z_t)_{t\in\T}$ is a solution of BSDE$(\xi,g)$ such that $Y_\cdot$ is bounded, then $\ps$, for each $t\in \T$, we have
\begin{equation}\label{eq:3.15}
|Y_t|\leq \psi(|Y_t|,\ \mu_{\beta,\gamma}(t))\leq C\E\left[\left.\psi\left(|\xi|+\int_0^T\alpha_t {\rm d}t,\ \mu_{\beta,\gamma}(T)\right)\right|\F_t\right]+C,
\end{equation}
where $C$ is a positive constant depending only on $(\beta,\gamma,T)$.
\end{pro}

\begin{proof}
Define
$$
\bar Y_t:=|Y_t|+\int_0^t \alpha_s {\rm d}s\ \ \ \
{\rm and}\ \ \ \ \bar Z_t:={\rm sgn}(Y_t)Z_t,\ \ \ \ t\in \T.
$$
It follows from It\^{o}-Tanaka's formula that
$$
\bar Y_t=\bar Y_T+\int_t^T \left({\rm sgn}(Y_s)g(s,Y_s,Z_s)-\alpha_s\right){\rm d}s-\int_t^T \bar Z_s \cdot {\rm d}B_s-\int_t^T {\rm d}L_s, \ \ \ t\in\T,
$$
where $L_\cdot$ denotes the local time of $Y_\cdot$ at the origin. Now, we apply It\^{o}-Tanaka's formula to the process $\varphi(s, \bar Y_s)$, where the function $\varphi(\cdot,\cdot)$ is defined in \eqref{eq:3.10}, to derive, in view of assumption \ref{A:A1},
$$
\begin{array}{lll}
\Dis {\rm d}\varphi(s,\bar Y_s)
&=&\Dis \varphi_x(s,\bar Y_s)
\left(-{\rm sgn}(Y_s)g(s,Y_s,Z_s)+\alpha_s\right){\rm d}s+\varphi_x(s,\bar Y_s)\bar Z_s \cdot {\rm d}B_s+\varphi_x(s,\bar Y_s){\rm d}L_s\vspace{0.1cm}\\
&&\Dis +{1\over 2}\varphi_{xx}(s,\bar Y_s)|Z_s|^2{\rm d}s+\varphi_s(s,\bar Y_s){\rm d}s
\vspace{0.2cm}\\
&\geq &\Dis \left[-\varphi_x(s,\bar Y_s)\left(\beta |Y_s|\sqrt{\ln |Y_s|}{\bf 1}_{|Y_s|>1}+\gamma |Z_s|\right)+{1\over 2}\varphi_{xx}(s,\bar Y_s)|Z_s|^2+\varphi_s(s,\bar Y_s)\right]{\rm d}s\vspace{0.1cm}\\
&&\Dis +\varphi_x(s,\bar Y_s)\bar Z_s \cdot {\rm d}B_s,\ \ s\in \T.
\end{array}
$$
Furthermore, in view of the fact that $|Y_s|\sqrt{\ln |Y_s|}{\bf 1}_{|Y_s|>1}\leq \bar Y_s\sqrt{\ln \bar Y_s}{\bf 1}_{\bar Y_s>1}$, thanks to \cref{pro:3.2} we obtain that\vspace{-0.1cm}
\begin{equation}\label{eq:3.16}
{\rm d}\varphi(s,\bar Y_s)\geq \varphi_x(s,\bar Y_s)\bar Z_s \cdot {\rm d}B_s,\ \ s\in (0,T].
\end{equation}
Let us consider, for each integer $n\geq 1$ and each $t\in (0,T]$,  the following stopping time
$$
\tau_n^t:=\inf\left\{s\in [t,T]: \int_t^s \left[\varphi_x(r,\bar Y_r)\right]^2|\bar Z_r|^2{\rm d}r\geq n \right\}\wedge T.
$$
It follows from the definition of $\tau_n^t$ and the inequality \eqref{eq:3.16} that for each $n\geq 1$,
$$
\varphi(t,\bar Y_t)\leq \E\left[\left. \varphi(\tau_n^t,\bar Y_{\tau_n^t}) \right|\F_t\right],\ \ \ t\in (0,T].
$$
By \cref{pro:3.4}, there exists a constant $K>0$ depending only on $(\beta, \gamma, T)$ such that
$$
\psi(\bar Y_t,\ \mu_{\beta,\gamma}(t))\leq \varphi(t,\bar Y_t)\leq \E\left[\left. \varphi(\tau_n^t,\bar Y_{\tau_n^t}) \right|\F_t\right]\leq K \E\left[\left. \psi(\bar Y_{\tau_n^t},\ \mu_{\beta,\gamma}(\tau_n^t))\right|\F_t\right]+K,\ \ t\in (0,T].
$$
And, by virtue of \cref{lem:2.6} we have that for each $n\geq 1$,
$$
|Y_t|\leq\psi\left(|Y_t|,\ \mu_{\beta,\gamma}(t)\right)\leq K \E\left[\left. \psi\left(|Y_{\tau_n^t}|+\int_0^{\tau_n^t} \alpha_s {\rm d}s,\ \mu_{\beta,\gamma}(\tau_n^t)\right)\right|\F_t\right]+K,\ \ \ t\in (0,T],
$$
from which the desired inequality \eqref{eq:3.15} follows for $t\in (0,T]$ by sending $n$ to infinity. Finally, in view of the continuity of $Y_\cdot$ and the martingale in the right side of \eqref{eq:3.15} with respect to the time variable $t$, the inequality \eqref{eq:3.15} holds still true for $t=0$. The proof is complete.
\end{proof}

\begin{rmk}\label{rmk:3.6}
From the above proof, it is not hard to check that in \cref{pro:3.5}, if the $|\xi|$ and $|Y_t|$ are respectively replaced with $\xi^+$ and $Y^+_t$, and Assumption \ref{A:A1} is replaced with the following one
\begin{enumerate}
\renewcommand{\theenumi}{(A1')}
\renewcommand{\labelenumi}{\theenumi}
\item\label{A:A1'} There exist two nonnegative constants $\beta\geq 0$ and $\gamma\geq 0$ with $\beta+\gamma>0$, and an $\R_+$-valued progressively measurable process $(\alpha_t)_{t\in\T}$  such that $\as$,
    \[
    \RE\ (y,z)\in \R_+\times\R^d,\ \ \
    g(\omega,t,y,z)\leq \alpha_t(\omega)+\beta |y|\sqrt{\ln |y|}{\bf 1}_{|y|>1}+\gamma |z|.\vspace{-0.2cm}
    \]
\end{enumerate}
then \cref{pro:3.5} still holds. For this, one only needs to respectively use $Y_\cdot^+$, ${\bf 1}_{Y_\cdot>0}$ and $\frac{1}{2}L_\cdot$ instead of $|Y_\cdot|$, ${\rm sgn} (Y_\cdot)$ and $L_\cdot$ in the above proof.\vspace{0.1cm}
\end{rmk}

Now, we prove \cref{thm:3.3}.

\begin{proof}[Proof of \cref{thm:3.3}]
For each pair of positive integers $n,p\geq 1$, let $\xi^{n,p}$ and $g^{n,p}$ be defined in \eqref{eq:2.1}, and $(Y^{n,p}_\cdot,Z^{n,p}_\cdot)$ be the minimal (maximal) bounded solution of \eqref{eq:2.2}. It is easy to verify that the generator $g^{n,p}$ satisfies Assumption \ref{A:A1} with $\alpha_\cdot$ being replaced with $\alpha_\cdot \wedge (n\vee p)$, where $n\vee p$ denotes the maximum of $n$ and $p$. Then, applying \cref{pro:3.5} to BSDE \eqref{eq:2.2} yields the existence of a constant $C>0$ depending only on $(\beta,\gamma,T)$ such that $\ps$, for each pair of $n,p\geq 1$,
\begin{equation}\label{eq:3.17}
\begin{array}{lll}
|Y^{n,p}_t|&\leq & \Dis \psi(|Y^{n,p}_t|,\ \mu_{\beta,\gamma}(t))\leq C\E\left[\left.\psi\left(|\xi^{n,p}|+\int_0^T[\alpha_t\wedge (n\vee p)]{\rm d}t,\ \mu_{\beta,\gamma}(T)\right)\right|\F_t\right]+C\vspace{0.2cm}\\
&\leq& \Dis C\E\left[\left.\psi\left(|\xi|+\int_0^T \alpha_t {\rm d}t,\ \mu_{\beta,\gamma}(T)\right)\right|\F_t\right]+C\ =:X_t,\ \ \ \ t\in\T.
\end{array}
\end{equation}
Thus, in view of Assumption~\eqref{eq:3.11}, we have found an $\R_+$-valued, progressively measurable and continuous process $(X_t)_{t\in\T}$ such that\vspace{-0.1cm}
$$
\as,\ \ \RE\ n,p\geq 1,\ \ \ |Y^{n,p}_\cdot|\leq X_\cdot.\vspace{-0.1cm}
$$
Now, we can apply \cref{Pro:2.1--stability} to obtain the existence of a progressively measurable process $(Z_t)_{t\in\T} $ such that $(Y_\cdot:=\inf_p\sup_n Y^{n,p}_\cdot, \ Z_\cdot)$ is a solution to BSDE$(\xi,g)$.

Furthermore, sending $n$ and $p$ to infinity in \eqref{eq:3.17} yields \eqref{eq:3.12}, and then $\left(\psi\left(|Y_t|,\mu_{\beta,\gamma}(t)\right)\right)_{t\in\T}$ belongs to class (D). The existence part is then proved.

Finally, the uniqueness of the desired solution under Assumptions \ref{A:UN1} and \ref{A:UN2} is a direct consequence of \cref{pro:2.5}. And, the uniqueness under Assumption \ref{A:UN3} is a consequence of the following \cref{pro:3.7}. The proof is complete.\vspace{0.2cm}
\end{proof}

\begin{pro}\label{pro:3.7}
Let the functions $\psi(\cdot,\cdot)$ and $\mu_{\beta,\gamma}(\cdot)$ be respectively defined in \eqref{eq:2.12} and \eqref{eq:3.8}, $\xi$ and $\xi'$ be two terminal conditions, $g$ and $g'$ be two generators which are continuous in the variables $(y,z)$, and $(Y_t, Z_t)_{t\in\T}$ and $(Y'_t, Z'_t)_{t\in\T}$ be respectively a solution to BSDE $(\xi, g)$ and BSDE $(\xi', g')$ such that both $\left(\psi\left(|Y_t|,\ \mu_{\beta,\gamma}(t)\right)\right)_{t\in\T}$ and $\left(\psi\left(|Y'_t|,\ \mu_{\beta,\gamma}(t)\right)\right)_{t\in\T}$
belong to class (D).

Assume that $\ps$, $\xi\leq \xi'$. If $g$ (resp. $g'$) satisfies Assumptions~\ref{A:UN3} and~\ref{A:A1} with parameters $(\alpha_\cdot,\beta,\gamma)$ such that $\psi\left(\int_0^T \alpha_t{\rm d}t, \ \mu_{\beta,\gamma}(T)\right)\in L^1$, and $\as$,
\begin{equation}\label{eq:3.18}
g(t,Y'_t,Z'_t)\leq g'(t,Y'_t,Z'_t)\ \ \ ({\rm resp.}\  \ g(t,Y_t,Z_t)\leq g'(t,Y_t,Z_t)\ ),
\end{equation}
then $\ps$, $Y_t\leq Y'_t$ for all $t\in\T$.
\end{pro}

\begin{proof}
We first consider the case that the generator $g$ is convex in the state variables $(y,z)$, satisfies Assumption \ref{A:A1} with parameters $(\alpha_\cdot, \beta, \gamma)$ such that $\psi\left(\int_0^T \alpha_t{\rm d}t, \ \mu_{\beta,\gamma}(T)\right)\in L^1$, and $\as$,
$$g(t,Y'_t,Z'_t)\leq g'(t,Y'_t,Z'_t).$$

To use the convexity condition of the generator $g$, we use the $\theta$-technique developed in for example \cite{BriandHu2008PTRF}. For each fixed $\theta\in (0,1)$, define
\begin{equation}\label{eq:3.19}
\Delta^\theta U_\cdot:=\frac{Y_\cdot-\theta Y'_\cdot}{1-\theta}\ \  {\rm and} \ \ \Delta^\theta V_\cdot:=\frac{Z_\cdot-\theta Z'_\cdot}{1-\theta}.
\end{equation}
Then the pair $(\Delta^\theta U_\cdot,\Delta^\theta V_\cdot)$ satisfies the following BSDE:\vspace{0.1cm}
\begin{equation}\label{eq:3.20}
  \Delta^\theta U_t=\Delta^\theta U_T +\int_t^T \Delta^\theta g (s,\Delta^\theta U_s,\Delta^\theta V_s) {\rm d}s-\int_t^T \Delta^\theta V_s \cdot {\rm d}B_s, \ \ \ \ t\in\T,\vspace{0.1cm}
\end{equation}
where $\ass$, for each $(y,z)\in \R\times\R^d$,
\begin{equation}\label{eq:3.21}
\begin{array}{lll}
\Dis \Delta^\theta g(s,y,z)&:=& \Dis \frac{1}{1-\theta}\left[\  g(s,(1-\theta)y+\theta Y'_s,(1-\theta)z+\theta Z'_s)-\theta g(s, Y'_s, Z'_s)\ \right]\vspace{0.2cm}\\
&& \Dis +\frac{\theta}{1-\theta}\left[\  g(s,Y'_s, Z'_s)-g'(s,Y'_s, Z'_s)\ \right].
\end{array}
\end{equation}
It follows from the assumptions that $\ass$,
\begin{equation}\label{eq:3.22}
\RE\ (y,z)\in \R_+\times \R^d,\ \ \  \Delta^\theta g(s,y,z)\leq g(s,y,z)\leq \alpha_s+\beta |y|\sqrt{\ln |y|}{\bf 1}_{|y|>1}+\gamma |z|,
\end{equation}
which means that the generator $\Delta^\theta g$ satisfies Assumption \ref{A:A1'} defined in \cref{rmk:3.6}.

On the other hand, since both processes $\left(\psi\left(|Y_t|,\mu_{\beta,\gamma}(t)\right)\right)_{t\in\T}$ and $\left(\psi\left(|Y'_t|,\mu_{\beta,\gamma}(t)\right)\right)_{t\in\T}$
belong to class (D), and $\psi\left(\int_0^T \alpha_t{\rm d}t, \ \mu_{\beta,\gamma}(T)\right)\in L^1$, by virtue of \cref{lem:2.6},  we conclude that the process
\[
\psi\left(\left(\Delta^\theta U_t\right)^+ +\int_0^t\alpha_s {\rm d}s,\  \mu_{\beta,\gamma}(t)\right),\ \ \ t\in\T
\]
belongs to class (D) for each $\theta\in (0,1)$. Thus, for BSDE \eqref{eq:3.20}, by virtue of \eqref{eq:3.22}, \cref{rmk:3.6} and the proof of \cref{pro:3.5}, we derive that there exists a constant $C>0$ depending on $(\beta,\gamma,T)$ such that
\begin{equation}\label{eq:3.23}
\left(\Delta^\theta U_t\right)^+\leq \psi\left(\left(\Delta^\theta U_t\right)^+,\ \mu_{\beta,\gamma}(t)\right)\leq C\E\left[\left.\psi\left(\left(\Delta^\theta U_T\right)^+ +\int_0^T\alpha_s {\rm d}s,\  \mu_{\beta,\gamma}(T)\right)\right|\F_t\right]+C, \ \ \ t\in\T.
\end{equation}

Furthermore, since
\begin{equation}\label{eq:3.24}
\left(\Delta^\theta U_T\right)^+=\frac{(\xi-\theta \xi')^+}{1-\theta}=\frac{\left[\xi-\theta \xi+\theta(\xi-\xi')\right]^+}{1-\theta}\leq \xi^+,
\end{equation}
we have from \eqref{eq:3.23} that
$$
\left(Y_t-\theta Y'_t\right)^+ \leq (1-\theta)
\left(C\E\left[\left.\psi\left(\xi^+ +\int_0^T\alpha_s {\rm d}s,\  \mu_{\beta,\gamma}(T)\right)\right|\F_t\right]+C
\right), \quad t\in\T.
$$
Thus, the desired conclusion follows in the limit as  $\theta\To 1$. \vspace{0.1cm}

For the case that the generator $g$ is concave in the state variables $(y,z)$, we need to respectively use the $\theta Y_\cdot-Y'_\cdot$ and $\theta Z_\cdot-Z'_\cdot$ to replace $Y_\cdot-\theta Y'_\cdot$ and $Z_\cdot-\theta Z'_\cdot$ in \eqref{eq:3.19} . In this case, the generator $\Delta^\theta g$ in \eqref{eq:3.21} should be replaced with
\[
\begin{array}{lll}
\Dis \Delta^\theta g(s,y,z)&:=& \Dis \frac{1}{1-\theta}\left[\  \theta g(s, Y_s, Z_s) -g(s,-(1-\theta)y+\theta Y_s,-(1-\theta)z+\theta Z_s)\ \right]\vspace{0.2cm}\\
&& \Dis +\frac{1}{1-\theta}\left[\  g(s,Y'_s, Z'_s)-g'(s,Y'_s, Z'_s)\ \right].
\end{array}
\]
Since $g$ is concave in $(y,z)$, we have $\ass$,
$$
\RE\ (y,z)\in \R\times \R^d,\ \ g(s,-(1-\theta)y+\theta Y_s,-(1-\theta)z+\theta Z_s)\geq \theta g(s,Y_s,Z_s)+(1-\theta)g(t,-y,-z),
$$
and then, \eqref{eq:3.22} can be replaced by
\[
\RE\ (y,z)\in \R_+\times \R^d,\ \ \ \Delta^\theta g(s,y,z)\leq -g(s,-y,-z)\leq \alpha_s+\beta |y|\sqrt{\ln |y|}{\bf 1}_{|y|>1}+\gamma |z|,
\]
which means that the generator $\Delta^\theta g$ still satisfies Assumption~\ref{A:A1'}. Consequently, \eqref{eq:3.23} still holds. Moreover, we use
\[
\left(\Delta^\theta U_T\right)^+=\frac{(\theta \xi-\xi')^+}{1-\theta}=\frac{\left[\theta \xi- \xi+(\xi-\xi')\right]^+}{1-\theta}\leq (-\xi)^+=\xi^-,\vspace{0.1cm}
\]
instead of \eqref{eq:3.24}, and it follows from \eqref{eq:3.23} that
$$
\left(\theta Y_t-Y'_t\right)^+ \leq (1-\theta)
\left(C\E\left[\left.\psi\left(\xi^- +\int_0^T\alpha_s {\rm d}s,\  \mu_{\beta,\gamma}(T)\right)\right|\F_t\right]+C
\right), \ \ \ t\in\T.
$$
Thus, the desired conclusion follows in the limit as $\theta\To 1$.\vspace{0.1cm}

Finally, in the same way as above, we can prove the desired conclusion under assumptions that the generator $g'$ satisfies Assumptions~\ref{A:UN3} and \ref{A:A1} with parameters $(\alpha_\cdot, \beta, \gamma)$ such that $\psi\left(\int_0^T \alpha_t{\rm d}t, \ \mu_{\beta,\gamma}(T)\right)\in L^1$, and $\as$, $g(t,Y_t,Z_t)\leq g'(t,Y_t,Z_t)$. The proof is then complete.
\end{proof}

\begin{ex}\label{exp:3.8}
Let $\beta>0$, $\gamma>0$ and $k\geq 0$. For each $(\omega,t,y,z)\in \Omega\times\T\times\R\times\R^d$, define
$$
g_1(\omega,t,y,z):=|B_t(\omega)|-ke^y+\beta |y|\sqrt{|\ln |y||}{\bf 1}_{|y|\leq {1\over \sqrt{e}}}+{\beta\over \sqrt{2e}}{\bf 1}_{|y|>{1\over \sqrt{e}}}+\gamma \left(|z|-k\sqrt{|z|}\right)
$$
and
$$
g_2(\omega,t,y,z):=|B_t(\omega)|+ke^{-y}+\beta |y|\sqrt{\ln |y|}{\bf 1}_{|y|\geq 1}+\gamma |z|.
$$
It is easy to verify that the generator $g_1$ satisfies Assumptions \ref{A:EX1}, \ref{A:A1}, \ref{A:UN1} and \ref{A:UN2}, and the generator $g_2$ satisfies Assumptions \ref{A:EX1}, \ref{A:A1} and \ref{A:UN3}.
\end{ex}

\section{The case of $L\exp\left(\mu (\ln (1+L))^p\right)$-integrability terminal condition with $p\in (1/2,1)$}
\label{sec:4-mainresult--4}
\setcounter{equation}{0}

In this section, we always assume that the generator $g$ satisfies the following assumption,  which is strictly weaker than \eqref{eq:1.2} with $\lambda \in (0,1/2)$ and $\delta=\lambda +1/2$.
\begin{enumerate}
\renewcommand{\theenumi}{(A2)}
\renewcommand{\labelenumi}{\theenumi}
\item \label{A:A2} There exist three constants $\beta\geq 0$, $\gamma > 0$ and $\lambda\in (0,1/2)$, and an $\R_+$-valued progressively measurable process $(\alpha_t)_{t\in\T}$  such that $\as$,
    $$
    \RE\ (y,z)\in \R\times\R^d,\ \ \
    {\rm sgn}(y)g(\omega,t,y,z)\leq \alpha_t(\omega)+\beta|y|(\ln |y|)^{\lambda+\frac{1}{2}}{\bf 1}_{|y|>1}+\gamma |z||\ln |z||^{\lambda}.
    $$
  \end{enumerate}

For each $\eps>0$, by letting $k=1+\eps$ in \cref{Pro:2.2Inequality} we know that there exists a constant $C_{\lambda,\eps}>0$ depending only on parameters $(\lambda,\eps)$ such that \eqref{eq:2.5} holds. Then, our objective is to search for a strictly increasing and strictly convex function $\phi$ such that for each $(s,x)\in \T\times \R_+$, with $\beta\geq0$ and $\gamma>0$,\vspace{0.1cm}
\begin{equation}\label{eq:4.1}
-\beta\phi_x(s, x) x(\ln x)^{\lambda+\frac{1}{2}} {\bf 1}_{x>1}-{\gamma^2\over 2}\frac{\left(\phi_x(s, x)\right)^2}
{\phi_{xx}(s, x)}\left[(1+\eps)\left|\ln \frac{\gamma\phi_x(s, x)}{\phi_{xx}(s, x)}
\right|^{2\lambda}+C_{\lambda,\eps}\right]+\phi_s(s, x)\geq 0,\vspace{0.1cm}
\end{equation}
which is just \eqref{eq:2.11} with $\lambda\in (0,1/2)$, $\delta=\lambda+1/2$, $\bar\delta=0$, $k=1+\eps$ and $C_{k,\lambda}=C_{\lambda, \eps}$.\vspace{0.2cm}

Now, let $\eps>0$, $\beta\geq0$, $\gamma>0$, $\lambda\in (0,1/2)$ and $\mu_s,\nu_s:\T\To\R_+$ be two increasing and continuously differential functions with $\mu_0=\eps$ and $\nu_0=0$. Define\vspace{-0.1cm}
\begin{equation}\label{eq:4.2}
\Dis k_{\eps}:=\exp \left({\frac{1+\eps}{2\eps}}\right)+\left(\frac{1}{2\eps^2}
\right)^{\frac{1}{2\eps}}+\frac{\mu_T}{\gamma}+{\rm e}.\vspace{-0.1cm}
\end{equation}
We choose the following function
$$
\phi(s,x;\eps):=(x+k_\eps)\exp\left(\mu_s \left(\ln \left(x+k_{\eps}\right)\right)^{\lambda+\frac{1}{2}}+\nu_s\right), \ \ \ (s,x)\in \T\times \R_+
$$
to explicitly solve the inequality \eqref{eq:4.1}. For each $(s,x)\in \T\times \R_+$, a simple computation gives
\begin{equation}\label{eq:4.3}
\phi_x(s,x;\eps)=\phi(s,x;\eps){\left(\lambda+\frac{1}{2}\right)\mu_s+\left(\ln (x+k_\eps)\right)^{\frac{1}{2}-\lambda} \over (x+k_\eps)\left(\ln (x+k_\eps)\right)^{\frac{1}{2}-\lambda}}>0,\vspace{0.1cm}
\end{equation}
\begin{equation}\label{eq:4.4}
\phi_{xx}(s,x;\eps)=\phi(s,x;\eps){\left(\lambda+\frac{1}{2}\right)\mu_s\left[
\ln (x+k_\eps)+\left(\lambda+\frac{1}{2}\right)\mu_s \left(\ln (x+k_\eps)\right)^{\lambda+\frac{1}{2}}-\left(\frac{1}{2}-\lambda\right)\right]\over (x+k_\eps)^2\left(\ln (x+k_\eps)\right)^{\frac{3}{2}-\lambda}}>0
\end{equation}
and\vspace{-0.1cm}
\begin{equation}\label{eq:4.5}
\phi_s(s,x;\eps)=\phi(s,x;\eps)\left(\mu'_s\left(\ln \left(x+k_{\eps}\right)\right)^{\lambda+\frac{1}{2}}+\nu'_s \right)>0.\vspace{0.1cm}
\end{equation}
Furthermore, from the definition of $k_\eps$ it can be directly verified that for each $(s,x)\in \T\times \R_+$,
\begin{equation}\label{eq:4.6}
\begin{array}{ll}
& \Dis \ln (x+k_\eps)+\left(\lambda+\frac{1}{2}\right)\mu_s \left(\ln (x+k_\eps)\right)^{\lambda+\frac{1}{2}}-\left(\frac{1}{2}-\lambda\right)\vspace{0.2cm}
\\
\geq & \Dis \frac{1}{1+\eps}\left[
\ln (x+k_\eps)+\left(\lambda+\frac{1}{2}\right)\mu_s \left(\ln (x+k_\eps)\right)^{\lambda+\frac{1}{2}}\right]
\end{array}
\end{equation}
and
\begin{equation}\label{eq:4.7}
\left(\ln (x+k_\eps)\right)^{\frac{1}{2}-\lambda}\leq \sqrt{\ln (x+k_\eps)}\leq (x+k_\eps)^\eps.
\end{equation}
It follows from \eqref{eq:4.4} and \eqref{eq:4.6} that \vspace{0.1cm}
$$
\phi_{xx}(s,x;\eps)\geq\phi(s,x;\eps)
{\left(\lambda+\frac{1}{2}\right)\mu_s\left[\left(\lambda+\frac{1}{2}\right)\mu_s+
\left(\ln (x+k_\eps)\right)^{\frac{1}{2}-\lambda}\right]\over(1+\eps)(x+k_\eps)^2\left(\ln (x+k_\eps)\right)^{1-2\lambda}},\ \ (s,x)\in \T\times \R_+.\vspace{0.2cm}
$$
Then, for each $(s,x)\in \T\times \R_+$, in view of~\eqref{eq:4.3}, we have\vspace{0.1cm}
\begin{equation}\label{eq:4.8}
{\gamma^2\over 2}\frac{\left(\phi_x(s,x;\eps)\right)^2}{\phi_{xx}(s,x;\eps)}\leq {\gamma^2(1+\eps)\over 2}\phi(s,x;\eps)\left(1+\frac{\left(\ln (x+k_\eps)\right)^{\frac{1}{2}-\lambda}}{\left(\lambda+\frac{1}{2}\right)\mu_s} \right)
\end{equation}
and, in view of \eqref{eq:4.2}, \eqref{eq:4.4} and \eqref{eq:4.7},\vspace{0.1cm}
\[
\left\{
\begin{array}{l}
\Dis \frac{\gamma\phi_x(s,x;\eps)}{\phi_{xx}(s,x;\eps)}\geq {\gamma\over \left(\lambda+\frac{1}{2}\right)\mu_s}(x+k_\eps)\left(\ln (x+k_\eps)\right)^{\frac{1}{2}-\lambda}\geq {\gamma \over \mu_T}k_\eps\geq 1;\vspace{0.3cm}\\
\Dis \frac{\gamma\phi_x(s,x;\eps)}{\phi_{xx}(s,x;\eps)}\leq {\gamma (1+\eps)\over \left(\lambda+\frac{1}{2}\right)\mu_s}(x+k_\eps)\left(\ln (x+k_\eps)\right)^{\frac{1}{2}-\lambda}\leq {2\gamma (1+\eps)\over \eps}(x+k_\eps)^{1+\eps},\vspace{0.1cm}
\end{array}
\right.
\]
which yields the following
\begin{equation}\label{eq:4.9}
\left|\ln \frac{\gamma\phi_x(s,x;\eps)}{\phi_{xx}(s,x;\eps)}\right|^{2\lambda}\leq
 \left|\ln {2\gamma (1+\eps)\over \eps}\right|^{2\lambda}+ (1+\eps)^{2\lambda}\left( \ln (x+k_\eps)\right)^{2\lambda}.\vspace{0.3cm}
\end{equation}

In the sequel, observe that\vspace{-0.1cm}
$$
\RE\ x\in\R_+,\ \ x(\ln x)^{\lambda+\frac{1}{2}} {\bf 1}_{x>1}\leq (x+k_\eps)\left(\ln (x+k_\eps)\right)^{\lambda+\frac{1}{2}}.
$$
We substitute \eqref{eq:4.3}, \eqref{eq:4.8}, \eqref{eq:4.9} and \eqref{eq:4.5} into the left side of \eqref{eq:4.1} with $\phi(s, x;\eps)$ instead of $\phi(s,x)$ to get \vspace{0.1cm}
$$
\begin{array}{ll}
&\Dis -\beta\phi_x(s,x;\eps) x(\ln x)^{\lambda+\frac{1}{2}} {\bf 1}_{x>1}-{\gamma^2\over 2}\frac{\left(\phi_x(s, x;\eps)\right)^2}
{\phi_{xx}(s, x;\eps)}\left[(1+\eps)\left|\ln \frac{\gamma\phi_x(s, x;\eps)}{\phi_{xx}(s, x;\eps)}
\right|^{2\lambda}+C_{\lambda,\eps}\right]+\phi_s(s, x;\eps)\vspace{0.2cm}\\
\geq & \Dis -\beta \phi(s, x;\eps) \left[\left(\lambda+\frac{1}{2}\right)\mu_s\left(\ln (x+k_\eps)\right)^{2\lambda}+\left(\ln (x+k_\eps)\right)^{\lambda+\frac{1}{2}}\right]\vspace{0.2cm}\\
& \Dis
-{\gamma^2(1+\eps)\over 2}\phi(s,x;\eps)\left(1+\frac{\left(\ln (x+k_\eps)\right)^{\frac{1}{2}-\lambda}}{\left(\lambda+\frac{1}{2}\right)\mu_s} \right)\left[(1+\eps)^{2\lambda+1}\left( \ln (x+k_\eps)\right)^{2\lambda}+\bar C_{\lambda,\eps}\right]
\vspace{0.2cm}\\
& \Dis +\phi(s,x;\eps)\left(\mu'_s\left(\ln \left(x+k_{\eps}\right)\right)^{\lambda+\frac{1}{2}}+\nu'_s \right), \ \ \ \ (s,x)\in \T\times \R_+\vspace{0.1cm}
\end{array}
$$
where
$$
\bar C_{\lambda,\eps}:=(1+\eps)\left|\ln {2\gamma (1+\eps)\over \eps}\right|^{2\lambda}+C_{\lambda,\eps}.\vspace{0.2cm}
$$
Furthermore, in view of the facts that $\mu_\cdot\geq \eps$ and $\lambda\in (0,1/2)$, using Young's inequality,  we see that there is a constant $C_{\beta,\lambda,\eps}^1>0$ depending only on $(\beta,\lambda,\eps)$, and a constant $C_{\gamma,\lambda,\eps}^2>0$ depending only on $(\gamma,\lambda,\eps,\bar C_{\lambda,\eps})$ such that for each $(s,x)\in \T\times \R_+$,
$$
\beta\left(\lambda+\frac{1}{2}\right)\left(\ln (x+k_\eps)\right)^{2\lambda}\leq \eps\left(\ln (x+k_\eps)\right)^{\lambda+\frac{1}{2}}+C_{\beta,\lambda,\eps}^1
$$
and
$$
\begin{array}{ll}
&\Dis {\gamma^2(1+\eps)\over 2}\left(1+\frac{\left(\ln (x+k_\eps)\right)^{\frac{1}{2}-\lambda}}{\left(\lambda+\frac{1}{2}\right)\mu_s} \right)\left[(1+\eps)^{2\lambda+1}\left(\ln (x+k_\eps)\right)^{2\lambda}+\bar C_{\lambda,\eps}\right]\vspace{0.2cm}\\
\leq & \Dis \left( \frac{\gamma^2(1+\eps)^{2\lambda+2}}{(2\lambda+1)\mu_s}+\eps\right)\left( \ln (x+k_\eps)\right)^{\lambda+\frac{1}{2}}+C_{\gamma,\lambda,\eps}^2.\vspace{0.1cm}
\end{array}
$$
Then, with the notation $v:=\left(\ln (x+k_\eps)\right)^{\lambda+\frac{1}{2}}$ we have for each $(s,x)\in \T\times \R_+$,\vspace{0.1cm}
$$
\begin{array}{ll}
&\Dis -\beta\phi_x(s,x;\eps) x(\ln x)^{\lambda+\frac{1}{2}} {\bf 1}_{x>1}-{\gamma^2\over 2}\frac{\left(\phi_x(s, x;\eps)\right)^2}
{\phi_{xx}(s, x;\eps)}\left[(1+\eps)\left|\ln \frac{\gamma\phi_x(s, x;\eps)}{\phi_{xx}(s, x;\eps)}
\right|^{2\lambda}+C_{\lambda,\eps}\right]+\phi_s(s, x;\eps)\vspace{0.3cm}\\
\geq & \Dis \phi(s, x;\eps)\left[\left(-\eps \mu_s-\beta
-\frac{\gamma^2(1+\eps)^{2\lambda+2}}{(2\lambda+1)\mu_s}-\eps+\mu'_s \right)v+ \left(-C_{\beta,\lambda,\eps}^1\mu_s- C_{\gamma,\lambda,\eps}^2+\nu'_s\right)\right].\vspace{0.1cm}
\end{array}
$$
Thus, \eqref{eq:4.1} holds if the functions $\mu_s, \nu_s\in \T$ satisfies\vspace{0.2cm}
\begin{equation}\label{eq:4.10}
\mu'_s=\eps\mu_s
+\frac{\gamma^2(1+\eps)^{2\lambda+2}}{2\lambda+1}\frac{1}{\mu_s}+\eps+\beta,\ \ \ s\in [0,T]
\end{equation}
and
\begin{equation}\label{eq:4.11}
\nu_s=C_{\beta,\lambda,\eps}^1\int_0^s \mu_r {\rm d}r+C_{\gamma,\lambda,\eps}^2 s,\ \ \ s\in [0,T].\vspace{0.2cm}
\end{equation}

It is not very hard to verify that for each $\beta\geq 0$, $\gamma > 0$, $\lambda\in (0,1/2)$ and $\eps>0$, there exists a unique strictly increasing and continuous function with $\eps$ at the origin satisfying \eqref{eq:4.10}. We denote this unique solution by $\bar \mu_{\beta,\gamma,\lambda,\eps}(\cdot)$, i.e., $\bar \mu_{\beta,\gamma,\lambda,\eps}(0)=\eps$ and
\begin{equation}\label{eq:4.12}
\bar \mu'_{\beta,\gamma,\lambda,\eps}(s)=\eps
\bar\mu_{\beta,\gamma,\lambda,\eps}(s)
+\frac{\gamma^2(1+\eps)^{2\lambda+2}}{2\lambda+1}
\frac{1}{\bar \mu_{\beta,\gamma,\lambda,\eps}(s)}+\eps+\beta,\ \ \ s\in [0,T].
\end{equation}
We also denote, in view of \eqref{eq:4.11},
\begin{equation}\label{eq:4.13}
\bar \nu_{\beta,\gamma,\lambda,\eps}(s):=C_{\beta,\lambda,\eps}^1\int_0^s \bar \mu_{\beta,\gamma,\lambda,\eps}(r){\rm d}r+C_{\gamma,\lambda,\eps}^2 s,\ \ \ s\in [0,T].\vspace{0.1cm}
\end{equation}

Moreover, it can be directly checked that as $\eps\To 0^+$, the function $\bar \mu_{\beta,\gamma,\lambda,\eps}(\cdot)$ tends decreasingly to the unique solution $\bar \mu_{\beta,\gamma,\lambda}^0(\cdot)$ of the following ODE:
\begin{equation}\label{eq:4.14}
\bar \mu_{\beta,\gamma,\lambda}^0(0)=0\ \ \ {\rm and}\ \ \ \left(\bar \mu_{\beta,\gamma,\lambda}^0(s)\right)'=
\frac{\gamma^2}{2\lambda+1}
\frac{1}{\bar \mu_{\beta,\gamma,\lambda}^0(s)}+\beta,\ \ \ s\in (0,T].\vspace{0.1cm}
\end{equation}

\begin{rmk}\label{rmk:4.1}
In general, $\bar \mu_{\beta,\gamma,\lambda}^0(\cdot)$ does not have an explicit expression. However, it is easy to check that $\bar \mu_{\beta,\gamma,\lambda}^0(s)=\frac{\gamma}{\sqrt{\lambda+{1\over 2}}}\sqrt{s}$ in the case of $\beta=0$. In addition, when $\lambda=0$, we have $\bar \mu_{\beta,\gamma,\lambda}^0(s)=\sqrt{2}\mu_{\beta,\gamma}(s)$, where $\mu_{\beta,\gamma}(s)$ is defined in \eqref{eq:3.8}. \vspace{0.2cm}
\end{rmk}

In summary, we have

\begin{pro}\label{pro:4.2}
For each $\beta\geq 0$, $\gamma > 0$, $\lambda\in (0,1/2)$ and $\eps>0$, define the function
\begin{equation}\label{eq:4.15}
\bar \varphi(s,x;\eps):=(x+\bar k_{\eps})\exp\left(\bar \mu_{\beta,\gamma,\lambda,\eps}(s)\left(\ln \left(x+\bar k_{\eps}\right)\right)^{\lambda+\frac{1}{2}}+\bar \nu_{\beta,\gamma,\lambda,\eps}(s)\right), \ \ \ (s,x)\in [0,T]\times \R_+
\end{equation}
with $\bar \mu_{\beta,\gamma,\lambda,\eps}(\cdot)$ and $\bar \nu_{\beta,\gamma,\lambda,\eps}(\cdot)$ being respectively defined in \eqref{eq:4.12} and \eqref{eq:4.13}, and
\[
\bar k_{\eps}:=\exp \left({\frac{1+\eps}{2\eps}}\right)+\left(\frac{1}{2\eps^2}
\right)^{\frac{1}{2\eps}}+\frac{\bar \mu_{\beta,\gamma,\lambda,\eps}(T)}{\gamma}+\rm {e},
\]
which is $k_\eps$ in \eqref{eq:4.2} with $\bar \mu_{\beta,\gamma,\lambda,\eps}(T)$ instead of $\mu_T$. Then we have
\begin{itemize}
\item [(i)] $\bar \varphi(\cdot,\cdot;\eps)\in \mathcal {C}^{1,2}\left([0,T]\times \R_+\right)$;
\item [(ii)] $\bar \varphi(\cdot,\cdot;\eps)$ satisfies the inequality \eqref{eq:2.3} with $\lambda\in (0,1/2)$, $\delta=\lambda+1/2$ and $\bar\delta=0$. i.e., for each $(s,x,z)\in [0,T]\times \R_+\times\R^d$, we have\vspace{-0.1cm}
$$
-\bar \varphi_x(s,x;\eps)\left(\beta x(\ln x)^{\lambda+{1\over 2}}{\bf 1}_{x>1}+\gamma |z||\ln |z||^{\lambda}\right)+{1\over 2}\bar \varphi_{xx}(s,x;\eps)|z|^2+\bar \varphi_s(s,x;\eps)\geq 0.
$$
\end{itemize}
\end{pro}

Now, for $\lambda\in (0,1/2)$, we define the function
\begin{equation}\label{eq:4.16}
\bar \psi(x,\mu):=x\exp\left(\mu\left(\ln (1+x)\right)^{\lambda+{1\over 2}}\right),\ \ (x,\mu)\in \R_+\times \R_+.
\end{equation}
The following existence and uniqueness theorem is one main result of this section.

\begin{thm}\label{thm:4.3}
Let the functions $\bar \mu_{\beta,\gamma,\lambda,\eps}(\cdot)$, $\bar \mu_{\beta,\gamma,\lambda}^0(\cdot)$ and $\bar \psi(\cdot,\cdot)$ be respectively defined in \eqref{eq:4.12}, \eqref{eq:4.14}, and \eqref{eq:4.16}, $\xi$ be a terminal condition and $g$ be a generator which is continuous in the state variables $(y,z)$. If $g$ satisfies Assumptions \ref{A:EX1} and \ref{A:A2} with parameters $(\alpha_\cdot, \beta, \gamma, \lambda)$, and there exists a positive constant $\mu>\bar \mu_{\beta,\gamma,\lambda}^0(T)$ such that\vspace{-0.1cm}
\begin{equation}\label{eq:4.17}
\E\left[\bar \psi\left(|\xi|+\int_0^T \alpha_t {\rm d}t,\ \mu\right)\right]<+\infty,\vspace{0.1cm}
\end{equation}
then BSDE$(\xi,g)$ admits a solution $(Y_t,Z_t)_{t\in\T}$ such that for some $\eps>0$, $\left(\bar\psi\left(|Y_t|,\ \bar \mu_{\beta,\gamma,\lambda,\eps}(t)\right)\right)_{t\in\T}$ belongs to class (D), and $\ps$,
\begin{equation}\label{eq:4.18}
|Y_t|\leq \bar \psi(|Y_t|,\ \bar \mu_{\beta,\gamma,\lambda,\eps}(t))\leq C\E\left[\left.\bar \psi\left(|\xi|+\int_0^T\alpha_t {\rm d}t,\ \mu\right)\right|\F_t\right]+C,\ \ \ t\in\T,
\end{equation}
where $C$ is a positive constant depending only on $(\beta,\gamma,\lambda,\eps,T)$, and $\eps$ is the unique positive constant satisfying $\bar \mu_{\beta,\gamma,\lambda,\eps}(T)=\mu$.

Moreover, if $g$ satisfies either Assumptions \ref{A:UN1} and \ref{A:UN2} or Assumption \ref{A:UN3}, then BSDE $(\xi,g)$ admits a unique solution $(Y_t,Z_t)_{t\in\T}$ such that $\left(\bar\psi\left(|Y_t|,\ \bar \mu_{\beta,\gamma,\lambda,\eps}(t)\right)\right)_{t\in\T}$ belongs to class (D).
\end{thm}

To prove this theorem, we need the following two propositions. First, similar to \cref{pro:3.4},  we can prove the following

\begin{pro}\label{pro:4.4}
Let the functions $\bar \mu_{\beta,\gamma,\lambda,\eps}(s)$, $\bar \varphi(s,x;\eps)$ and $\bar \psi(x,\mu)$ be respectively defined on \eqref{eq:4.12}, \eqref{eq:4.15} and \eqref{eq:4.16}. Then, there exists a universal constant $K>0$ depending on $(\beta,\gamma,\lambda,\eps,T)$ such that\vspace{-0.1cm}
$$
\RE\ (s,x)\in [0,T]\times\R_+,\ \ \ \bar \psi(x,\bar \mu_{\beta,\gamma,\lambda,\eps}(s))\leq \bar \varphi(s,x;\eps)\leq K\bar \psi(x,\bar \mu_{\beta,\gamma,\lambda,\eps}(s))+K.
$$
\end{pro}

Similar to the proof of \cref{pro:3.5},  we can prove the following

\begin{pro}\label{pro:4.5}
Let the functions $\bar \mu_{\beta,\gamma,\lambda,\eps}(\cdot)$ and $\bar \psi(\cdot,\cdot)$ be respectively defined in \eqref{eq:4.12} and \eqref{eq:4.16}, $\xi$ be a terminal condition and $g$ be a generator which is continuous in $(y,z)$. If $g$ satisfies Assumption \ref{A:A2} with parameters $(\alpha_\cdot, \beta, \gamma, \lambda)$, $|\xi|+\int_0^T \alpha_t {\rm d}t$ is a bounded random variable, and $(Y_t,Z_t)_{t\in\T}$ is a solution of BSDE$(\xi,g)$ such that $Y_\cdot$ is bounded, then for each $\eps>0$, $\ps$,
\begin{equation}\label{eq:4.19}
|Y_t|\leq \bar\psi(|Y_t|,\ \bar \mu_{\beta,\gamma,\lambda,\eps}(t))\leq C\E\left[\left.\bar\psi\left(|\xi|+\int_0^T\alpha_t {\rm d}t,\ \bar \mu_{\beta,\gamma,\lambda,\eps}(T)\right)\right|\F_t\right]+C,\ \ \ t\in \T,
\end{equation}
where $C$ is a positive constant depending only on $(\beta,\gamma,\lambda,\eps,T)$.
\end{pro}

\begin{proof}
Define
$$
\bar Y_t:=|Y_t|+\int_0^t \alpha_s {\rm d}s\ \ \ \
{\rm and}\ \ \ \ \bar Z_t:={\rm sgn}(Y_t)Z_t,\ \ \ \ t\in \T.
$$
Using It\^{o}-Tanaka's formula, we have
$$
\bar Y_t=\bar Y_T+\int_t^T \left({\rm sgn}(Y_s)g(s,Y_s,Z_s)-\alpha_s\right){\rm d}s-\int_t^T \bar Z_s \cdot {\rm d}B_s-\int_t^T {\rm d}L_s, \ \ \ t\in\T,
$$
with $L_\cdot$ being the local time of $Y$ at the origin. Now, we fix $\eps>0$ and apply It\^{o}-Tanaka's formula to the process $\bar\varphi(s, \bar Y_s;\eps)$ (see~\eqref{eq:4.15} for the function $\bar\varphi(\cdot,\cdot;\eps)$), to derive, in view of \ref{A:A2},
$$
\begin{array}{lll}
\Dis {\rm d}\bar\varphi(s,\bar Y_s;\eps)
&=&\Dis \bar\varphi_x(s,\bar Y_s;\eps)
\left(-{\rm sgn}(Y_s)g(s,Y_s,Z_s)+\alpha_s\right){\rm d}s+\bar\varphi_x(s,\bar Y_s;\eps)\bar Z_s \cdot {\rm d}B_s+\bar\varphi_x(s,\bar Y_s;\eps){\rm d}L_s\vspace{0.1cm}\\
&&\Dis +{1\over 2}\bar\varphi_{xx}(s,\bar Y_s;\eps)|Z_s|^2{\rm d}s+\bar\varphi_s(s,\bar Y_s;\eps){\rm d}s\vspace{0.2cm}\\
&\geq &\Dis \bigg[-\bar\varphi_x(s,\bar Y_s;\eps)\left(\beta |Y_s|(\ln |Y_s|)^{\lambda+{1\over 2}}{\bf 1}_{|Y_s|>1}+\gamma |Z_s||\ln |Z_s||^{\lambda}\right)\vspace{0.1cm}\\
&&\Dis
\ \ +{1\over 2}\bar\varphi_{xx}(s,\bar Y_s;\eps)|Z_s|^2+\varphi_s(s,\bar Y_s;\eps)\bigg]{\rm d}s+\bar\varphi_x(s,\bar Y_s;\eps)\bar Z_s \cdot {\rm d}B_s, \ \ s\in\T.
\end{array}
$$
Furthermore, in view of the fact that $|Y_s|(\ln |Y_s|)^{\lambda+{1\over 2}}{\bf 1}_{|Y_s|>1}\leq \bar Y_s(\ln \bar Y_s)^{\lambda+{1\over 2}}{\bf 1}_{\bar Y_s>1}$ and \cref{pro:4.2},  we have
\begin{equation}\label{eq:4.20}
{\rm d}\bar\varphi(s,\bar Y_s;\eps)\geq \bar\varphi_x(s,\bar Y_s;\eps)\bar Z_s \cdot {\rm d}B_s,\ \ s\in \T.
\end{equation}
Let us consider, for each integer $n\geq 1$ and each $t\in \T$,  the following stopping time
$$
\tau_n^t:=\inf\left\{s\in [t,T]: \int_t^s \left[\bar \varphi_x(r,\bar Y_r;\eps)\right]^2|\bar Z_r|^2{\rm d}r\geq n \right\}\wedge T.
$$
It follows from the definition of $\tau_n^t$ and the inequality \eqref{eq:4.20} that for each $n\geq 1$,
$$
\bar\varphi(t,\bar Y_t;\eps)\leq \E\left[\left. \bar\varphi(\tau_n^t,\bar Y_{\tau_n^t};\eps) \right|\F_t\right],\ \ \ t\in \T.
$$
Thus, by \cref{pro:4.4}, there exists a constant $K>0$ depending only on $(\beta, \gamma, \lambda,\eps,T)$ such that
$$
\bar\psi(\bar Y_t,\ \bar \mu_{\beta,\gamma,\lambda,\eps}(t))\leq \bar\varphi(t,\bar Y_t;\eps)\leq \E\left[\left. \bar\varphi(\tau_n^t,\bar Y_{\tau_n^t};\eps) \right|\F_t\right]\leq K \E\left[\left. \bar\psi(\bar Y_{\tau_n^t},\ \bar \mu_{\beta,\gamma,\lambda,\eps}(\tau_n^t))\right|\F_t\right]+K,\ \ t\in\T.
$$
And, since $\bar\psi(x,\mu)$ is increasing in $x$, we have that for each $n\geq 1$,
$$
|Y_t|\leq \bar\psi\left(|Y_t|,\ \bar \mu_{\beta,\gamma,\lambda,\eps}(t)\right)\leq K \E\left[\left. \bar\psi\left(|Y_{\tau_n^t}|+\int_0^{\tau_n^t} \alpha_s {\rm d}s,\ \bar \mu_{\beta,\gamma,\lambda,\eps}(\tau_n^t)\right)\right|\F_t\right]+K,\ \ \ t\in \T,
$$
from which the desired inequality \eqref{eq:4.19} follows by sending $n$ to infinity. The proof is completed.
\end{proof}

\begin{rmk}\label{rmk:4.6}
From the preceding proof, we can check that Assertions of  \cref{pro:4.5} are still true if $(|\xi|, |Y_t|)$ is  replaced with $(\xi^+, Y^+_t)$, and Assumption~\ref{A:A2} is replaced with the following one~\ref{A:A2'}:
\begin{enumerate}
\renewcommand{\theenumi}{(A2')}
\renewcommand{\labelenumi}{\theenumi}
\item\label{A:A2'} There exist three constants $\beta\geq 0$, $\gamma > 0$ and $\lambda\in (0,1/2)$, and an $\R_+$-valued progressively measurable process $(\alpha_t)_{t\in\T}$  such that $\as$,
    \[
    \RE\ (y,z)\in \R_+\times\R^d,\ \ \
    g(\omega,t,y,z)\leq \alpha_t(\omega)+\beta|y|(\ln |y|)^{\lambda+\frac{1}{2}}{\bf 1}_{|y|>1}+\gamma |z||\ln |z||^{\lambda}.\vspace{-0.2cm}
    \]
\end{enumerate}
For this,  it is sufficient to use $(Y_\cdot^+, {\bf 1}_{Y_\cdot>0}, \frac{1}{2}L_\cdot)$ instead of $(|Y_\cdot|,  {\rm sgn} (Y_\cdot),  L_\cdot)$ in the proof.\vspace{0.1cm}
\end{rmk}

Now, we prove~\cref{thm:4.3}.

\begin{proof}[Proof of \cref{thm:4.3}]
For each pair of positive integers $n,p\geq 1$, let $\xi^{n,p}$ and $g^{n,p}$ be defined in \eqref{eq:2.1}, and $(Y^{n,p}_\cdot,Z^{n,p}_\cdot)$ be the minimal (maximal) bounded solution of \eqref{eq:2.2}. It is easy to check that the generator $g^{n,p}$ satisfies Assumption \ref{A:A2} with $\alpha_\cdot$ being replaced with $\alpha_\cdot \wedge (n\vee p)$.

Now, we assume that there exists a positive constant $\mu>\bar \mu_{\beta,\gamma,\lambda}^0(T)$ such that \eqref{eq:4.17} holds. From the definitions of $\bar \mu_{\beta,\gamma,\lambda,\eps}(\cdot)$ and $\bar \mu_{\beta,\gamma,\lambda}^0(\cdot)$ in \eqref{eq:4.12} and \eqref{eq:4.14}, it is not very difficult to find a (unique) positive constant $\eps>0$ satisfying $\bar \mu_{\beta,\gamma,\lambda,\eps}(T)=\mu$. Then, applying \cref{pro:4.5} with this $\eps$ to BSDE \eqref{eq:2.2} yields the existence of a $C>0$ depending only on $(\beta,\gamma,\lambda,\eps,T)$ such that $\ps$, for all $n,p\geq 1$,\vspace{0.1cm}
\begin{equation}\label{eq:4.21}
\begin{array}{lll}
|Y^{n,p}_t|&\leq & \Dis \bar\psi(|Y^{n,p}_t|,\ \bar \mu_{\beta,\gamma,\lambda,\eps}(t))\leq C\E\left[\left.\bar\psi\left(|\xi^{n,p}|+\int_0^T[\alpha_t\wedge (n\vee p)]{\rm d}t,\ \bar \mu_{\beta,\gamma,\lambda,\eps}(T)\right)\right|\F_t\right]+C\vspace{0.2cm}\\
&\leq& \Dis C\E\left[\left.\bar\psi\left(|\xi|+\int_0^T \alpha_t {\rm d}t,\ \bar \mu_{\beta,\gamma,\lambda,\eps}(T)\right)\right|\F_t\right]+C\vspace{0.3cm}\\
&=& \Dis C\E\left[\left.\bar\psi\left(|\xi|+\int_0^T \alpha_t {\rm d}t,\ \mu\right)\right|\F_t\right]+C\ =:X_t,\ \ \ \ t\in\T.\vspace{-0.1cm}
\end{array}
\end{equation}
Thus, in view of \eqref{eq:4.17}, we have found an $\R_+$-valued, progressively measurable and continuous process $(X_t)_{t\in\T}$ such that\vspace{-0.1cm}
$$
\as,\ \ \RE\ n,p\geq 1,\ \ \ |Y^{n,p}_\cdot|\leq X_\cdot.\vspace{0.1cm}
$$
Now, we can apply \cref{Pro:2.1--stability} to obtain the existence of a progressively measurable process $(Z_t)_{t\in\T} $ such that $(Y_\cdot:=\inf_p\sup_n Y^{n,p}_\cdot, \ Z_\cdot)$ is a solution to BSDE$(\xi,g)$.

Furthermore, sending $n$ and $p$ to infinity in \eqref{eq:4.21} yields the desired inequality \eqref{eq:4.18}, and then the process $\left(\bar\psi\left(|Y_t|,\bar\mu_{\beta,\gamma,\lambda,\eps}(t)\right)\right)_{t\in\T}$ belongs to class (D). The existence part is then proved.

As for the uniqueness, it is easy to verify that if $\left(\bar\psi\left(|Y_t|,\bar\mu_{\beta,\gamma,\lambda,\eps}(t)\right)\right)_{t\in\T}$ belongs to class (D), then both $\left(\psi\left(|Y_t|,\bar\mu_{\beta,\gamma,\lambda,\eps}(t)\right)
\right)_{t\in\T}$ and $\left(\psi\left(|Y_t|,\bar\mu_{\beta,\gamma,\lambda,\eps}(t)-\eps\right)
\right)_{t\in\T}$ also belong to class (D), where the function $\psi(x,\mu)$ is defined in \eqref{eq:2.12}. Then, the uniqueness of the desired solution under Assumptions \ref{A:UN1} and \ref{A:UN2} follows from \cref{pro:2.5} immediately.

Finally, the uniqueness under Assumption \ref{A:UN3} is a consequence of the following \cref{pro:4.7}. The proof is complete.
\end{proof}

\begin{pro}\label{pro:4.7}
Let the functions $\bar \mu_{\beta,\gamma,\lambda,\eps}(\cdot)$ and $\bar \psi(\cdot,\cdot)$ be respectively defined in \eqref{eq:4.12} and \eqref{eq:4.16}, $\xi$ and $\xi'$ be two terminal conditions, $g$ and $g'$ be two generators which are continuous in the variables $(y,z)$, and $(Y_t, Z_t)_{t\in\T}$ and $(Y'_t, Z'_t)_{t\in\T}$ be respectively a solution to BSDE$(\xi, g)$ and BSDE$(\xi', g')$ such that both $\left(\bar\psi\left(|Y_t|,\ \bar \mu_{\beta,\gamma,\lambda,\eps}(t)\right)\right)_{t\in\T}$ and $\left(\bar\psi\left(|Y'_t|,\ \bar \mu_{\beta,\gamma,\lambda,\eps}(t)\right)\right)_{t\in\T}$
belong to class (D) for some $\eps>0$.

Assume that $\ps$, $\xi\leq \xi'$. If $g$ (resp. $g'$) satisfies Assumptions~\ref{A:UN3} and~\ref{A:A2} with parameters $(\alpha_\cdot,\beta,\gamma,\lambda)$ such that $\bar\psi\left(\int_0^T \alpha_t{\rm d}t, \ \bar \mu_{\beta,\gamma,\lambda,\eps}(T)\right)\in L^1$, and $\as$,
\begin{equation}\label{eq:5.18}
g(t,Y'_t,Z'_t)\leq g'(t,Y'_t,Z'_t)\ \ \ ({\rm resp.}\  \ g(t,Y_t,Z_t)\leq g'(t,Y_t,Z_t)\ ),
\end{equation}
then $\ps$, for each $t\in\T$, $Y_t\leq Y'_t$.
\end{pro}

\begin{proof}
Note that the function $\bar \psi(x,\mu)$ has similar properties to the function $\psi(x,\mu)$ defined in \eqref{eq:2.12}. By a same analysis as in the proof of \cref{pro:3.7} and with the help of \cref{rmk:4.6} and \cref{pro:4.5}, one can derive the desired conclusion, whose proof is omitted here.
\end{proof}

\begin{ex}\label{exp:4.8}
Let $\lambda\in (0,1/2)$, $\beta>0$, $\gamma>0$ and $k\geq 0$. For $(\omega,t,y,z)\in \Omega\times\T\times\R\times\R^d$, define
$$
g(\omega,t,y,z):=|B_t(\omega)|+ky^2{\bf 1}_{y\leq 0}+\beta |y|(\ln |y|)^{\lambda+{1\over 2}}{\bf 1}_{|y|> 1}+\gamma |z|(\ln |z|)^{\lambda}{\bf 1}_{|z|>1}.
$$
It is easy to verify that this generator $g$ satisfies Assumptions \ref{A:EX1}, \ref{A:A2} and \ref{A:UN3}.
\end{ex}

\section{The case of $L^{(1+\mu)}$-integrability terminal condition}
\label{sec:5-mainresult--5}
\setcounter{equation}{0}

In this section, we always assume that the generator $g$ satisfies the following assumption,  which is strictly weaker than \eqref{eq:1.2} with $\lambda=1/2$ and $\delta=1$.
\begin{enumerate}
\renewcommand{\theenumi}{(A3)}
\renewcommand{\labelenumi}{\theenumi}
\item \label{A:A3} There exist two constants $\beta\geq 0$ and $\gamma > 0$, and an $\R_+$-valued progressively measurable process $(\alpha_t)_{t\in\T}$  such that $\as$,
    $$
    \RE\ (y,z)\in \R\times\R^d,\ \ \
    {\rm sgn}(y)g(\omega,t,y,z)\leq \alpha_t(\omega)+\beta|y|\ln |y|{\bf 1}_{|y|>1}+\gamma |z|\sqrt{|\ln |z||}.
    $$
  \end{enumerate}

As in section 4, for $\lambda=1/2$ and $\eps>0$, by letting $k=1+\eps$ in \cref{Pro:2.2Inequality} we know that there exists a constant $C_{\eps}>0$ depending only on $\eps$ such that \eqref{eq:2.5} holds. Then, our objective is to search for a strictly increasing and convex function $\phi$ such that for each $(s,x)\in \T\times \R_+$, with $\beta\geq 0$ and $\gamma>0$,\vspace{0.1cm}
\begin{equation}\label{eq:5.1}
-\beta\phi_x(s, x) x\ln x {\bf 1}_{x>1}-{\gamma^2\over 2}\frac{\left(\phi_x(s, x)\right)^2}
{\phi_{xx}(s, x)}\left[(1+\eps)\left|\ln \frac{\gamma\phi_x(s, x)}{\phi_{xx}(s, x)}
\right|+C_{\eps}\right]+\phi_s(s, x)\geq 0,\vspace{0.1cm}
\end{equation}
which is just \eqref{eq:2.11} with $\lambda=1/2$, $\delta=1$, $\bar\delta=0$, $k=1+\eps$ and $C_{k,\lambda}=C_{\eps}$.\vspace{0.2cm}

Now, let $\eps>0$, $\beta\geq 0$, $\gamma>0$, and $\mu_s,\nu_s:\T\To\R_+$ be two increasing and continuously differential functions with $\mu_0=\eps$ and $\nu_0=1$. We choose the following function
\begin{equation}\label{eq:5.2}
\phi(s,x):=\nu_s (1+x)^{1+\mu_s}, \ \ \ (s,x)\in \T\times \R_+
\end{equation}
to explicitly solve the inequality \eqref{eq:5.1}. For each $(s,x)\in \T\times \R_+$, a simple computation gives
\begin{equation}\label{eq:5.3}
\phi_x(s,x)=\phi(s,x)\frac{1+\mu_s}{1+x}>0,
\end{equation}
\begin{equation}\label{eq:5.4}
\phi_{xx}(s,x)=\phi(s,x)\frac{\mu_s(1+\mu_s)}{(1+x)^2}>0
\end{equation}
and\vspace{-0.1cm}
\begin{equation}\label{eq:5.5}
\phi_s(s,x)=\phi(s,x)\left(\mu'_s\ln (1+x)+\frac{\nu'_s}{\nu_s}\right).\vspace{0.1cm}
\end{equation}
Note that\vspace{-0.1cm}
$$\RE\ x\in\R_+,\ \ x\ln x {\bf 1}_{x>1}\leq (1+x)\ln (1+x).$$
We substitute \eqref{eq:5.3}, \eqref{eq:5.4} and \eqref{eq:5.5} into the left side of \eqref{eq:5.1} to get that for each $(s,x)\in \T\times \R_+$,\vspace{0.1cm}
$$
\begin{array}{ll}
&\Dis -\beta\phi_x(s, x) x\ln x {\bf 1}_{x>1}-{\gamma^2\over 2}\frac{\left(\phi_x(s, x)\right)^2}
{\phi_{xx}(s, x)}\left[(1+\eps)\left|\ln \frac{\gamma\phi_x(s, x)}{\phi_{xx}(s, x)}
\right|+C_{\eps}\right]+\phi_s(s, x)\vspace{0.3cm}\\
\geq & \Dis -\beta \phi(s, x)(1+\mu_s)\ln (1+x)-{\gamma^2(1+\mu_s)\over 2\mu_s}\phi(s, x)\left[(1+\eps)\left|\ln \frac{\gamma(1+x)}{\mu_s}\right|+C_{\eps}\right]\\
& +\phi(s,x)\left(\mu'_s\ln (1+x)+\frac{\nu'_s}{\nu_s}\right)\vspace{0.3cm}\\
\geq & \Dis \phi(s, x) \left[\left(-\beta(1+\mu_s)-{\gamma^2(1+\mu_s)(1+\eps)\over 2\mu_s}+\mu'_s\right)\ln (1+x)+\left(-\gamma^2(1+\mu_s)\bar C_\eps+\frac{\nu'_s}{\nu_s}\right)\right]\vspace{0.1cm}
\end{array}
$$
where
$$
\bar C_{\eps}:=\frac{(1+\eps)|\ln \gamma-\ln \eps|+C_{\eps}}{2\eps}.\vspace{0.1cm}
$$
Thus, \eqref{eq:5.1} holds if the functions $\mu_s, \nu_s\in \T$ satisfies\vspace{0.2cm}
\begin{equation}\label{eq:5.6}
\mu'_s=\beta\mu_s+{\gamma^2(1+\eps)\over 2\mu_s}+{\gamma^2(1+\eps)\over 2}+\beta,\ \ \ s\in [0,T]
\end{equation}
and
\begin{equation}\label{eq:5.7}
\nu_s=\exp\left(\gamma^2 \bar C_\eps \int_0^s(1+\mu_r){\rm d}r \right) ,\ \ \ s\in [0,T].\vspace{0.3cm}
\end{equation}

It is not very hard to verify that for each $\beta\geq 0$, $\gamma > 0$ and $\eps>0$, there exists a unique strictly increasing and continuous function with $\eps$ at the origin satisfying \eqref{eq:5.6}. We denote this unique function by $\tilde  \mu_{\beta,\gamma,\eps}(\cdot)$, i.e., $\tilde  \mu_{\beta,\gamma,\eps}(0)=\eps$ and
\begin{equation}\label{eq:5.8}
\tilde  \mu'_{\beta,\gamma,\eps}(s)=\beta
\tilde \mu_{\beta,\gamma,\eps}(s)
+\frac{\gamma^2(1+\eps)}{2}
\frac{1}{\tilde\mu_{\beta,\gamma,\eps}(s)}+\frac{\gamma^2(1+\eps)}{2}+\beta,\ \ \ s\in [0,T].
\end{equation}
We also denote, in view of \eqref{eq:5.7},
\begin{equation}\label{eq:5.9}
\tilde  \nu_{\beta,\gamma,\eps}(s):=\exp\left(\gamma^2 \bar C_\eps \int_0^s(1+\tilde\mu_{\beta,\gamma,\eps}(r)){\rm d}r\right),\ \ \ s\in [0,T].\vspace{0.1cm}
\end{equation}

Furthermore, it can be directly checked that as $\eps\To 0$, the function $\tilde  \mu_{\beta,\gamma,\eps}(\cdot)$ tends decreasingly to the unique solution $\tilde  \mu_{\beta,\gamma}^0(\cdot)$ of the following ODE:
\begin{equation}\label{eq:5.10}
\tilde\mu_{\beta,\gamma}^0(0)=0\ \ \ {\rm and}\ \ \ \left(\tilde  \mu_{\beta,\gamma}^0(s)\right)'=\beta\tilde\mu_{\beta,\gamma}^0(s)+
\frac{\gamma^2}{2}
\frac{1}{\tilde\mu_{\beta,\gamma}^0(s)}+\frac{\gamma^2}{2}+\beta,\ \ \ s\in (0,T].\vspace{0.2cm}
\end{equation}

\begin{rmk}\label{rmk:5.1}
In general, $\tilde  \mu_{\beta,\gamma}^0(\cdot)$ does not have an explicit expression. \vspace{0.2cm}
\end{rmk}

In summary, we have proved the following proposition.

\begin{pro}\label{pro:5.2}
For each $\beta\geq 0$, $\gamma > 0$ and $\eps>0$, define the function
\begin{equation}\label{eq:5.11}
\tilde  \varphi(s,x;\eps):=\tilde\nu_{\beta,\gamma,\eps}(s) (1+x)^{1+\tilde\mu_{\beta,\gamma,\eps}(s)}, \ \ \ (s,x)\in [0,T]\times \R_+
\end{equation}
with $\tilde\mu_{\beta,\gamma,\eps}(\cdot)$ and $\tilde\nu_{\beta,\gamma,\eps}(\cdot)$ being respectively defined in \eqref{eq:5.8} and \eqref{eq:5.9}. Then we have
\begin{itemize}
\item [(i)] $\tilde\varphi(\cdot,\cdot;\eps)\in \mathcal {C}^{1,2}\left([0,T]\times \R_+\right)$;
\item [(ii)] $\tilde  \varphi(\cdot,\cdot;\eps)$ satisfies the inequality in \eqref{eq:2.3} with $\lambda=1/2$, $\delta=1$ and $\bar\delta=0$. i.e., for each $(s,x,z)\in [0,T]\times \R_+\times\R^d$, we have\vspace{-0.1cm}
$$
-\tilde  \varphi_x(s,x;\eps)\left(\beta x\ln x{\bf 1}_{x>1}+\gamma |z|\sqrt{|\ln |z||}\right)+{1\over 2}\tilde  \varphi_{xx}(s,x;\eps)|z|^2+\tilde  \varphi_s(s,x;\eps)\geq 0.
$$
\end{itemize}
\end{pro}

The following existence and uniqueness theorem is one of the main results of this section.

\begin{thm}\label{thm:5.3}
Let the functions $\tilde  \mu_{\beta,\gamma,\eps}(\cdot)$ and $\tilde  \mu_{\beta,\gamma}^0(\cdot)$ be respectively defined in \eqref{eq:5.8} and \eqref{eq:5.10}, $\xi$ be a terminal condition and $g$ be a generator which is continuous in the state variables $(y,z)$. If $g$ satisfies Assumptions \ref{A:EX1} and \ref{A:A3} with parameters $(\alpha_\cdot, \beta, \gamma)$, and there exists a positive constant $\mu>\tilde \mu_{\beta,\gamma}^0(T)$ such that\vspace{-0.1cm}
\begin{equation}\label{eq:5.12}
\E\left[\left(|\xi|+\int_0^T \alpha_t {\rm d}t\right)^{1+\mu}\right]<+\infty,\vspace{0.1cm}
\end{equation}
then BSDE $(\xi,g)$ admits a solution $(Y_t,Z_t)_{t\in\T}$ such that for some $\eps>0$, $\left(|Y_t|^{1+\tilde\mu_{\beta,\gamma,\eps}(t)}\right)_{t\in\T}$ belongs to class (D), $Z_\cdot\in \mcal^p$ for some $p>1$ and $\ps$,
\begin{equation}\label{eq:5.13}
\begin{array}{lll}
\Dis |Y_t|&\leq &\Dis \left(1+|Y_t|\right)^{1+\tilde\mu_{\beta,\gamma,\eps}(t)}\leq C\E\left[\left.\left(1+|\xi|+\int_0^T\alpha_t {\rm d}t\right)^{1+\tilde\mu_{\beta,\gamma,\eps}(T)}\right|\F_t\right]\vspace{0.2cm}\\
&\leq &\Dis 2^\mu C\left(\E\left[\left.\left(|\xi|+\int_0^T\alpha_t {\rm d}t\right)^{1+\mu}\right|\F_t\right]+1\right),\ \ \ \ t\in\T,
\end{array}
\end{equation}
where $C$ is a positive constant depending only on $(\beta,\gamma,\eps,T)$, and $\eps$ is the unique positive constant satisfying $\tilde\mu_{\beta,\gamma,\eps}(T)=\mu$.

Moreover, assume that $g=g_1+g_2$. If $g_1$ satisfies Assumptions \ref{A:UN1} and \ref{A:UN2} and $g_2$ satisfies Assumption \ref{A:UN3}, then BSDE $(\xi,g)$ admits a unique solution $(Y_t,Z_t)_{t\in\T}$ such that the process $\left(|Y_t|^{1+\tilde\mu_{\beta,\gamma,\eps}(t)}\right)_{t\in\T}$ belongs to class (D) and $Z_\cdot\in \mcal^p$ for some $p>1$.
\end{thm}

In order to prove this theorem, we need the following proposition.

\begin{pro}\label{pro:5.4}
Let the function $\tilde \mu_{\beta,\gamma,\eps}(\cdot)$ be defined in \eqref{eq:5.8}, $\xi$ be a terminal condition and $g$ be a generator which is continuous in $(y,z)$. If $g$ satisfies Assumption \ref{A:A3} with parameters $(\alpha_\cdot, \beta, \gamma)$, $|\xi|+\int_0^T \alpha_t {\rm d}t$ be a bounded random variable, and $(Y_t,Z_t)_{t\in\T}$ is a solution of BSDE $(\xi,g)$ such that $Y_\cdot$ is bounded, then for each $\eps>0$, $\ps$,
\begin{equation}\label{eq:5.14}
|Y_t|\leq \left(1+|Y_t|\right)^{1+\tilde\mu_{\beta,\gamma,\eps}(t)}\leq C\E\left[\left.\left(1+|\xi|+\int_0^T\alpha_t {\rm d}t\right)^{1+\tilde\mu_{\beta,\gamma,\eps}(T)}\right|\F_t\right],\ \ \ t\in \T,
\end{equation}
where $C$ is a positive constant depending only on $(\beta,\gamma,\eps,T)$.
\end{pro}

\begin{proof}
Define
$$
\tilde  Y_t:=|Y_t|+\int_0^t \alpha_s {\rm d}s\ \ \ \
{\rm and}\ \ \ \ \tilde  Z_t:={\rm sgn}(Y_t)Z_t,\ \ \ \ t\in \T.
$$
It follows from It\^{o}-Tanaka's formula that
$$
\tilde  Y_t=\tilde  Y_T+\int_t^T \left({\rm sgn}(Y_s)g(s,Y_s,Z_s)-\alpha_s\right){\rm d}s-\int_t^T \tilde  Z_s \cdot {\rm d}B_s-\int_t^T {\rm d}L_s, \ \ \ t\in\T,
$$
where $L_\cdot$ denotes the local time of $Y_\cdot$ at the origin. Now, we fix $\eps>0$ and apply It\^{o}-Tanaka's formula to the process $\tilde \varphi(s, \tilde  Y_s;\eps)$, where the function $\tilde \varphi(\cdot,\cdot;\eps)$ is defined in \eqref{eq:5.11}, to derive, in view of \ref{A:A3},
$$
\begin{array}{lll}
\Dis {\rm d}\tilde \varphi(s,\tilde  Y_s;\eps)
&=&\Dis \tilde \varphi_x(s,\tilde  Y_s;\eps)
\left(-{\rm sgn}(Y_s)g(s,Y_s,Z_s)+\alpha_s\right){\rm d}s+\tilde \varphi_x(s,\tilde  Y_s;\eps)\tilde  Z_s \cdot {\rm d}B_s+\tilde \varphi_x(s,\tilde  Y_s;\eps){\rm d}L_s\vspace{0.1cm}\\
&&\Dis +{1\over 2}\tilde \varphi_{xx}(s,\tilde  Y_s;\eps)|Z_s|^2{\rm d}s+\tilde \varphi_s(s,\tilde  Y_s;\eps){\rm d}s\vspace{0.2cm}\\
&\geq &\Dis \bigg[-\tilde \varphi_x(s,\tilde  Y_s;\eps)\left(\beta |Y_s|\ln |Y_s|{\bf 1}_{|Y_s|>1}+\gamma |Z_s|\sqrt{|\ln |Z_s||}\right)\vspace{0.1cm}\\
&&\Dis
\ \ +{1\over 2}\tilde \varphi_{xx}(s,\tilde  Y_s;\eps)|Z_s|^2+\varphi_s(s,\tilde  Y_s;\eps)\bigg]{\rm d}s+\tilde \varphi_x(s,\tilde  Y_s;\eps)\tilde  Z_s \cdot {\rm d}B_s,\ \ s\in\T.
\end{array}
$$
Furthermore, in view of the fact $|Y_s|\ln |Y_s|{\bf 1}_{|Y_s|>1}\leq \tilde  Y_s\ln \tilde  Y_s {\bf 1}_{\tilde  Y_s>1}$ and~\cref{pro:5.2},  we have \vspace{-0.1cm}
\begin{equation}\label{eq:5.15}
{\rm d}\tilde \varphi(s,\tilde  Y_s;\eps)\geq \tilde \varphi_x(s,\tilde  Y_s;\eps)\tilde  Z_s \cdot {\rm d}B_s,\ \ s\in \T.\vspace{-0.1cm}
\end{equation}
Let us consider, for each integer $n\geq 1$ and each $t\in \T$,  the following stopping time
$$
\tau_n^t:=\inf\left\{s\in [t,T]: \int_t^s \left[\tilde  \varphi_x(r,\tilde  Y_r;\eps)\right]^2|\tilde  Z_r|^2{\rm d}r\geq n \right\}\wedge T.
$$
In view of inequality~\eqref{eq:5.15}, we have
$$
\tilde \varphi(t,\tilde  Y_t;\eps)\leq \E\left[\left. \tilde \varphi(\tau_n^t,\tilde  Y_{\tau_n^t};\eps) \right|\F_t\right],\quad  t\in \T, \quad n\geq 1.\vspace{-0.2cm}
$$
Thus, we have\vspace{-0.1cm}
$$
\begin{array}{lll}
\Dis |Y_t|&\leq & \Dis \left(1+|Y_t|\right)^{1+\tilde\mu_{\beta,\gamma,\eps}(t)}\leq \tilde\nu_{\beta,\gamma,\eps}(t)\left(1+\tilde Y_t\right)^{1+\tilde\mu_{\beta,\gamma,\eps}(t)}\vspace{0.2cm}\\
&\leq& \Dis \tilde\nu_{\beta,\gamma,\eps}(T)\E\left[\left.\left(1+|Y_{\tau_n^t}|
+\int_0^{\tau_n^t}\alpha_t {\rm d}t\right)^{1+\tilde\mu_{\beta,\gamma,\eps}(T)}\right|\F_t\right],\ \ \ t\in \T,
\end{array}
$$
which yields the desired inequality~\eqref{eq:5.14} in the limit as  $n\to \infty$. The proof is completed.
\end{proof}

\begin{rmk}\label{rmk:5.5} Assertions of \cref{pro:5.4} are still  true if in \cref{pro:5.4},  $(|\xi|, |Y_t|)$ is replaced with $(\xi^+, Y^+_t)$, and Assumption~\ref{A:A3} is replaced with the following one~\ref{A:A3'}:
\begin{enumerate}
\renewcommand{\theenumi}{(A3')}
\renewcommand{\labelenumi}{\theenumi}
\item\label{A:A3'} There exist two constants $\beta\geq 0$ and $\gamma > 0$, and an $\R_+$-valued progressively measurable process $(\alpha_t)_{t\in\T}$  such that $\as$,
    \[
    \RE\ (y,z)\in \R_+\times\R^d,\ \ \
    g(\omega,t,y,z)\leq \alpha_t(\omega)+\beta |y| \ln |y| {\bf 1}_{|y|>1}+\gamma |z|\sqrt{|\ln |z||}.\vspace{-0.2cm}
    \]
\end{enumerate}
In fact, it is sufficient to use $(Y_\cdot^+, {\bf 1}_{Y_\cdot>0} Y_\cdot, \frac{1}{2}L_\cdot)$ instead of $(|Y_\cdot|, {\rm sgn} (Y_\cdot),  L_\cdot)$ in the  proof.\vspace{0.1cm}
\end{rmk}

Now, we prove~\cref{thm:5.3}.\vspace{-0.1cm}

\begin{proof}[Proof of \cref{thm:5.3}]
For each pair of positive integers $n,p\geq 1$, let $\xi^{n,p}$ and $g^{n,p}$ be defined in \eqref{eq:2.1}, and $(Y^{n,p}_\cdot,Z^{n,p}_\cdot)$ be the minimal (maximal) bounded solution of \eqref{eq:2.2}. It is easy to check that the generator $g^{n,p}$ satisfies assumption \ref{A:A3} with $\alpha_\cdot$ being replaced with $\alpha_\cdot \wedge (n\vee p)$.

Now, we assume that there exists a positive constant $\mu>\tilde  \mu_{\beta,\gamma}^0(T)$ such that \eqref{eq:5.12} holds. From the definitions of $\tilde\mu_{\beta,\gamma,\eps}(\cdot)$ and $\tilde \mu_{\beta,\gamma}^0(\cdot)$ in \eqref{eq:5.8} and \eqref{eq:5.10}, it is not very difficult to find a (unique) positive constant $\eps>0$ satisfying $\tilde\mu_{\beta,\gamma,\eps}(T)= \mu$. Then, applying \cref{pro:5.4} with this $\eps$ to BSDE \eqref{eq:2.2} yields the existence of a constant $C>0$ depending only on $(\beta,\gamma,\eps,T)$ such that for all $n,p\geq 1$ and $t\in\T$,\vspace{0.1cm}
\begin{equation}\label{eq:5.16}
\hspace*{-0.2cm}\begin{array}{lll}
|Y^{n,p}_t|&\leq & \Dis \left(1+|Y^{n,p}_t|\right)^{1+\tilde  \mu_{\beta,\gamma,\eps}(t)}\leq C\E\left[\left.\left(1+|\xi^{n,p}|+\int_0^T[\alpha_t\wedge (n\vee p)]{\rm d}t\right)^{1+\tilde\mu_{\beta,\gamma,\eps}(T)}\right|\F_t\right]\vspace{0.2cm}\\
&\leq& \Dis C\E\left[\left.\left(1+|\xi|+\int_0^T\alpha_t {\rm d}t\right)^{1+\tilde\mu_{\beta,\gamma,\eps}(T)}\right|\F_t\right]= C\E\left[\left.\left(1+|\xi|+\int_0^T\alpha_t {\rm d}t\right)^{1+\mu}\right|\F_t\right].\vspace{0.1cm}
\end{array}
\end{equation}
Thus, in view of \eqref{eq:5.12}, we have found an $\R_+$-valued, progressively measurable and continuous process $(X_t)_{t\in\T}$ such that\vspace{-0.2cm}
$$
\as,\ \ \RE\ n,p\geq 1,\ \ \ |Y^{n,p}_t|\leq  X_t\ :=\ C\E\left[\left.\left(1+|\xi|+\int_0^T\alpha_t {\rm d}t\right)^{1+\mu}\right|\F_t\right].
$$
Now, we see from \cref{Pro:2.1--stability} that there is a progressively measurable process $(Z_t)_{t\in\T} $ such that $(Y_\cdot:=\inf_p\sup_n Y^{n,p}_\cdot, \ Z_\cdot)$ is a solution to BSDE $(\xi, g)$.

Furthermore, sending $n$ and $p$ to infinity in \eqref{eq:5.16} yields the desired inequality \eqref{eq:5.13}, and then the process $\left(|Y_t|^{1+\tilde\mu_{\beta,\gamma,\eps}(t)}\right)_{t\in\T}$ belongs to class (D). Thus, to complete the proof of the existence part of \cref{thm:5.3}, it remains to show that $Z_\cdot\in \mcal^p$ for some $p>1$. In view of the inequality  $\tilde\mu_{\beta,\gamma,\eps}(\cdot)\geq \eps$, we see from~\cite[Lemma 6.1]{BriandDelyonHu2003SPA} and \eqref{eq:5.13} that
\[
\begin{array}{lll}
\Dis \E\left[\sup\limits_{t\in\T}|Y_t|^{1+{\eps\over 2}}\right]&\leq &\Dis \E\left[\sup\limits_{t\in\T}\left[\left(1+
|Y_t|\right)^{1+\tilde\mu_{\beta,\gamma,\eps}(t)}\right]^{\frac{1+{\eps\over 2}}{1+\eps}}\right]\leq \E\left[\left(\sup\limits_{t\in\T} X_t\right)^{\frac{1+{\eps\over 2}}{1+\eps}}\right]\vspace{0.2cm}\\
&\leq & \Dis \frac{2(1+\eps)}{\eps} \left(\E\left[X_T\right]\right)^{\frac{1+{\eps\over 2}}{1+\eps}}<+\infty,
\end{array}
\]
which implies that $Y_\cdot\in \s^{1+{\eps\over 2}}$.

Now, we take $p:=\frac{1+\eps/2}{1+\eps/4}>1$ to prove $Z_\cdot\in \mcal^p$. First, by letting $\lambda=1/2$, $k=1$ and $y/2$ instead of the variable $y$ in \cref{Pro:2.2Inequality}, and using both inequalities
$$\sqrt{|\ln y|}\leq \sqrt{|\ln y/2|}+\sqrt{\ln 2}$$ and
$$2xy\sqrt{\ln 2}\leq (4\ln 2)x^2+y^2/4,$$
 we  see that there is a universal constant $K>0$ such that
\begin{equation}\label{eq:5.17}
\RE\ x,y>0,\quad  2xy\sqrt{|\ln y|}\leq 2x^2(|\ln x|+K)+\frac{3}{4}y^2.
\end{equation}
Then, we pick a sufficiently large constant $k_\eps>1$ depending only on $\eps$ such that
\begin{equation}\label{eq:5.18}
\RE\ x\in \R_+,\ \ \ \ln (x+k_\eps)\leq (x+k_\eps)^{\frac{\eps}{2}}.\vspace{-0.1cm}
\end{equation}

For each $n\geq 1$, define the following stopping time\vspace{0.1cm}
$$
\tau_n:=\inf\left\{s\in [0,T]: \int_0^s (|Y_r|+k_\eps)^2|Z_r|^2{\rm d}r\geq n \right\}\wedge T.\vspace{0.1cm}
$$
Applying It\^{o}-Tanaka's formula to $(|Y_s|+k_\eps)^2$, we have
\begin{equation}\label{eq:5.19}
\begin{array}{ll}
&\Dis (|Y_0|+k_\eps)^2+\int_0^{\tau_n}|Z_s|^2{\rm d}s \vspace{0.2cm}\\
\leq & \Dis (|Y_{\tau_n}|+k_\eps)^2 +2\int_0^{\tau_n}(|Y_s|+k_\eps){\rm sgn}(Y_s)g(s,Y_s,Z_s){\rm d}s+2\int_0^{\tau_n}(|Y_s|+k_\eps){\rm sgn}(Y_s)Z_s\cdot {\rm d}B_s.
\end{array}
\end{equation}
It follows from assumption \ref{A:A3} together with \eqref{eq:5.17} and \eqref{eq:5.18} that
$$
\begin{array}{ll}
&\Dis 2(|Y_s|+k_\eps){\rm sgn}(Y_s)g(s,Y_s,Z_s)\vspace{0.2cm}\\
\leq &\Dis 2\alpha_t(|Y_s|+k_\eps)+2\beta (|Y_s|+k_\eps)^2\ln (|Y_s|+k_\eps)+2\gamma (|Y_s|+k_\eps)|Z_s|\sqrt{|\ln |Z_s||}\vspace{0.2cm}\\
\leq &\Dis 2\alpha_t(|Y_s|+k_\eps)+2(\beta+\gamma^2) (|Y_s|+k_\eps)^{2+{\eps\over 2}}+2\gamma^2 K (|Y_s|+k_\eps)^2+{3\over 4}|Z_s|^2.
\end{array}
$$
Then, coming back to \eqref{eq:5.19}, we have for each $n\geq 1$,
$$
\begin{array}{lll}
\Dis \int_0^{\tau_n}|Z_s|^2{\rm d}s
&\leq & \Dis 4\left(1+2\gamma^2 K\right)\sup\limits_{s\in\T}(|Y_s|+k_\eps)^2 +8\sup\limits_{s\in\T}(|Y_s|+k_\eps)\int_0^T \alpha_s{\rm d}s\vspace{0.2cm}\\
&& \Dis +8(\beta+\gamma^2) \sup\limits_{s\in\T}(|Y_s|+k_\eps)^{2+\frac{\eps}{2}}
+8\int_0^{\tau_n}(|Y_s|+k_\eps){\rm sgn}(Y_s)Z_s\cdot {\rm d}B_s,
\end{array}
$$
and then there exists a positive constant $c_{p}>0$ depending only on $(p,\beta,\gamma,K)$ such that\vspace{-0.2cm}
\begin{equation}\label{eq:5.20}
\begin{array}{ll}
&\Dis \E\left[\left(\int_0^{\tau_n}|Z_s|^2{\rm d}s\right)^{{p\over 2}}\right] \\
\leq & \Dis c_p\E\left[\sup\limits_{s\in\T}(|Y_s|+k_\eps)^p\right] + c_p\E\left[\left(\sup\limits_{s\in\T}(|Y_s|+k_\eps)\int_0^T\alpha_s{\rm d}s\right)^{{p\over 2}}\right]\vspace{0.2cm}\\
& \Dis +c_p \E\left[\sup\limits_{s\in\T}(|Y_s|+k_\eps)^{1+{\eps\over 2}}\right] +c_p \E\left[\left(\int_0^{\tau_n}(|Y_s|+k_\eps){\rm sgn}(Y_s)Z_s\cdot {\rm d}B_s\right)^{{p\over 2}}\right].
\end{array}
\end{equation}

Moreover, by the BDG inequality and Young's inequality, we can find another constant $d_p>0$ depending only on $(p,\beta,\gamma,K)$ such that
$$
c_p \E\left[\left(\int_0^{\tau_n}(|Y_s|+k_\eps){\rm sgn}(Y_s)Z_s\cdot {\rm d}B_s\right)^{{p\over 2}}\right]\leq d_p \E\left[\sup\limits_{s\in\T}(|Y_s|+k_\eps)^p\right]+{1\over 2}\E\left[\left(\int_0^{\tau_n}|Z_s|^2{\rm d}s\right)^{{p\over 2}}\right].
$$
Substituting the above inequality to \eqref{eq:5.20} and using H\"{o}lder's inequality imply that for each $n\geq 1$,
$$
\begin{array}{lll}
\Dis {1\over 2}\E\left[\left(\int_0^{\tau_n}|Z_s|^2{\rm d}s\right)^{{p\over 2}}\right]&\leq & \Dis (c_p+d_p)\E\left[\sup\limits_{s\in\T}(|Y_s|+k_\eps)^p\right] +c_p \E\left[\sup\limits_{s\in\T}(|Y_s|+k_\eps)^{1+{\eps\over 2}}\right] \vspace{0.2cm}\\
&& \Dis +c_p\left(\E\left[\sup\limits_{s\in\T}(|Y_s|+k_\eps)^p\right]\right)^{{1\over 2}}\left(\E\left[\left(\int_0^T\alpha_s{\rm d}s\right)^p\right]\right)^{{1\over 2}}.
\end{array}
$$
Thus, in view of \eqref{eq:5.12} and the facts that $p<1+{\eps\over 2}<1+\mu$ and $Y_\cdot\in \s^{1+{\eps\over 2}}$, using Fatou's Lemma, we have $Z_\cdot\in \mcal^p$. The existence part is then proved.

Finally, the uniqueness part of \cref{thm:5.3} is a direct consequence of the following \cref{pro:5.6}. The proof is complete.
\end{proof}

\begin{pro}\label{pro:5.6}
Let the function $\tilde  \mu_{\beta,\gamma,\eps}(\cdot)$  be respectively defined in \eqref{eq:5.8}, $\xi$ and $\xi'$ be two terminal conditions, $g$ and $g'$ be two generators which are continuous in the variables $(y,z)$, and $(Y_t, Z_t)_{t\in\T}$ and $(Y'_t, Z'_t)_{t\in\T}$ be respectively a solution to BSDE $(\xi, g)$ and BSDE $(\xi', g')$ satisfying both $\left(|Y_t|^{1+\tilde\mu_{\beta,\gamma,\eps_0}(t)}\right)_{t\in\T}$ and $\left(|Y'_t|^{1+\tilde\mu_{\beta,\gamma,\eps_0}(t)}\right)_{t\in\T}$
belong to class (D) for some $\eps_0>0$ and there exists an $X\in L^1$ such that
\begin{equation}\label{eq:5.21}
(1+|Y_t|+|Y'_t|)^{1+\tilde\mu_{\beta,\gamma,\eps_0}(t)}
\leq \E[X|\F_t],\ \ \ t\in \T.
\end{equation}

Assume that $\ps$, $\xi\leq \xi'$. If $g=g_1+g_2$ (resp. $g'=g'_1+g'_2$), $g_1$ (resp. $g'_1$) satisfies Assumptions \ref{A:UN1} and \ref{A:UN2}, $g_2$ (resp. $g'_2$) satisfies Assumptions \ref{A:UN3} and \ref{A:A3} with parameters $(\alpha_\cdot,\beta,\gamma)$ such that $\ps$, $\int_0^T \alpha_t{\rm d}t<+\infty$, and $\as$,
$$
g(t,Y'_t,Z'_t)\leq g'(t,Y'_t,Z'_t)\ \ \ ({\rm resp.}\  \ g(t,Y_t,Z_t)\leq g'(t,Y_t,Z_t)\ ),
$$
then $\ps$, for each $t\in\T$, $Y_t\leq Y'_t$.
\end{pro}

\begin{proof}
We only prove the case that the generator $g=g_1+g_2$, $g_1$ satisfies \ref{A:UN1} and \ref{A:UN2}, $g_2$ is convex in the state variables $(y,z)$ and satisfies \ref{A:A3} with parameters $(\alpha_\cdot,\beta,\gamma)$, and $\as$,
$$g(t,Y'_t,Z'_t)\leq g'(t,Y'_t,Z'_t).\vspace{-0.1cm}$$
The other cases can be proved in the same way.

We first use the $\theta$-technique. For each fixed $\theta\in (0,1)$, define
\begin{equation}\label{eq:5.22}
\Delta^\theta U_\cdot:=\frac{Y_\cdot-\theta Y'_\cdot}{1-\theta}\ \  {\rm and} \ \ \Delta^\theta V_\cdot:=\frac{Z_\cdot-\theta Z'_\cdot}{1-\theta}.
\end{equation}
Then the pair $(\Delta^\theta U_\cdot,\Delta^\theta V_\cdot)$ satisfies the following BSDE:\vspace{0.1cm}
\begin{equation}\label{eq:5.23}
  \Delta^\theta U_t=\Delta^\theta U_T +\int_t^T \frac{g(s,Y_s, Z_s)-\theta g'(s,Y'_s, Z'_s)}{1-\theta}{\rm d}s-\int_t^T \Delta^\theta V_s \cdot {\rm d}B_s, \ \ \ \ t\in\T.\vspace{0.1cm}
\end{equation}
Denote $\Delta^\theta g_1(s):=g_1(s,Y_s,Z_s)-\theta g_1(s, Y'_s, Z'_s)$. Observe from the assumptions that $\ass$,
$$
\begin{array}{lll}
&&\Dis {\bf 1}_{\Delta^\theta U_s>0}(g(s,Y_s, Z_s)-\theta g'(s,Y'_s, Z'_s))\vspace{0.2cm}\\
&=& \Dis {\bf 1}_{\Delta^\theta U_s>0}[g(s,Y_s, Z_s)-\theta g(s,Y'_s, Z'_s)+\theta (g(s,Y'_s, Z'_s)-g'(s,Y'_s, Z'_s))]\vspace{0.2cm}\\
&\leq & \Dis {\bf 1}_{\Delta^\theta U_s>0}\left\{\Delta^\theta g_1(s)+\left[g_2\left(s,\theta Y'_s+(1-\theta)\Delta^\theta U_s,\theta Z'_s+(1-\theta)\Delta^\theta V_s\right)-\theta g_2(s, Y'_s, Z'_s)\right]\right\}\vspace{0.2cm}\\
&\leq & \Dis {\bf 1}_{\Delta^\theta U_s>0}\left[\Delta^\theta g_1(s) +(1-\theta)g_2\left(s,\Delta^\theta U_s,\Delta^\theta V_s\right)\right]\vspace{0.2cm}\\
&\leq & \Dis {\bf 1}_{\Delta^\theta U_s>0}\Delta^\theta g_1(s) +(1-\theta)\left(\alpha_s+\beta \left(\Delta^\theta U_s\right)^+\ln (\left(\Delta^\theta U_s\right)^+){\bf 1}_{\left(\Delta^\theta U_s\right)^+>1}+\gamma |\Delta^\theta V_s|\sqrt{\left|\ln |\Delta^\theta V_s|\right|}\right)
\end{array}
$$
and then
\begin{equation}\label{eq:5.24}
\hspace*{-0.4cm}\begin{array}{lll}
&&\Dis {\bf 1}_{\Delta^\theta U_s>0}\frac{g(s,Y_s, Z_s)-\theta g'(s,Y'_s, Z'_s)}{1-\theta}\vspace{0.2cm}\\
&\leq & \Dis \frac{{\bf 1}_{\Delta^\theta U_s>0}\Delta^\theta g_1(s)}{1-\theta}+
\alpha_s+\beta \left(\Delta^\theta U_s\right)^+\ln (\left(\Delta^\theta U_s\right)^+){\bf 1}_{\left(\Delta^\theta U_s\right)^+>1}+\gamma |\Delta^\theta V_s|\sqrt{\left|\ln |\Delta^\theta V_s|\right|}.
\end{array}
\end{equation}

Now, we pick $\eps=\eps_0/2$. Apply It\^{o}-Tanaka's formula to the process $\tilde \varphi(s, \left(\Delta^\theta U_s\right)^+;\eps)$, where the function $\tilde \varphi(\cdot,\cdot;\eps)$ is defined in \eqref{eq:5.11}, and in view of \eqref{eq:5.23}, \eqref{eq:5.24} and \cref{pro:5.2}, we use a similar argument to that in the proof of \cref{pro:5.4} to derive that for each $s\in\T$,
\begin{equation}\label{eq:5.25}
\begin{array}{lll}
\Dis {\rm d}\tilde\varphi(s,\left(\Delta^\theta U_s\right)^+;\eps)
&\geq & \Dis -\tilde\varphi_x(s,\left(\Delta^\theta U_s\right)^+;\eps)\left(\frac{{\bf 1}_{\Delta^\theta U_s>0}\Delta^\theta g_1(s)}{1-\theta}+
\alpha_s(\omega)\right){\rm d}s\vspace{0.2cm}\\
&& \Dis
+\tilde\varphi_x(s,\left(\Delta^\theta U_s\right)^+;\eps){\bf 1}_{\Delta^\theta U_s>0}\Delta^\theta V_s \cdot {\rm d}B_s.\vspace{0.1cm}
\end{array}
\end{equation}
For each $t\in \T$ and each integer $n\geq 1$, define the following stopping time
$$
\tau_n^t:=\inf\left\{s\in [t,T]: |Y_s|+|Y'_s|+\int_t^s (|Z_r|^2+|Z'_r|^2){\rm d}r\geq n \right\}\wedge T.
$$
By virtue of \eqref{eq:5.25} and the definition of $\tilde\varphi$, we deduce that for each $t\in\T$, \vspace{0.1cm}
$$
\begin{array}{lll}
\Dis \left(\left(\Delta^\theta U_t\right)^+\right)^{1+\tilde\mu_{\beta,\gamma,\eps}(t)} &\leq & \Dis \tilde\nu_{\beta,\gamma,\eps}(t)\left(1+\left(\Delta^\theta U_t\right)^+\right)^{1+\tilde\mu_{\beta,\gamma,\eps}(t)}\vspace{0.2cm}\\
&\leq & \Dis \tilde\nu_{\beta,\gamma,\eps}(T)\left(1+\left(\Delta^\theta U_{\tau_n^t}\right)^+\right)^{1+\tilde\mu_{\beta,\gamma,\eps}(\tau_n^t)}\vspace{0.2cm}\\
&&\Dis +\int_t^{\tau_n^t} \tilde\nu_{\beta,\gamma,\eps}(s)\left(1+\tilde\mu_{\beta,\gamma,\eps}(s)\right)
\left(1+\left(\Delta^\theta U_s\right)^+\right)^{\tilde\mu_{\beta,\gamma,\eps}(s)}
\left(\frac{\Delta^\theta g_1(s)}{1-\theta}+
\alpha_s(\omega)\right){\rm d}s\vspace{0.2cm}\\
&&\Dis -\int_t^{\tau_n^t} \tilde\nu_{\beta,\gamma,\eps}(s)\left(1+\tilde\mu_{\beta,\gamma,\eps}(s)\right)
\left(1+\left(\Delta^\theta U_s\right)^+\right)^{\tilde\mu_{\beta,\gamma,\eps}(s)}\Delta^\theta V_s \cdot {\rm d}B_s,\vspace{0.2cm}
\end{array}
$$
and then in view of \eqref{eq:5.22}, for each $\theta\in (0,1)$ and $n\geq 1$,\vspace{-0.1cm}
\begin{equation}\label{eq:5.26}
\begin{array}{ll}
&\Dis \left((Y_t-\theta Y'_t)^+\right)^{1+\tilde\mu_{\beta,\gamma,\eps}(t)}  \vspace{0.2cm}\\
\leq & \Dis K_T\left((1-\theta)+(Y_{\tau_n^t}-\theta Y'_{\tau_n^t})^+\right)^{1+\tilde\mu_{\beta,\gamma,\eps}(\tau_n^t)} \vspace{0.2cm}\\
& \Dis+\int_t^{\tau_n^t}
K_s\left((1-\theta)+(Y_s-\theta Y'_s)^+\right)^{\tilde\mu_{\beta,\gamma,\eps}(s)}\left[\Delta^\theta g_1(s)+
(1-\theta)\alpha_s(\omega)\right]{\rm d}s\vspace{0.2cm}\\
& \Dis -\int_t^{\tau_n^t}
K_s\left((1-\theta)+(Y_s-\theta Y'_s)^+\right)^{\tilde\mu_{\beta,\gamma,\eps}(s)}(Z_s-\theta Z'_s) \cdot {\rm d}B_s, \ \ \ t\in \T,
\end{array}
\end{equation}
where $K_s:=\tilde\nu_{\beta,\gamma,\eps}(s)\left(1+\tilde\mu_{\beta,\gamma,\eps}(s)\right),\ s\in\T$.

Moreover, define $\hat U_\cdot:=Y_\cdot-Y'_\cdot$ and $\hat V_\cdot:=Z_\cdot-Z'_\cdot$. In view of the definition of $\Delta^\theta g_1(s)$ and Assumptions \ref{A:UN1} and \ref{A:UN2} of $g_1$, sending $\theta\To 1$  in \eqref{eq:5.26} implies that for each $n\geq 1$,\vspace{0.1cm}
\begin{equation}\label{eq:5.27}
\hspace*{-0.2cm}\begin{array}{lll}
\Dis \left(\hat U_t^+\right)^{1+\tilde\mu_{\beta,\gamma,\eps}(t)}
&\leq& \Dis K_T\left(\hat U_{\tau_n^t}^+\right)^{1+\tilde\mu_{\beta,\gamma,\eps}(\tau_n^t)} +\int_t^{\tau_n^t}
K_s\left(\hat U_s^+\right)^{\tilde\mu_{\beta,\gamma,\eps}(s)}\left(\rho(\hat U_s^+)+\kappa(|\hat V_s|)\right){\rm d}s\vspace{0.2cm}\\
&& \Dis -\int_t^{\tau_n^t}
K_s\left(\hat U_s^+\right)^{\tilde\mu_{\beta,\gamma,\eps}(s)}\hat V_s \cdot {\rm d}B_s, \ \ \ t\in \T.
\end{array}
\end{equation}
In view of \cref{rmk:2.4}, since $\kappa(\cdot)$ is a continuous function of linear growth with $\kappa(0)=0$, it follows from the proof of Theorem 1 in \cite{FanJiangDavison2010CRA} that for each $m\geq 1$,
\begin{equation}\label{eq:5.28}
\RE\ x\in\R_+,\ \ \kappa(x)\leq (m+2A)x+\kappa\left(\frac{2A}{m+2A}\right).
\end{equation}
Let $\mathbb{P}_m$ be the probability on $(\Omega,\F)$ which is equivalent to $\mathbb{P}$ defined by
$$
\frac{{\rm d}\mathbb{P}_m}{{\rm d}\mathbb{P}}:=\exp\left\{(m+2A)\int_0^T \frac{\hat V_s}{|\hat V_s|}{\bf 1}_{|\hat V_s|>0}\cdot {\rm d}B_s+{1\over 2}(m+2A)^2 \int_0^T {\bf 1}_{|\hat V_s|>0} {\rm d}s\right\}.
$$
Note that ${\rm d}\mathbb{P}_m/{\rm d}\mathbb{P}$ has moments of all orders. Let $\E_m[\eta|\F_t]$ represent the conditional mathematical expectation of the random variable $\eta$ with respect to $\F_t$ under the probability measure $\mathbb{P}_m$. Then, in view of \eqref{eq:5.27} and \eqref{eq:5.28}, by Girsanov's theorem we deduce that for each $n,m\geq 1$,\vspace{0.1cm}
\begin{equation}\label{eq:5.29}
\begin{array}{lll}
\Dis \left(\hat U_t^+\right)^{1+\tilde\mu_{\beta,\gamma,\eps}(t)}
&\leq& \Dis \kappa\left(\frac{2A}{m+2A}\right)K_T T+K_T
\E_m\left[\left(\hat U_{\tau_n^t}^+\right)^{1+\tilde\mu_{\beta,\gamma,\eps}(\tau_n^t)}\bigg|\F_t\right]
\vspace{0.2cm}\\
&& \Dis +K_T\E_m\left[\int_t^{\tau_n^t}
\left(\hat U_s^+\right)^{\tilde\mu_{\beta,\gamma,\eps}(s)}\rho(\hat U_s^+){\rm d}s\bigg|\F_t\right],\ \ \ t\in \T.
\end{array}
\end{equation}

On the other hand, denote
$$
p:=\inf\limits_{s\in\T}\frac{1+\tilde\mu_{\beta,\gamma,\eps_0}(s)}
{1+\tilde\mu_{\beta,\gamma,\eps}(s)}\ \ \ {\rm and}\ \ P_t:=\left(\hat U_t^+\right)^{1+\tilde\mu_{\beta,\gamma,\eps}(t)},\ \ t\in\T.\vspace{0.2cm}
$$
Since $\eps=\eps_0/2$, from the definition of $\tilde\mu_{\beta,\gamma,\eps}(s)$ it is not hard to conclude that $p>1$. We pick a $q\in (1,p)$. It follows from~\cite[Lemma 6.1]{BriandDelyonHu2003SPA} and \eqref{eq:5.21} that
\[
\begin{array}{lll}
\Dis \E\left[\sup\limits_{t\in\T}|P_t|^q\right]&\leq &\Dis \E\left[\sup\limits_{t\in\T}\left[\left(1+|Y_t|+
|Y'_t|\right)^{1+\tilde\mu_{\beta,\gamma,\eps_0}(t)}
\right]^{\frac{q}{p}}\right]\leq \E\left[\left(\sup\limits_{t\in\T} \E[X|\F_t]\right)^{\frac{q}{p}}\right]\vspace{0.2cm}\\
&\leq & \Dis \frac{p}{p-q} \left(\E\left[X\right]\right)^{\frac{q}{p}}<+\infty,
\end{array}
\]
which means that $P_\cdot\in \s^q$, and then, by virtue of H\"{o}lder's inequality, \begin{equation}\label{eq:5.30}
\RE\ m\geq 1,\ \ \ \E_m\left[\sup\limits_{t\in\T}|P_t|\right]
=\E\left[\left(\sup\limits_{t\in\T}|P_t|\right)\frac{{\rm d}\mathbb{P}_m}{{\rm d}\mathbb{P}}\right]<+\infty.
\end{equation}
Thus, in view of $\hat U_T^+=0$ and the fact that $\rho(\cdot)$ is of linear growth, we can send $n$ to infinity in \eqref{eq:5.29} to obtain that for each $m\geq 1$,
\begin{equation}\label{eq:5.31}
P_t \leq \kappa\left(\frac{2A}{m+2A}\right)K_T T+K_T\E_m\left[\int_t^T
P_s^{\frac{\tilde\mu_{\beta,\gamma,\eps}(s)}
{1+\tilde\mu_{\beta,\gamma,\eps}(s)}}\rho\left(P_s^{\frac{1}
{1+\tilde\mu_{\beta,\gamma,\eps}(s)}}\right){\rm d}s\bigg|\F_t\right],\ \ \ t\in \T.\vspace{0.2cm}
\end{equation}

Finally, define the function
$$
\bar\rho (x):=\sup\limits_{t\in\T}\left[x^{\frac{\tilde\mu_{\beta,\gamma,\eps}(t)}
{1+\tilde\mu_{\beta,\gamma,\eps}(t)}}\rho\left(x^{\frac{1}
{1+\tilde\mu_{\beta,\gamma,\eps}(t)}}\right)\right],\ \ \ x\in\R_+.
$$
It is clear that $\bar\rho(0)=0$. Since $\rho(x)$ is a nondecreasing concave function on $\R_+$ with $\rho(0)=0$, it follows from \cite[Lemma 6.1]{FanJiang2013AMSE} that $\rho(x)/x, x>0$ is a non-increasing function, and then,
$$
\bar\rho (x)=x\sup\limits_{t\in\T}\left[\frac{\rho\left(x^{\frac{1}
{1+\tilde\mu_{\beta,\gamma,\eps}(t)}}\right)}{x^{\frac{1}
{1+\tilde\mu_{\beta,\gamma,\eps}(t)}}}\right]= x\ \frac{\rho\left(x^{\frac{1}
{1+\eps}}\right)}{x^{\frac{1}{1+\eps}}}=x^{\frac{\eps}
{1+\eps}}\rho\left(x^{\frac{1}{1+\eps}}\right) ,\ \ \  0<x\leq 1,
$$
and
$$
\bar\rho (x)=x\sup\limits_{t\in\T}\left[\frac{\rho\left(x^{\frac{1}
{1+\tilde\mu_{\beta,\gamma,\eps}(t)}}\right)}{x^{\frac{1}
{1+\tilde\mu_{\beta,\gamma,\eps}(t)}}}\right]=x\ \frac{\rho\left(x^{\frac{1}
{1+\tilde\mu_{\beta,\gamma,\eps}(T)}}\right)}
{x^{\frac{1}{1+\tilde\mu_{\beta,\gamma,\eps}(T)}}}
=x^{\frac{\tilde\mu_{\beta,\gamma,\eps}(T)}
{1+\tilde\mu_{\beta,\gamma,\eps}(T)}}
\rho\left(x^{\frac{1}{1+\tilde\mu_{\beta,\gamma,\eps}(T)}}\right) ,\ \ \  x\geq 1.
$$
Then, $\bar\rho(\cdot)$ is well defined on $\R_+$, nondecreasing, continuous and of linear growth, and
\begin{equation}\label{eq:5.32}
\int_{0^+}\frac{{\rm d}x}{\bar\rho(x)}= \int_{0^+}\frac{{\rm d}x}{x^{\frac{\eps}
{1+\eps}}\rho\left(x^{\frac{1}{1+\eps}}\right)}=(1+\eps)\int_{0^+}\frac{{\rm d}u}{\rho(u)}=+\infty.
\end{equation}
Now, coming back to \eqref{eq:5.31}, we have that for each $m\geq 1$,
\begin{equation}\label{eq:5.33}
P_t \leq \kappa\left(\frac{2A}{m+2A}\right)K_T T+K_T\E_m\left[\int_t^T
\bar\rho(P_s){\rm d}s\bigg|\F_t\right],\ \ \ t\in \T.\vspace{0.2cm}
\end{equation}
Then, in view of \eqref{eq:5.30}, \eqref{eq:5.32}, \eqref{eq:5.33} and the fact that
$$
\lim\limits_{m\To\infty} \kappa\left(\frac{2A}{m+2A}\right)=0,
$$
applying  \cite[Lemma 2.1]{Fan2016SPA} implies that $\ps$, for each $t\in\T$, $P_t=0$, that is, $Y_t\leq Y'_t$. The proof is complete.
\end{proof}

\begin{ex}\label{exp:5.7}
Let $\beta>0$, $\gamma>0$ and $k\geq 0$. For each $(\omega,t,y,z)\in \Omega\times\T\times\R\times\R^d$, define
$$
g_1(\omega,t,y,z):=-ke^y+\beta |y|\ln |y|{\bf 1}_{|y|\leq 1}+\gamma |z|\sqrt{|\ln |z||}{\bf 1}_{|z|\leq 1},
$$
$$
g_2(\omega,t,y,z):=|B_t(\omega)|+ke^{-y}+\beta |y|\ln |y|{\bf 1}_{|y|> 1}+\gamma |z|\sqrt{\ln |z|}{\bf 1}_{|z|>1}
$$
and
$$
g(\omega,t,y,z):=|B_t(\omega)|+k\left(e^{-y}-e^y\right)+\beta |y|\ln |y|+\gamma |z|\sqrt{|\ln |z||}.
$$
We can check that $g$ satisfies Assumptions \ref{A:EX1} and \ref{A:A3} with parameters $(|B_\cdot|+\beta/e+\gamma/\sqrt{e}+2k,\beta,\gamma)$, $g=g_1+g_2$, $g_1$ satisfies Assumptions \ref{A:UN1} and \ref{A:UN2}, and $g_2$ satisfies Assumptions \ref{A:UN3} and \ref{A:A3} with parameters $(|B_\cdot|+k,\beta,\gamma)$. Then the generator $g$ satisfies all assumptions in \cref{thm:5.3}.
\end{ex}

\begin{rmk}\label{rem:5.8}
In comparison with the related results of \cite{BahlaliElAsri2012BSM} and \cite{BahlaliKebiri2017Stochastics}, weaker assumptions on both the terminal condition and the generator are required in \cref{thm:5.3} to ensure the existence of the unbounded solution of BSDEs. In fact, the terminal condition in \cref{thm:5.3} may only need to belong to $L^p$ for some $1<p\leq 2$, whereas $p>2$ in \cite{BahlaliElAsri2012BSM} and \cite{BahlaliKebiri2017Stochastics}, and the generator $g$ only need to satisfy Assumptions \ref{A:EX1} and \ref{A:A3}, whereas the condition \eqref{eq:1.2} with $\delta=1$ and $\lambda=1/2$ in \cite{BahlaliElAsri2012BSM} and \cite{BahlaliKebiri2017Stochastics}. To ensure the uniqueness of the unbounded solution, a  locally monotone condition is required in \cite{BahlaliElAsri2012BSM} and \cite{BahlaliKebiri2017Stochastics}, which is different from the assumptions in \cref{thm:5.3}. However, we especially mention that the generator $g$ in \cref{exp:5.7} is covered by \cref{thm:5.3}, but is not covered by the results of \cite{BahlaliElAsri2012BSM} and \cite{BahlaliKebiri2017Stochastics} when $k>0$. In addition, it is noted that a general comparison theorem for the bounded solutions of BSDEs (see \cref{pro:5.6}) is also established in this section.
\end{rmk}

\section{The case of $L\exp\left(\mu (\ln (1+L))^p\right)$-integrability terminal condition with $p>1$}
\label{sec:6-mainresult--6}
\setcounter{equation}{0}

In this section, we always assume that the generator $g$ satisfies the following assumption,  which is strictly weaker than \eqref{eq:1.2} with $\delta=1$ and $\lambda>1/2$.
\begin{enumerate}
\renewcommand{\theenumi}{(A4)}
\renewcommand{\labelenumi}{\theenumi}
\item \label{A:A4} There exist three constants $\beta\geq 0$, $\gamma > 0$ and $\lambda>1/2$, and an $\R_+$-valued progressively measurable process $(\alpha_t)_{t\in\T}$  such that $\as$,
    $$
    \RE\ (y,z)\in \R\times\R^d,\ \ \
    {\rm sgn}(y)g(\omega,t,y,z)\leq \alpha_t(\omega)+\beta|y|\ln |y|{\bf 1}_{|y|>1}+\gamma |z||\ln |z||^{\lambda}.
    $$
  \end{enumerate}

For each $\eps>0$, by letting $k=1+\eps$ in \cref{Pro:2.2Inequality} we know that there exists a positive constant $C_{\lambda,\eps}>0$ depending only on parameters $(\lambda,\eps)$ such that \eqref{eq:2.5} holds. Then, our objective is to search for a positive constant $\bar\delta>0$ and a strictly increasing and strictly convex function $\phi$ such that for each $(s,x)\in \T\times \R_+$, with $\beta\geq0$ and $\gamma>0$,\vspace{0.1cm}
\begin{equation}\label{eq:6.1}
\begin{array}{l}
\Dis -{\gamma^2\over 2}\frac{\left(\phi_x(s, x)\right)^2}
{\phi_{xx}(s, x)-\bar\delta}\left[(1+\eps)4^{(\lambda-1)^+}\left|\ln \frac{\gamma\phi_x(s, x)}{\phi_{xx}(s, x)-\bar\delta}
\right|^{2\lambda}+C_{\lambda,\eps}\right]\vspace{0.2cm}\\
\Dis \ \ \ \ -\beta\phi_x(s, x) x\ln x {\bf 1}_{x>1}+\phi_s(s, x)\geq 0,\vspace{0.1cm}
\end{array}
\end{equation}
which is just \eqref{eq:2.11} with $\lambda>1/2$, $\delta=1$, $k=1+\eps$ and $C_{k,\lambda}=C_{\lambda, \eps}$.\vspace{0.2cm}

Now, let $\eps>0$, $\beta\geq0$, $\gamma>0$, $\lambda>1/2$ and $\mu_s,\nu_s:\T\To\R_+$ be two increasing and continuously differential functions with $\mu_0=\eps$ and $\nu_0=0$. Define\vspace{-0.1cm}
\begin{equation}\label{eq:6.2}
\Dis \tilde k:={\rm e}+\left(\frac{2\lambda\mu_T}{\gamma}\right)^2.\vspace{-0.1cm}
\end{equation}
We choose the following function
$$
\phi(s,x):=(x+\tilde k)\exp\left(\mu_s \left(\ln (x+\tilde k)\right)^{2\lambda}+\nu_s\right), \ \ \ (s,x)\in \T\times \R_+
$$
to explicitly solve the inequality \eqref{eq:6.1}. For each $(s,x)\in \T\times \R_+$, a simple computation gives
\begin{equation}\label{eq:6.3}
\phi_x(s,x)=\phi(s,x){2\lambda \mu_s\left(\ln (x+\tilde k)\right)^{2\lambda-1}+1\over x+\tilde k}>0,\vspace{0.1cm}
\end{equation}
\begin{equation}\label{eq:6.4}
\phi_{xx}(s,x)=\phi(s,x){2\lambda\mu_s\left(\ln (x+\tilde k)\right)^{2\lambda-1}\left[
\ln (x+\tilde k)+2\lambda \mu_s \left(\ln (x+\tilde k)\right)^{2\lambda}+(2\lambda-1)\right]\over (x+\tilde k)^2 \ln (x+\tilde k)}>0
\end{equation}
and\vspace{-0.1cm}
\begin{equation}\label{eq:6.5}
\phi_s(s,x)=\phi(s,x)\left(\mu'_s\left(\ln \left(x+\tilde k\right)\right)^{2\lambda}+\nu'_s \right)>0.\vspace{0.1cm}
\end{equation}
Furthermore, for each $(s,x)\in \T\times \R_+$, in view of $\mu_s\geq \mu_0=\eps$, we have
\begin{equation}\label{eq:6.6}
\phi(s,x){2\lambda(2\lambda-1)\mu_s\left(\ln (x+\tilde k)\right)^{2\lambda-1}\over (x+\tilde k)^2\ln (x+\tilde k)}\geq \Dis \exp\left(\eps \left(\ln (x+\tilde k)\right)^{2\lambda}\right)
{2\lambda(2\lambda-1)\eps\left(\ln (x+\tilde k)\right)^{2\lambda-1}\over (x+\tilde k)\ln (x+\tilde k)}.\vspace{0.2cm}
\end{equation}
Observe that the function in the right hand side of \eqref{eq:6.6} is positive and continuous on $\R_+$, and tends to infinity as $x\To +\infty$. It follows that there exists a constant $\delta_{\lambda,\eps}>0$ depending only on $(\lambda,\eps)$ such that
$$
\phi(s,x){2\lambda(2\lambda-1)\mu_s\left(\ln (x+\tilde k)\right)^{2\lambda-1}\over (x+\tilde k)^2\ln (x+\tilde k)}\geq \delta_{\lambda,\eps},\ \ \ (s,x)\in \T\times \R_+.\vspace{0.2cm}
$$
It then follows from \eqref{eq:6.4} that\vspace{0.1cm}
\begin{equation}\label{eq:6.7}
\phi_{xx}(s,x)-\delta_{\lambda,\eps}\geq\phi(s,x)
{2\lambda\mu_s\left(\ln (x+\tilde k)\right)^{2\lambda-1}\left[2\lambda\mu_s\left(\ln (x+\tilde k)\right)^{2\lambda-1}+1\right]\over (x+\tilde k)^2},\ \ (s,x)\in \T\times \R_+.\vspace{0.2cm}
\end{equation}
For each $(s,x)\in \T\times \R_+$, by \eqref{eq:6.3} and \eqref{eq:6.7} we have\vspace{0.1cm}
\begin{equation}\label{eq:6.8}
{\gamma^2\over 2}\frac{\left(\phi_x(s,x)\right)^2}{\phi_{xx}(s,x)-\delta_{\lambda,\eps}}\leq {\gamma^2 \over 2}\phi(s,x)\left(1+\frac{1}{2\lambda\mu_s\left(\ln (x+\tilde k)\right)^{2\lambda-1}}\right)\leq {\gamma^2 \over 2}\phi(s,x)\left(1+\frac{1}{2\lambda\mu_s}\right).
\end{equation}
And, from the definition of $\tilde k$ it can be directly verified that for each $(s,x)\in \T\times \R_+$,
\begin{equation}\label{eq:6.9}
\Dis \ln (x+\tilde k)+2\lambda\mu_s \left(\ln (x+\tilde k)\right)^{2\lambda}+\left(2\lambda-1\right)
\leq \Dis \frac{1+\eps}{\eps}\left[
\ln (x+\tilde k)+2\lambda\mu_s \left(\ln (x+\tilde k)\right)^{2\lambda}\right].
\end{equation}
In view of \eqref{eq:6.2}, \eqref{eq:6.3}, \eqref{eq:6.4}, \eqref{eq:6.9} and \eqref{eq:6.7}, we have\vspace{0.1cm}
\[
\left\{
\begin{array}{l}
\Dis \frac{\gamma\phi_x(s,x)}{\phi_{xx}(s,x)-\delta_{\lambda,\eps}}\geq {\gamma\eps\over 2\lambda(1+\eps)\mu_s}\frac{x+\tilde k}{\left(\ln (x+\tilde k)\right)^{2\lambda-1}}\geq {\gamma \over 2\lambda \mu_T}\sqrt{x+\tilde k}\geq {\gamma \over 2\lambda \mu_T}\sqrt{\tilde k}\geq 1;\vspace{0.3cm}\\
\Dis \frac{\gamma\phi_x(s,x)}{\phi_{xx}(s,x)-\delta_{\lambda,\eps}}\leq {\gamma\over 2\lambda \mu_s}\frac{x+\tilde k}{\left(\ln (x+\tilde k)\right)^{2\lambda-1}}\leq {\gamma\over 2\lambda \eps}(x+\tilde k),\vspace{0.1cm}
\end{array}
\right.
\]
which yields the following
\begin{equation}\label{eq:6.10}
\left|\ln \frac{\gamma\phi_x(s,x)}{\phi_{xx}(s,x)-\delta_{\lambda,\eps}}\right|^{2\lambda}\leq
 2^{2\lambda -1}\left|\ln {\gamma \over 2\lambda\eps}\right|^{2\lambda}+ 2^{2\lambda -1}\left( \ln (x+\tilde k)\right)^{2\lambda}.\vspace{0.3cm}
\end{equation}

In the sequel, observe that
$$x\ln x {\bf 1}_{x>1}\leq (x+\tilde k)\ln (x+\tilde k),\ \ \ x\in\R_+.$$
We substitute \eqref{eq:6.8}, \eqref{eq:6.10}, \eqref{eq:6.3} and \eqref{eq:6.5} into the left side of \eqref{eq:6.1} with $\delta_{\lambda,\eps}$ instead of $\bar\delta$ to get \vspace{0.1cm}
$$
\begin{array}{ll}
&\Dis -{\gamma^2\over 2}\frac{\left(\phi_x(s, x)\right)^2}
{\phi_{xx}(s, x)-\delta_{\lambda,\eps}}\left[ 4^{(\lambda-1)^+}(1+\eps)\left|\ln \frac{\gamma\phi_x(s, x)}{\phi_{xx}(s, x)-\delta_{\lambda,\eps}}
\right|^{2\lambda}+C_{\lambda,\eps}\right]\vspace{0.2cm}\\
&\Dis \ \ \ \ -\beta\phi_x(s,x) x \ln x {\bf 1}_{x>1}+\phi_s(s, x)
\vspace{0.2cm}\\
\geq & \Dis
-{\gamma^2\over 2}\phi(s,x)\left(1+\frac{1}{2\lambda\mu_s} \right)\left[(1+\eps)k_{\lambda}\left( \ln (x+\tilde k)\right)^{2\lambda}+\bar C_{\lambda,\eps}\right]
\vspace{0.2cm}\\
& \Dis -\beta \phi(s, x) \left[2\lambda\mu_s\left(\ln (x+\tilde k)\right)^{2\lambda}+\ln (x+\tilde k)\right]\vspace{0.2cm}\\
&\Dis +\phi(s,x)\left(\mu'_s\left(\ln \left(x+\tilde k\right)\right)^{2\lambda}+\nu'_s \right), \ \ \ \ (s,x)\in \T\times \R_+\vspace{0.1cm}
\end{array}
$$
where
$$
k_{\lambda}:= 2^{2(\lambda-1)^+ +2\lambda-1}\ \ {\rm and }\ \
\bar C_{\lambda,\eps}:=(1+\eps)k_{\lambda}\left|\ln {\gamma \over 2\lambda \eps}\right|^{2\lambda}+C_{\lambda,\eps}.\vspace{0.2cm}
$$
Furthermore, in view of $2\lambda>1$, by Young's inequality it is not very hard to verify that there exists a constant $\Tilde C_{\beta,\lambda,\eps}>0$ depending only on $(\beta,\lambda,\eps)$ such that for each $(s,x)\in \T\times \R_+$,
$$
\beta \ln (x+\tilde k)\leq \eps\left(\ln (x+\tilde k)\right)^{2\lambda}+\Tilde C_{\beta,\lambda,\eps}
$$
Then, for each $(s,x)\in \T\times \R_+$, we have \vspace{0.1cm}
$$
\begin{array}{ll}
&\Dis -{\gamma^2\over 2}\frac{\left(\phi_x(s, x)\right)^2}
{\phi_{xx}(s, x)-\delta_{\lambda,\eps}}\left[ 4^{(\lambda-1)^+}(1+\eps)\left|\ln \frac{\gamma\phi_x(s, x)}{\phi_{xx}(s, x)-\delta_{\lambda,\eps}}
\right|^{2\lambda}+C_{\lambda,\eps}\right]\vspace{0.2cm}\\
&\Dis \ \ \ \ -\beta\phi_x(s,x) x \ln x {\bf 1}_{x>1}+\phi_s(s, x)
\vspace{0.3cm}\\
\geq & \Dis \phi(s, x)\left\{\left[
-\frac{\gamma^2(1+\eps)k_\lambda}{2}\left(1+\frac{1}{2\lambda\mu_s} \right)-2\beta\lambda \mu_s-\eps+\mu'_s \right]\left(\ln (x+\tilde k)\right)^{2\lambda}\right. \vspace{0.2cm}\\
&\Dis \hspace*{1.5cm} \left.+ \left[-\frac{\gamma^2}{2}\left(1+\frac{1}{2\lambda\mu_s} \right)\bar C_{\lambda,\eps} -\Tilde C_{\beta,\lambda,\eps}+\nu'_s\right]\right\}.\vspace{0.1cm}
\end{array}
$$
Thus, \eqref{eq:6.1} holds if the functions $\mu_s, \nu_s\in \T$ satisfies\vspace{0.2cm}
\begin{equation}\label{eq:6.11}
\mu'_s=\frac{\gamma^2(1+\eps)k_\lambda}{2}
\left(1+\frac{1}{2\lambda\mu_s} \right)+2\beta\lambda \mu_s+\eps\ \ \ s\in [0,T]
\end{equation}
and
\begin{equation}\label{eq:6.12}
\nu_s=\frac{\gamma^2}{4\lambda}\bar C_{\lambda,\eps} \int_0^s \frac{1}{\mu_r} {\rm d}r+\left(\frac{\gamma^2}{2}\bar C_{\lambda,\eps}+\Tilde C_{\beta,\lambda,\eps}\right)s,\ \ \ s\in [0,T].\vspace{0.2cm}
\end{equation}

It is not very hard to verify that for each $\beta\geq 0$, $\gamma > 0$, $\lambda>1/2$ and $\eps>0$, there exists a unique strictly increasing and continuous function with $\eps$ at the origin satisfying \eqref{eq:6.11}. We denote this unique solution by $\breve \mu_{\beta,\gamma,\lambda,\eps}(\cdot)$, i.e., $\breve \mu_{\beta,\gamma,\lambda,\eps}(0)=\eps$ and
\begin{equation}\label{eq:6.13}
\breve\mu'_{\beta,\gamma,\lambda,\eps}(s)
=\frac{\gamma^2(1+\eps)k_\lambda}{2}
\left(1+\frac{1}{2\lambda \breve\mu_{\beta,\gamma,\lambda,\eps}(s)} \right)+2\beta\lambda \breve\mu_{\beta,\gamma,\lambda,\eps}(s)+\eps,\ \ \ s\in [0,T].
\end{equation}
We also denote, in view of \eqref{eq:6.12},
\begin{equation}\label{eq:6.14}
\breve \nu_{\beta,\gamma,\lambda,\eps}(s):=\frac{\gamma^2}{4\lambda}\breve C_{\lambda,\eps} \int_0^s \frac{1}{\breve\mu_{\beta,\gamma,\lambda,\eps}(r)} {\rm d}r+\left(\frac{\gamma^2}{2}\breve C_{\lambda,\eps}+\Tilde C_{\beta,\lambda,\eps}\right)s,\ \ \ s\in [0,T].\vspace{0.1cm}
\end{equation}

Moreover, it can be directly checked that as $\eps\To 0^+$, the function $\breve \mu_{\beta,\gamma,\lambda,\eps}(\cdot)$ tends decreasingly to the unique solution $\breve \mu_{\beta,\gamma,\lambda}^0(\cdot)$ of the following ODE:
\begin{equation}\label{eq:6.15}
\breve \mu_{\beta,\gamma,\lambda}^0(0)=0\ \ \ {\rm and}\ \ \ \left(\breve \mu_{\beta,\gamma,\lambda}^0(s)\right)'=
\frac{\gamma^2k_\lambda}{2}
\left(1+\frac{1}{2\lambda \breve \mu_{\beta,\gamma,\lambda}^0(s)} \right)+2\beta\lambda \breve \mu_{\beta,\gamma,\lambda}^0(s),\ \ \ s\in (0,T].\vspace{0.3cm}
\end{equation}

In summary, we have proved the following proposition.

\begin{pro}\label{pro:6.1}
For each $\beta\geq 0$, $\gamma > 0$, $\lambda>1/2$ and $\eps>0$, define the function
\begin{equation}\label{eq:6.16}
\breve \varphi(s,x;\eps):=(x+\breve k)\exp\left(\breve \mu_{\beta,\gamma,\lambda,\eps}(s)\left(\ln \left(x+\breve k\right)\right)^{2\lambda}+\breve \nu_{\beta,\gamma,\lambda,\eps}(s)\right), \ \ \ (s,x)\in [0,T]\times \R_+
\end{equation}
with $\breve \mu_{\beta,\gamma,\lambda,\eps}(\cdot)$ and $\breve \nu_{\beta,\gamma,\lambda,\eps}(\cdot)$ being respectively defined in \eqref{eq:6.13} and \eqref{eq:6.14}, and
\[
\breve k:={\rm e}+\left(\frac{2\lambda\breve \mu_{\beta,\gamma,\lambda,\eps}(T)}{\gamma}\right)^2,
\]
which is $\tilde k$ in \eqref{eq:6.2} with $\breve \mu_{\beta,\gamma,\lambda,\eps}(T)$ instead of $\mu_T$. Then we have
\begin{itemize}
\item [(i)] $\breve \varphi(\cdot,\cdot;\eps)\in \mathcal {C}^{1,2}\left([0,T]\times \R_+\right)$;
\item [(ii)] There exists a positive constant $\delta_{\lambda,\eps}>0$ depending only on $(\lambda,\eps)$ such that $\breve \varphi(\cdot,\cdot;\eps)$ satisfies the inequality \eqref{eq:2.4}, i.e., for each $(s,x,z)\in [0,T]\times \R_+\times\R^d$, we have\vspace{-0.1cm}
$$
-\breve \varphi_x(s,x;\eps)\left(\beta x\ln x{\bf 1}_{x>1}+\gamma |z||\ln |z||^{\lambda}\right)+{1\over 2}\left(\breve\varphi_{xx}(s,x;\eps)-\delta_{\lambda,\eps}\right)|z|^2+\breve \varphi_s(s,x;\eps)\geq 0.
$$
\end{itemize}
\end{pro}

Now, for $\lambda>1/2$, we define the function
\begin{equation}\label{eq:6.17}
\breve\psi(x,\mu):=x\exp\left(\mu\left(\ln (1+x)\right)^{2\lambda}\right),\ \ (x,\mu)\in \R_+\times \R_+.
\end{equation}

The following existence and uniqueness theorem is one main result of this section.

\begin{thm}\label{thm:6.2}
Let the functions $\breve \mu_{\beta,\gamma,\lambda,\eps}(\cdot)$, $\breve \mu_{\beta,\gamma,\lambda}^0(\cdot)$ and $\breve \psi(\cdot,\cdot)$ be respectively defined in \eqref{eq:6.13}, \eqref{eq:6.15}, and \eqref{eq:6.17}, $\xi$ be a terminal condition and $g$ be a generator which is continuous in the state variables $(y,z)$. If $g$ satisfies Assumptions \ref{A:EX1} and \ref{A:A4} with parameters $(\alpha_\cdot, \beta, \gamma, \lambda)$, and there exists a positive constant $\mu>\breve \mu_{\beta,\gamma,\lambda}^0(T)$ such that\vspace{-0.1cm}
\begin{equation}\label{eq:6.18}
\E\left[\breve \psi\left(|\xi|+\int_0^T \alpha_t {\rm d}t,\ \mu\right)\right]<+\infty,\vspace{0.1cm}
\end{equation}
then BSDE$(\xi,g)$ admits a solution $(Y_t,Z_t)_{t\in\T}$ such that $\left(\breve\psi\left(|Y_t|+\int_0^t \alpha_s {\rm d}s,\ \breve \mu_{\beta,\gamma,\lambda,\eps}(t)\right)\right)_{t\in\T}$ belongs to class (D) for some $\eps>0$ and $Z_\cdot\in \mcal^2$, where $\eps$ is the unique positive constant satisfying $\breve \mu_{\beta,\gamma,\lambda,\eps}(T)=\mu$. And, there exists a constant $\delta_{\lambda,\eps}>0$ depending only on $(\lambda,\eps)$ such that $\ps$, for each $t\in\T$,
\begin{equation}\label{eq:6.19}
\begin{array}{lll}
|Y_t|&\leq& \Dis  \breve \psi(|Y_t|+\int_0^t \alpha_s {\rm d}s,\ \breve \mu_{\beta,\gamma,\lambda,\eps}(t))+\frac{\delta_{\lambda,\eps}}{2}\E\left[\left.
\int_t^T |Z_s|^2{\rm d}s\right| \F_t\right]\vspace{0.2cm}\\
&\leq & \Dis C\E\left[\left.\breve \psi\left(|\xi|+\int_0^T\alpha_t {\rm d}t,\ \mu\right)\right|\F_t\right]+C,
\end{array}
\end{equation}
where $C>0$ is a positive constant depending only on $(\beta,\gamma,\lambda,\eps,T)$.

Moreover, if the inequality \eqref{eq:6.18} holds for some $\mu>2q3^{2\lambda-1} \breve \mu_{\beta,\gamma,\lambda}^0(T)$ with $q>1$, then BSDE $(\xi,g)$ admits a solution $(Y_t,Z_t)_{t\in\T}$ such that the process
$$\left(\breve\psi\left(|Y_t|+\int_0^t \alpha_s {\rm d}s,\ 2q3^{2\lambda-1} \breve \mu_{\beta,\gamma,\lambda,\eps}(t)\right)\right)_{t\in\T}$$
belongs to class (D), where $\eps$ is the unique positive constant satisfying $2q3^{2\lambda-1} \breve \mu_{\beta,\gamma,\lambda,\eps}(T)=\mu$. And, if the generator $g$ further satisfies Assumption \ref{A:UN3}, then the solution satisfying the preceding condition is also unique.
\end{thm}

In order to prove this theorem, we need the following two propositions.

\begin{pro}\label{pro:6.3}
Let functions $\breve \mu_{\beta,\gamma,\lambda,\eps}(s)$, $\breve \varphi(s,x;\eps)$ and $\breve \psi(x,\mu)$ be respectively defined on \eqref{eq:6.13}, \eqref{eq:6.16} and \eqref{eq:6.17}. Then, there exists a universal constant $K>0$ depending on $(\beta,\gamma,\lambda,\eps,T)$ such that\vspace{-0.1cm}
\begin{equation}\label{eq:6.20}
\RE\ (s,x)\in [0,T]\times\R_+,\ \ \ \breve \psi(x,\breve \mu_{\beta,\gamma,\lambda,\eps}(s))\leq \breve \varphi(s,x;\eps)\leq K\breve \psi(x,\breve \mu_{\beta,\gamma,\lambda,\eps}(s))+K.
\end{equation}
\end{pro}

\begin{proof}
The first inequality in \eqref{eq:6.20} is obvious. We now prove the second inequality. In fact,
$$
\begin{array}{lll}
\Dis {\breve\varphi(s,x;\eps) \over\breve\psi(x,\breve\mu_{\beta,\gamma,\lambda,\eps}(s))+1}&=&\Dis \frac {(x+\breve k)\exp\left(\breve\mu_{\beta,\gamma,\lambda,\eps}(s)\left(\ln \left(x+\breve k\right)\right)^{2\lambda}+\breve\nu_{\beta,\gamma,\lambda,\eps}(s)\right)}
{x\exp\left(\breve\mu_{\beta,\gamma,\lambda,\eps}(s)\left(\ln (1+x)\right)^{2\lambda}\right)+1}\vspace{0.2cm}\\
&\leq & \Dis {x+\breve k\over x}\exp\left(\breve\mu_{\beta,\gamma,\lambda,\eps}(T)\left(\left(\ln \left(x+\breve k\right)\right)^{2\lambda}-\left(\ln \left(1+x\right)\right)^{2\lambda}\right)
+\breve\nu_{\beta,\gamma,\lambda,\eps}(T)\right)\vspace{0.2cm}\\
&=: & \Dis H_1(x;\beta,\gamma,\lambda,\eps,T),\ \ \ \ \ (s,x)\in [0,T]\times [1,+\infty),
\end{array}
$$
with $\breve k$ being defined in \cref{pro:6.1}. And, in the case of $x\in [0,1]$,
$$
\begin{array}{lll}
\Dis {\breve\varphi(s,x;\eps) \over\breve\psi(x,\breve\mu_{\beta,\gamma,\lambda,\eps}(s))+1}&\leq & \Dis (1+\breve k)\exp\left(\breve\mu_{\beta,\gamma,\lambda,\eps}(T)\left(\ln \left(1+\breve k\right)\right)^{2\lambda}+\breve\nu_{\beta,\gamma,\lambda,\eps}(T)\right)\vspace{0.2cm}\\
&=:& \Dis H_2(\beta,\gamma,\lambda,\eps,T),\ \ \ s\in\T.
\end{array}
$$
Hence, for all $x\in \R_+$, we have
\begin{equation}\label{eq:6.21}
{\breve\varphi(s,x;\eps) \over\breve\psi(x,\breve\mu_{\beta,\gamma,\lambda,\eps}(s))+1}\leq H_1(x;\beta,\gamma,\lambda,\eps,T){\bf 1}_{x\geq 1}+H_2(\beta,\gamma,\lambda,\eps,T){\bf 1}_{0\leq x<1},\ \ \ s\in\T.\vspace{0.1cm}
\end{equation}
Thus, in view of \eqref{eq:6.21} and the fact that the function $H_1(x;\beta, \gamma,\lambda,\eps, T)$ is continuous on $[1,+\infty)$ and tends to $\exp(\breve\nu_{\beta,\gamma,\lambda,\eps}(T))$ as $x\To +\infty$, the second inequality in \eqref{eq:6.20} follows immediately.
\end{proof}

\begin{pro}\label{pro:6.4}
Let functions $\breve \mu_{\beta,\gamma,\lambda,\eps}(\cdot)$ and $\breve \psi(\cdot,\cdot)$ be respectively defined in \eqref{eq:6.13} and \eqref{eq:6.17}, $\xi$ be a terminal condition and $g$ be a generator which is continuous in $(y,z)$. If $g$ satisfies Assumption \ref{A:A4} with parameters $(\alpha_\cdot, \beta, \gamma, \lambda)$, $|\xi|+\int_0^T \alpha_t {\rm d}t$ is a bounded random variable, and $(Y_t,Z_t)_{t\in\T}$ is a solution of BSDE $(\xi,g)$ such that $Y_\cdot$ is bounded, then for each $\eps>0$, there exists a constant $\delta_{\lambda,\eps}>0$ depending only on $(\lambda,\eps)$ such that $\ps$, for each $t\in\T$, we have
\begin{equation}\label{eq:6.22}
\begin{array}{lll}
|Y_t|&\leq & \Dis \breve \psi(|Y_t|+\int_0^t \alpha_s {\rm d}s,\ \breve\mu_{\beta,\gamma,\lambda,\eps}(t))+\frac{\delta_{\lambda,\eps}}{2}\E\left[\left.
\int_t^T |Z_s|^2{\rm d}s\right| \F_t\right]\vspace{0.2cm}\\
&\leq & \Dis C\E\left[\left.\breve \psi\left(|\xi|+\int_0^T\alpha_t {\rm d}t,\ \breve\mu_{\beta,\gamma,\lambda,\eps}(T)\right)\right|\F_t\right]+C,
\end{array}
\end{equation}
where $C$ is a positive constant depending only on $(\beta,\gamma,\lambda,\eps,T)$.
\end{pro}

\begin{proof}
Define
$$
\breve Y_t:=|Y_t|+\int_0^t \alpha_s {\rm d}s\ \ \ \
{\rm and}\ \ \ \ \breve Z_t:={\rm sgn}(Y_t)Z_t,\ \ \ \ t\in \T.
$$
Using It\^{o}-Tanaka's formula, we have
$$
\breve Y_t=\breve Y_T+\int_t^T \left({\rm sgn}(Y_s)g(s,Y_s,Z_s)-\alpha_s\right){\rm d}s-\int_t^T \breve Z_s \cdot {\rm d}B_s-\int_t^T {\rm d}L_s, \ \ \ t\in\T,
$$
where $L_\cdot$ is the local time of $Y_\cdot$ at the origin. Now, we fix $\eps>0$ and apply It\^{o}-Tanaka's formula to the process $\breve\varphi(s, \breve Y_s;\eps)$, where the function $\breve\varphi(\cdot,\cdot;\eps)$ is defined in \eqref{eq:6.16}, to derive, in view of \ref{A:A4},
$$
\begin{array}{lll}
\Dis {\rm d}\breve\varphi(s,\breve Y_s;\eps)
&=&\Dis \breve\varphi_x(s,\breve Y_s;\eps)
\left(-{\rm sgn}(Y_s)g(s,Y_s,Z_s)+\alpha_s\right){\rm d}s+\breve\varphi_x(s,\breve Y_s;\eps)\breve Z_s \cdot {\rm d}B_s+\breve\varphi_x(s,\breve Y_s;\eps){\rm d}L_s\vspace{0.1cm}\\
&&\Dis +{1\over 2}\breve\varphi_{xx}(s,\breve Y_s;\eps)|Z_s|^2{\rm d}s+\breve\varphi_s(s,\breve Y_s;\eps){\rm d}s\vspace{0.2cm}\\
&\geq &\Dis \bigg[-\breve\varphi_x(s,\breve Y_s;\eps)\left(\beta |Y_s|\ln |Y_s|{\bf 1}_{|Y_s|>1}+\gamma |Z_s||\ln |Z_s||^{\lambda}\right)\vspace{0.1cm}\\
&&\Dis
\ \ +{1\over 2}\breve\varphi_{xx}(s,\breve Y_s;\eps)|Z_s|^2+\varphi_s(s,\breve Y_s;\eps)\bigg]{\rm d}s+\breve\varphi_x(s,\breve Y_s;\eps)\breve Z_s \cdot {\rm d}B_s, \ \ s\in\T.
\end{array}
$$
Furthermore, from the inequality $|Y_s|\ln |Y_s|{\bf 1}_{|Y_s|>1}\leq \breve Y_s \ln \breve Y_s{\bf 1}_{\breve Y_s>1}$ and \cref{pro:6.1},  we see that there is a constant $\delta_{\lambda,\eps}>0$ depending only on $(\lambda,\eps)$ such that \vspace{-0.1cm}
\begin{equation}\label{eq:6.23}
{\rm d}\breve\varphi(s,\breve Y_s;\eps)\geq \frac{\delta_{\lambda,\eps}}{2}|Z_s|^2 {\rm d}s+\breve\varphi_x(s,\breve Y_s;\eps)\breve Z_s \cdot {\rm d}B_s,\ \ s\in \T.
\end{equation}
Let us consider, for each integer $n\geq 1$ and each $t\in \T$,  the following stopping time
$$
\tau_n^t:=\inf\left\{s\in [t,T]: \int_t^s \left[\breve \varphi_x(r,\breve Y_r;\eps)\right]^2|\breve Z_r|^2{\rm d}r\geq n \right\}\wedge T.
$$
In view of the inequality \eqref{eq:6.23}, we have that for each $n\geq 1$,
$$
\breve\varphi(t,\breve Y_t;\eps)+\frac{\delta_{\lambda,\eps}}{2}\E\left[\left.
\int_t^{\tau_n^t} |Z_s|^2{\rm d}s\right| \F_t\right]\leq \E\left[\left. \breve\varphi(\tau_n^t,\breve Y_{\tau_n^t};\eps) \right|\F_t\right],\ \ \ t\in \T.
$$
Thus, by \cref{pro:6.3}, there exists a constant $K>0$ depending only on $(\beta, \gamma, \lambda,\eps,T)$ such that
$$
\breve\psi(\breve Y_t,\ \breve \mu_{\beta,\gamma,\lambda,\eps}(t))\leq \breve\varphi(t,\breve Y_t;\eps)\leq \E\left[\left. \breve\varphi(\tau_n^t,\breve Y_{\tau_n^t};\eps) \right|\F_t\right]\leq K \E\left[\left. \breve\psi(\breve Y_{\tau_n^t},\ \breve \mu_{\beta,\gamma,\lambda,\eps}(\tau_n^t))\right|\F_t\right]+K,\ \ t\in\T.
$$
And, since $\breve\psi(x,\mu)$ is increasing in $x$, we have that for each $n\geq 1$,
$$
\begin{array}{lll}
|Y_t|&\leq& \Dis \breve\psi\left(|Y_t|+\int_0^t \alpha_s {\rm d}s,\ \breve \mu_{\beta,\gamma,\lambda,\eps}(t)\right)+\frac{\delta_{\lambda,\eps}}{2}\E\left[\left.
\int_t^{\tau_n^t} |Z_s|^2{\rm d}s\right| \F_t\right]\vspace{0.2cm}\\
&\leq& \Dis K \E\left[\left. \breve\psi\left(|Y_{\tau_n^t}|+\int_0^{\tau_n^t} \alpha_s {\rm d}s,\ \breve \mu_{\beta,\gamma,\lambda,\eps}(\tau_n^t)\right)\right|\F_t\right]+K,\ \ t\in\T,
\end{array}
$$
which gives the desired inequality \eqref{eq:6.22} in  the limit as $n$ to infinity. The proof is completed.
\end{proof}

\begin{rmk}\label{rmk:6.6}
Assertions of \cref{pro:6.4} are still true if $(|\xi|, |Y_t|)$ is replaced with $(\xi^+, Y^+_t)$, and Assumption \ref{A:A4} is replaced with the following one \ref{A:A4'}:
\begin{enumerate}
\renewcommand{\theenumi}{(A4')}
\renewcommand{\labelenumi}{\theenumi}
\item\label{A:A4'} There exist three constants $\beta\geq 0$, $\gamma > 0$ and $\lambda\in (0,1/2)$, and an $\R_+$-valued progressively measurable process $(\alpha_t)_{t\in\T}$  such that $\as$,
    \[
    \RE\ (y,z)\in \R_+\times\R^d,\ \ \
    g(\omega,t,y,z)\leq \alpha_t(\omega)+\beta|y|(\ln |y|)^{\lambda+\frac{1}{2}}{\bf 1}_{|y|>1}+\gamma |z||\ln |z||^{\lambda}.\vspace{-0.2cm}
    \]
\end{enumerate}
To show this, it is sufficient to use $(Y_\cdot^+,  {\bf 1}_{Y_\cdot>0}, \frac{1}{2}L_\cdot)$ instead of $(|Y_\cdot|,  {\rm sgn} (Y_\cdot), L_\cdot)$ in the  proof.\vspace{0.1cm}
\end{rmk}

Now, we prove \cref{thm:6.2}.

\begin{proof}[Proof of \cref{thm:6.2}]
For each pair of positive integers $n,p\geq 1$, let $\xi^{n,p}$ and $g^{n,p}$ be defined in \eqref{eq:2.1}, and $(Y^{n,p}_\cdot,Z^{n,p}_\cdot)$ be the minimal (maximal) bounded solution of \eqref{eq:2.2}. It is easy to check that the generator $g^{n,p}$ satisfies Assumption \ref{A:A4} with $\alpha_\cdot$ being replaced with $\alpha_\cdot \wedge (n\vee p)$.

Now, we assume that there exists a positive constant $\mu>\breve \mu_{\beta,\gamma,\lambda}^0(T)$ such that \eqref{eq:6.18} holds. From the definitions of $\breve \mu_{\beta,\gamma,\lambda,\eps}(\cdot)$ and $\breve \mu_{\beta,\gamma,\lambda}^0(\cdot)$ in \eqref{eq:6.13} and \eqref{eq:6.15}, it is not very difficult to find a (unique) constant $\eps>0$ satisfying $\breve \mu_{\beta,\gamma,\lambda,\eps}(T)=\mu$. Then, applying \cref{pro:6.4} with this $\eps$ to BSDE \eqref{eq:2.2} yields that there exists a constant $\delta_{\lambda,\eps}>0$ depending only on $(\lambda,\eps)$ such that $\ps$, for all $n,p\geq 1$,\vspace{0.1cm}
\begin{equation}\label{eq:6.24}
\begin{array}{lll}
|Y^{n,p}_t|&\leq & \Dis \breve\psi(|Y^{n,p}_t|+\int_0^t [\alpha_s\wedge (n\vee p)]{\rm d}s,\ \breve \mu_{\beta,\gamma,\lambda,\eps}(t))+\frac{\delta_{\lambda,\eps}}{2}\E\left[\left.
\int_t^T |Z_s^{n,p}|^2{\rm d}s\right| \F_t\right]\vspace{0.2cm}\\
&\leq& \Dis  C\E\left[\left.\breve\psi\left(|\xi^{n,p}|+\int_0^T[\alpha_t\wedge (n\vee p)]{\rm d}t,\ \breve \mu_{\beta,\gamma,\lambda,\eps}(T)\right)\right|\F_t\right]+C\vspace{0.2cm}\\
&\leq& \Dis C\E\left[\left.\breve\psi\left(|\xi|+\int_0^T \alpha_t {\rm d}t,\ \breve \mu_{\beta,\gamma,\lambda,\eps}(T)\right)\right|\F_t\right]+C\vspace{0.3cm}\\
&=& \Dis C\E\left[\left.\breve\psi\left(|\xi|+\int_0^T \alpha_t {\rm d}t,\ \mu\right)\right|\F_t\right]+C\ =:X_t,\ \ \ \ t\in\T,\vspace{0.2cm}
\end{array}
\end{equation}
where $C>0$ is a positive constant depending only on $(\beta,\gamma,\lambda,\eps,T)$. Thus, in view of \eqref{eq:6.24}, there is an $\R_+$-valued, progressively measurable and continuous process $(X_t)_{t\in\T}$ such that\vspace{-0.1cm}
$$
\as,\ \ \RE\ n,p\geq 1,\ \ \ |Y^{n,p}_\cdot|\leq X_\cdot.\vspace{-0.1cm}
$$
Now, we can apply \cref{Pro:2.1--stability} to obtain the existence of a progressively measurable process $(Z_t)_{t\in\T} $ such that $(Y_\cdot:=\inf_p\sup_n Y^{n,p}_\cdot, \ Z_\cdot)$ is a solution to BSDE$(\xi,g)$.

Furthermore, sending $n$ and $p$ to infinity in \eqref{eq:6.24} yields the desired inequality \eqref{eq:6.19}, and then the process $\left(\breve\psi\left(|Y_t|+\int_0^t \alpha_s {\rm d}s,\breve\mu_{\beta,\gamma,\lambda,\eps}(t)\right)\right)_{t\in\T}$ belongs to class (D). \vspace{0.2cm}

In the sequel, we assume that \eqref{eq:6.18} holds for some $\mu>2q3^{2\lambda-1} \breve \mu_{\beta,\gamma,\lambda}^0(T)$ with $q>1$. Then, there exists a unique constant $\eps>0$ such that $2q3^{2\lambda-1} \breve \mu_{\beta,\gamma,\lambda,\eps}(T)=\mu$. In view of the definition of $\breve \mu_{\beta,\gamma,\lambda,\eps}(\cdot)$ and by virtue of the analysis from \eqref{eq:6.1} to \eqref{eq:6.15}, it is not very difficult to verify that \cref{pro:6.1} still holds when the functions $\breve \mu_{\beta,\gamma,\lambda,\eps}(\cdot)$ in \eqref{eq:6.14}, \eqref{eq:6.16} and the definition of $\breve k$ are all replaced with $2q3^{2\lambda-1} \breve \mu_{\beta,\gamma,\lambda,\eps}(\cdot)$, and then BSDE$(\xi,g)$ admits a solution $(Y_t,Z_t)_{t\in\T}$ such that  the process
$$\left(\breve\psi\left(|Y_t|+\int_0^t \alpha_s {\rm d}s,\ 2q3^{2\lambda-1} \breve \mu_{\beta,\gamma,\lambda,\eps}(t)\right)\right)_{t\in\T}$$
belongs to class (D). \vspace{0.2cm}

Finally, we suppose further that the generator $g$ satisfies Assumption \ref{A:UN3}. The desired uniqueness is a direct consequence of the following \cref{pro:6.6}. The proof is complete.
\end{proof}

\begin{pro}\label{pro:6.6}
Let the functions $\breve \mu_{\beta,\gamma,\lambda,\eps}(\cdot)$ and $\breve \psi(\cdot,\cdot)$ be respectively defined in \eqref{eq:6.13} and \eqref{eq:6.17}, $\xi$ and $\xi'$ be two terminal conditions, $g$ and $g'$ be two generators which are continuous in the variables $(y,z)$, and $(Y_t, Z_t)_{t\in\T}$ and $(Y'_t, Z'_t)_{t\in\T}$ be respectively a solution to BSDE$(\xi, g)$ and BSDE$(\xi', g')$ such that for some $\eps>0$ and $q>1$, both
$$\left(\breve\psi\left(|Y_t|+\int_0^t \alpha_s {\rm d}s,\ 2q3^{2\lambda-1} \breve \mu_{\beta,\gamma,\lambda,\eps}(t)\right)\right)_{t\in\T}$$
and
$$\left(\breve\psi\left(|Y'_t|+\int_0^t \alpha_s {\rm d}s,\ 2q3^{2\lambda-1} \breve \mu_{\beta,\gamma,\lambda,\eps}(t)\right)\right)_{t\in\T}\vspace{0.1cm}$$
belong to class (D). Assume that $\ps$, $\xi\leq \xi'$. If $g$ (resp. $g'$) satisfies Assumptions~\ref{A:UN3} and~\ref{A:A4} with parameters $(\alpha_\cdot,\beta,\gamma,\lambda)$, and $\as$,
\begin{equation}\label{eq:6.25}
g(t,Y'_t,Z'_t)\leq g'(t,Y'_t,Z'_t)\ \ \ ({\rm resp.}\  \ g(t,Y_t,Z_t)\leq g'(t,Y_t,Z_t)\ ),
\end{equation}
then $\ps$, for each $t\in\T$, $Y_t\leq Y'_t$.
\end{pro}

\begin{proof}
We first consider the case that the generator $g$ is convex in the state variables $(y,z)$, satisfies Assumption \ref{A:A4} with parameters $(\alpha_\cdot,\beta,\gamma,\lambda)$, and $\as$,
$$g(t,Y'_t,Z'_t)\leq g'(t,Y'_t,Z'_t).$$

In order to utilize the convexity condition of the generator $g$, we use the $\theta$-technique developed in for example \cite{BriandHu2008PTRF}. For each fixed $\theta\in (0,1)$, define
\begin{equation}\label{eq:6.26}
\Delta^\theta U_\cdot:=\frac{Y_\cdot-\theta Y'_\cdot}{1-\theta}\ \  {\rm and} \ \ \Delta^\theta V_\cdot:=\frac{Z_\cdot-\theta Z'_\cdot}{1-\theta}.
\end{equation}
Then the pair $(\Delta^\theta U_\cdot,\Delta^\theta V_\cdot)$ satisfies the following BSDE:\vspace{0.1cm}
\begin{equation}\label{eq:6.27}
  \Delta^\theta U_t=\Delta^\theta U_T +\int_t^T \Delta^\theta g (s,\Delta^\theta U_s,\Delta^\theta V_s) {\rm d}s-\int_t^T \Delta^\theta V_s \cdot {\rm d}B_s, \ \ \ \ t\in\T,\vspace{0.1cm}
\end{equation}
where $\ass$, for each $(y,z)\in \R\times\R^d$,
\begin{equation}\label{eq:6.28}
\begin{array}{lll}
\Dis \Delta^\theta g(s,y,z)&:=& \Dis \frac{1}{1-\theta}\left[\  g(s,(1-\theta)y+\theta Y'_s,(1-\theta)z+\theta Z'_s)-\theta g(s, Y'_s, Z'_s)\ \right]\vspace{0.2cm}\\
&& \Dis +\frac{\theta}{1-\theta}\left[\  g(s,Y'_s, Z'_s)-g'(s,Y'_s, Z'_s)\ \right].
\end{array}
\end{equation}
It follows from the assumptions that $\ass$,
\begin{equation}\label{eq:6.29}
\RE\ (y,z)\in \R_+\times \R^d,\ \ \  \Delta^\theta g(s,y,z)\leq g(s,y,z)\leq \alpha_s+\beta |y|\ln |y|{\bf 1}_{|y|>1}+\gamma |z||\ln |z||^\lambda,
\end{equation}
which means that the generator $\Delta^\theta g$ satisfies Assumption \ref{A:A4'} defined in \cref{rmk:6.6}.

On the other hand, note that for each $x,y,z\geq 0$, $c>0$ and $\lambda>1/2$, we have
$$1+cx+cy\leq (1+c)(1+x)(1+y)$$
and
$$(x+y+z)^{2\lambda}\leq 3^{2\lambda-1}\left(x^{2\lambda}+y^{2\lambda}+z^{2\lambda}\right).$$
It follows that for each $x,y\geq 0$, $\theta\in (0,1)$ and $\lambda>1/2$,
$$
\left[\ln\left(1+\frac{(x-\theta y)^+}{1-\theta}\right)\right]^{2\lambda}\leq 3^{2\lambda-1}
\left[\left(\ln\left(1+\frac{1}{1-\theta}\right)\right)^{2\lambda}
+\left(\ln\left(1+x\right)\right)^{2\lambda}
+\left(\ln\left(1+y\right)\right)^{2\lambda}\right].
$$
Furthermore, note that for each $\mu>0$, $\lambda>1/2$, $p>1$ and $x>0$, we have
$$
\exp\left(\mu (\ln (1+x))^{2\lambda}\right)\geq x^p.$$
By virtue of the last two inequalities, Young's inequality, the definition of $\breve\psi$ and the assumptions of \cref{pro:6.6}, it is not very hard to verify
that the process
\[
\breve\psi\left(\left(\Delta^\theta U_t\right)^+ +\int_0^t\alpha_s {\rm d}s,\  \breve\mu_{\beta,\gamma,\lambda,\eps}(t)\right),\ \ \ t\in\T
\]
belongs to class (D) for each $\theta\in (0,1)$.
Thus, for BSDE \eqref{eq:6.27}, by virtue of \eqref{eq:6.29}, \cref{rmk:6.6} and the proof of \cref{pro:6.4}, we derive that there exists a $C>0$ depending on $(\beta,\gamma,\lambda,\eps, T)$ such that
\begin{equation}\label{eq:6.30}
\begin{array}{lll}
\left(\Delta^\theta U_t\right)^+&\leq& \Dis \breve\psi\left(\left(\Delta^\theta U_t\right)^+,\ \breve\mu_{\beta,\gamma,\lambda,\eps}(t)\right)\vspace{0.2cm}\\
&\leq& \Dis C\E\left[\left.\breve\psi\left(\left(\Delta^\theta U_T\right)^+ +\int_0^T\alpha_s {\rm d}s,\  \breve\mu_{\beta,\gamma,\lambda,\eps}(T)\right)\right|\F_t\right]+C, \ \ \ t\in\T.
\end{array}
\end{equation}

Furthermore, in view of the fact that
\begin{equation}\label{eq:6.31}
\left(\Delta^\theta U_T\right)^+=\frac{(\xi-\theta \xi')^+}{1-\theta}=\frac{\left[\xi-\theta \xi+\theta(\xi-\xi')\right]^+}{1-\theta}\leq \xi^+,
\end{equation}
it follows from \eqref{eq:6.30} that for each $\theta\in (0,1)$,
$$
\left(Y_t-\theta Y'_t\right)^+ \leq (1-\theta)
\left(C\E\left[\left.\breve\psi\left(\xi^+ +\int_0^T\alpha_s {\rm d}s,\  \breve\mu_{\beta,\gamma,\lambda,\eps}(T)\right)\right|\F_t\right]+C
\right), \ \ \ t\in\T.
$$
Thus, the desired conclusion follows by sending $\theta\To 1$ in above inequality. \vspace{0.1cm}

For the case that the generator $g$ is concave in the state variables $(y,z)$, we need to respectively use the $\theta Y_\cdot-Y'_\cdot$ and $\theta Z_\cdot-Z'_\cdot$ to replace $Y_\cdot-\theta Y'_\cdot$ and $Z_\cdot-\theta Z'_\cdot$ in \eqref{eq:6.26} . In this case, the generator $\Delta^\theta g$ in \eqref{eq:6.28} should be replaced with
\[
\begin{array}{lll}
\Dis \Delta^\theta g(s,y,z)&:=& \Dis \frac{1}{1-\theta}\left[\  \theta g(s, Y_s, Z_s) -g(s,-(1-\theta)y+\theta Y_s,-(1-\theta)z+\theta Z_s)\ \right]\vspace{0.2cm}\\
&& \Dis +\frac{1}{1-\theta}\left[\  g(s,Y'_s, Z'_s)-g'(s,Y'_s, Z'_s)\ \right].
\end{array}
\]
Since $g$ is concave in $(y,z)$, we have $\ass$,
$$
\RE\ (y,z)\in \R\times \R^d,\ \ g(s,-(1-\theta)y+\theta Y_s,-(1-\theta)z+\theta Z_s)\geq \theta g(s,Y_s,Z_s)+(1-\theta)g(t,-y,-z),
$$
and then, \eqref{eq:6.29} can be replaced by
\[
\RE\ (y,z)\in \R_+\times \R^d,\ \ \ \Delta^\theta g(s,y,z)\leq -g(s,-y,-z)\leq \alpha_s+\beta |y|\ln |y|{\bf 1}_{|y|>1}+\gamma |z||\ln |z||^\lambda,
\]
which means that the generator $\Delta^\theta g$ still satisfies Assumption \ref{A:A4'}. Consequently, \eqref{eq:6.30} still holds. Moreover, we use
\[
\left(\Delta^\theta U_T\right)^+=\frac{(\theta \xi-\xi')^+}{1-\theta}=\frac{\left[\theta \xi- \xi+(\xi-\xi')\right]^+}{1-\theta}\leq (-\xi)^+=\xi^-,\vspace{0.1cm}
\]
instead of \eqref{eq:6.31}, and it follows from \eqref{eq:6.30} that for each $\theta\in (0,1)$,
$$
\left(\theta Y_t-Y'_t\right)^+ \leq (1-\theta)
\left(C\E\left[\left.\breve\psi\left(\xi^- +\int_0^T\alpha_s {\rm d}s,\  \breve\mu_{\beta,\gamma,\lambda,\eps}(T)\right)\right|\F_t\right]+C
\right), \ \ \ t\in\T.
$$
Thus, the desired conclusion follows by sending $\theta\To 1$ in above inequality.\vspace{0.1cm}

Finally, in the same way as above, we can prove the desired conclusion under assumptions that the generator $g'$ satisfies \ref{A:UN3} and \ref{A:A4} with parameters $(\alpha_\cdot, \beta, \gamma)$, and $\as$, $g(t,Y_t,Z_t)\leq g'(t,Y_t,Z_t)$. The proof is then complete.\vspace{0.3cm}
\end{proof}

\begin{ex}\label{exp:6.7}
Let $\lambda>1/2$, $\beta>0$, $\gamma>0$ and $k\geq 0$. For $(\omega,t,y,z)\in \Omega\times\T\times\R\times\R^d$, define
$$
g(\omega,t,y,z):=|B_t(\omega)|+ky^2{\bf 1}_{y\leq 0}+\beta |y|\ln |y|{\bf 1}_{|y|> 1}+\gamma |z|(\ln |z|)^{\lambda}{\bf 1}_{|z|>1}.
$$
It is easy to verify that this generator $g$ satisfies Assumptions \ref{A:EX1}, \ref{A:A4} and \ref{A:UN3}.
\end{ex}

\section{Conclusion}
\label{sec:7-mainresult--7}
\setcounter{equation}{0}

The objective of this section is to unify the existence and uniqueness results obtained in previous four sections, namely, Theorems \ref{thm:3.3}, \ref{thm:4.3}, \ref{thm:5.3} and \ref{thm:6.2}.
First, we unify Assumptions \ref{A:A1}, \ref{A:A2}, \ref{A:A3} and \ref{A:A4} into the following one:
\begin{enumerate}
\renewcommand{\theenumi}{(EX2)}
\renewcommand{\labelenumi}{\theenumi}
\item \label{A:EX2} There exist three constants $\beta\geq 0$, $\gamma > 0$ and $\lambda\in [0,+\infty)$, and an $\R_+$-valued progressively measurable process $(\alpha_t)_{t\in\T}$  such that $\as$,
    $$
    \RE\ (y,z)\in \R\times\R^d,\ \ \
    {\rm sgn}(y)g(\omega,t,y,z)\leq \alpha_t(\omega)+\beta|y|(\ln |y|)^{\left(\lambda+1/2\right)\wedge 1}{\bf 1}_{|y|>1}+\gamma |z||\ln |z||^{\lambda}.
    $$
\end{enumerate}

For each $\lambda\in [0,+\infty)$, we define the function
\begin{equation}\label{eq:7.1}
\psi(x,\mu;\lambda):=x\exp\left(\mu\left(\ln (1+x)\right)^{\left(\lambda+{1\over 2}\right)\vee (2\lambda)}\right),\ \ (x,\mu)\in \R_+\times \R_+.
\end{equation}
It can be easily verified that
\begin{itemize}
\item[(i)] if $\lambda\in [0,1/2)$, then for each $p>1$,
$$ x\ln(1+x)< \psi(x,\mu;\lambda)=x\exp\left(\mu\left(\ln (1+x)\right)^{\lambda+{1\over 2}}\right)< x^p,\ \ (x,\mu)\in \R_+\times \R_+;\vspace{-0.3cm}$$
\item[(ii)] if $\lambda=1/2$, then
$$ x^{1+\mu}< \psi(x,\mu;\lambda)=x(1+x)^\mu< (1+x)^{1+\mu},\ \ (x,\mu)\in \R_+\times \R_+;\vspace{-0.3cm}$$
\item[(iii)]if $\lambda\in (1/2,+\infty)$, then for each $p>1$ and $\epsilon\in (0,1)$,
$$ x^p< \psi(x,\mu;\lambda)=x\exp\left(\mu\left(\ln (1+x)\right)^{2\lambda}\right)< \exp(x^\epsilon),\ \ (x,\mu)\in \R_+\times \R_+.\vspace{-0.2cm}$$
\end{itemize}
Moreover, for each $\beta\geq 0$, $\gamma > 0$, and $\eps>0$, let $\mu_{\beta,\gamma,\lambda,\eps}(\cdot)$ with $\lambda>0$ be the unique solution of the following ODE: $\mu(0)=\eps$ and \vspace{0.1cm}
\begin{equation}\label{eq:7.2}
\mu'(s)=
\left\{
\begin{array}{ll}
\eps \mu_s+
\Dis \frac{\gamma^2(1+\eps)^{2\lambda+2}}{2\lambda+1}\frac{1}{\mu_s}+\beta+\eps ,& \Dis \lambda\in (0,1/2);\vspace{0.3cm}\\
\beta \mu_s+
\Dis \frac{\gamma^2(1+\eps)}{2}\left(1+\frac{1}{\mu_s}\right)+\beta ,& \Dis \lambda=1/2;\vspace{0.3cm}\\
2\beta\lambda \mu_s+
\Dis \frac{\gamma^2(1+\eps)k_\lambda}{4\lambda}\left(1+
\frac{1}{\mu_s}\right)+\eps ,& \Dis \lambda\in (1/2,+\infty)
\end{array} \quad
 \hbox{ \rm for } s\in\T, \right.
\end{equation}
and $\mu_{\beta,\gamma,\lambda}^0(\cdot)$ with $\lambda\geq 0$ be the unique solution of the following ODE: $\mu(0)=0$ and \vspace{0.1cm}
\begin{equation}\label{eq:7.3}
\mu'(s)=
\left\{
\begin{array}{ll}
\Dis \frac{\gamma^2}{2\lambda+1}\frac{1}{\mu_s}+\beta ,& \Dis \lambda\in [0,1/2);\vspace{0.3cm}\\
\Dis \beta\mu_s+\frac{\gamma^2}{2}\left(1+\frac{1}{\mu_s}\right)+\beta ,& \Dis \lambda=1/2;\vspace{0.3cm}\\
\Dis 2\beta\lambda \mu_s+\frac{\gamma^2 k_\lambda}{4\lambda}
\left(1+\frac{1}{\mu_s}\right) ,& \Dis \lambda\in (1/2,+\infty)
\end{array} \quad
\hbox{ \rm for }  s\in (0,T]\right.
\end{equation}
with
$$k_\lambda:=2^{2(\lambda-1)^+ +2\lambda-1}.$$
It is not difficult to check that as $\eps\To 0^+$, $\mu_{\beta,\gamma,\lambda,\eps}(\cdot)$ tends decreasingly to $\mu_{\beta,\gamma,\lambda}^0(\cdot)$ on $\T$. \vspace{0.3cm}

In view of \cref{rmk:4.1}, Theorems \ref{thm:3.3}, \ref{thm:4.3}, \ref{thm:5.3} and \ref{thm:6.2} can be unified into the following one.

\begin{thm}\label{thm:7.1}
Let the functions $\psi(x,\mu;\lambda)$, $\mu_{\beta,\gamma,\lambda,\eps}(\cdot)$ and $\mu_{\beta,\gamma,\lambda}^0(\cdot)$ be respectively defined in \eqref{eq:7.1}, \eqref{eq:7.2} and \eqref{eq:7.3}, $\xi$ be a terminal condition and $g$ be a generator which is continuous in the state variables $(y,z)$ and satisfies Assumptions \ref{A:EX1} and \ref{A:EX2} with parameters $(\alpha_\cdot, \beta, \gamma)$. We have\vspace{0.2cm}

(i) If $\lambda=0$ and the terminal condition satisfies
\begin{equation}\label{eq:7.4}
\psi\left(|\xi|+\int_0^T \alpha_t {\rm d}t,\ \mu_{\beta,\gamma,\lambda}^0(T);\ \lambda\right)\in L^1,
\end{equation}
then BSDE $(\xi,g)$ admits a solution $(Y_t,Z_t)_{t\in\T}$ such that $\psi\left(|Y_t|,\mu_{\beta,\gamma,\lambda}^0(t);\lambda\right)$
belongs to class (D). \vspace{0.1cm}

Moreover, if $g$ also satisfies either Assumptions \ref{A:UN1} and \ref{A:UN2} or Assumption \ref{A:UN3}, then the solution satisfying the preceding condition is unique.\vspace{0.2cm}

(ii) If $\lambda=(0,1/2)$ and there exists a constant $\mu> \mu_{\beta,\gamma,\lambda}^0(T)$ such that
\begin{equation}\label{eq:7.5}
\psi\left(|\xi|+\int_0^T \alpha_t {\rm d}t,\ \mu;\ \lambda\right)\in L^1,
\end{equation}
then BSDE$(\xi,g)$ admits a solution $(Y_t,Z_t)_{t\in\T}$ such that for some $\eps>0$ satisfying $\mu_{\beta,\gamma,\lambda,\eps}(T)=\mu$,
$\psi\left(|Y_t|,\mu_{\beta,\gamma,\lambda,\eps}(t);\lambda\right)$
belongs to class (D).  \vspace{0.1cm}

Moreover, if $g$ also satisfies either Assumptions \ref{A:UN1} and \ref{A:UN2} or Assumption \ref{A:UN3}, then the solution satisfying the preceding condition is unique.\vspace{0.2cm}

(iii) If $\lambda=1/2$ and there exists a constant $\mu> \mu_{\beta,\gamma,\lambda}^0(T)$ such that the condition \eqref{eq:7.5} holds, then BSDE$(\xi,g)$ admits a solution $(Y_t,Z_t)_{t\in\T}$ such that for some $\eps>0$ satisfying $\mu_{\beta,\gamma,\lambda,\eps}(T)=\mu$,
$\psi\left(|Y_t|,\mu_{\beta,\gamma,\lambda,\eps}(t);\lambda\right)$
belongs to class (D) and $Z_\cdot\in \mcal^p$ for some $p>1$.  \vspace{0.1cm}

Moreover, if $g=g_1+g_2$, $g_1$ satisfies Assumptions \ref{A:UN1} and \ref{A:UN2} and $g_2$ satisfies Assumption \ref{A:UN3}, then the solution satisfying the preceding condition is unique.\vspace{0.2cm}

(iv) If $\lambda>1/2$ and there exists a constant $\mu> \mu_{\beta,\gamma,\lambda}^0(T)$ such that the condition \eqref{eq:7.5} holds, then BSDE $(\xi,g)$ admits a solution $(Y_t,Z_t)_{t\in\T}$ such that for some $\eps>0$ satisfying $\mu_{\beta,\gamma,\lambda,\eps}(T)=\mu$,
$\psi\left(|Y_t|+\int_0^t \alpha_s{\rm d}s,\mu_{\beta,\gamma,\lambda,\eps}(t);\lambda\right)$
belongs to class (D) and $Z_\cdot\in \mcal^2$. \vspace{0.2cm}

Moreover, if the condition \eqref{eq:7.5} holds for some $\mu>2q3^{2\lambda-1} \mu_{\beta,\gamma,\lambda}^0(T)$ with $q>1$, then BSDE$(\xi,g)$ admits a solution $(Y_t,Z_t)_{t\in\T}$ such that the process
$$\left(\breve\psi\left(|Y_t|+\int_0^t \alpha_s {\rm d}s,\ 2q3^{2\lambda-1} \mu_{\beta,\gamma,\lambda,\eps}(t);\lambda\right)\right)_{t\in\T}$$
belongs to class (D), where $\eps$ is the unique positive constant satisfying $2q3^{2\lambda-1} \mu_{\beta,\gamma,\lambda,\eps}(T)=\mu$. And, if $g$ further satisfies Assumption \ref{A:UN3}, then the solution satisfying the preceding condition is also unique.\vspace{0.2cm}
\end{thm}

Finally, we would like to mention that the integrability conditions~\eqref{eq:7.4} and \eqref{eq:7.5} are believed to be reasonably weakest possible for the existence of the solution, but it remains to be proved. To the best of our knowledge, all conclusions in \cref{thm:7.1} can not be obtained from any existing results.

\vspace{0.4cm}


\noindent {\bf References}\vspace{0.2cm}

\setlength{\bibsep}{2pt}

\begin{thebibliography}{29}
\expandafter\ifx\csname natexlab\endcsname\relax\def\natexlab#1{#1}\fi
\expandafter\ifx\csname url\endcsname\relax
  \def\url#1{\texttt{#1}}\fi
\expandafter\ifx\csname urlprefix\endcsname\relax\def\urlprefix{URL }\fi

\bibitem[{Bahlali(2019)}]{Bahlali2019Arxiv}
Bahlali K., 2019. Solving unbounded quadratic {BSDE}s by a domination method.
  \href{https://arxiv.org/pdf/1903.11325v1.pdf}{arXiv:1903.11325v1 [math.PR] 27
  Mar 2019}.

\bibitem[{Bahlali et~al.(2017)Bahlali, Eddahbi, and
  Ouknine}]{BahlaliEddahbiOuknine2017AoP}
Bahlali K., Eddahbi M., Ouknine Y., 2017. Quadratic {BSDE}
  with ${L}^2$-terminal data: {K}rylov's estimate, {I}t\^{o}-{K}rylov's formula
  and existence results. Ann. Probab. 45~(4), 2377--2397.

\bibitem[{Bahlali and El~Asri(2012)}]{BahlaliElAsri2012BSM}
Bahlali K., El~Asri B., 2012. Stochastic optimal control and {BSDE}s with
  logarithmic growth. Bull. Sci. Math. 136~(6), 617--637.

\bibitem[{Bahlali et~al.(2010)Bahlali, Essaky, and
  Hassani}]{BahlaliEssakyHassani2010CRM}
Bahlali K., Essaky E.~H., Hassani M., 2010. Multidimensional {BSDE}s with
  superlinear growth coefficient. {A}pplication to degenerate systems of
  semilinear {PDE}s. C. R. Math. Acad. Sci. Paris 348~(11-12), 677--682.

\bibitem[{Bahlali et~al.(2015)Bahlali, Essaky, and
  Hassani}]{BahlaliEssakyHassani2015SIAM}
Bahlali K., Essaky E.~H., Hassani M., 2015. Existence and
  uniqueness of multidimensional {BSDE}s and of systems of degenerate {PDE}s
  with superlinear growth generator. SIAM J. Math. Anal. 47~(6), 4251--4288.

\bibitem[{Bahlali et~al.(2015)Bahlali, Hakassou, and
  Ouknine}]{BahlaliHakassouOuknine2015Stochastics}
Bahlali K., Hakassou A., Ouknine Y., 2015. A class of
  stochastic differential equations with superlinear growth and non-{L}ipschitz
  coefficients. Stochastics 87~(5), 806--847.

\bibitem[{Bahlali et~al.(2017)Bahlali, Kebiri, Khelfallah, and
  Moussaoui}]{BahlaliKebiri2017Stochastics}
Bahlali K., Kebiri O., Khelfallah N., Moussaoui H., 2017. One
  dimensional {BSDE}s with logarithmic growth application to {PDE}s.
  Stochastics 89~(6-7), 1061--1081.

\bibitem[{Bahlali and Tangpi(2019)}]{BahlaliTangpi2019Arxiv}
Bahlali K., Tangpi L., 2019. {BSDE}s driven by $|z|^2/y$ and applications.
  \href{https://arxiv.org/pdf/1810.05664v2.pdf}{arXiv:1810.05664v2 [math.PR] 1
  Feb 2019}.

\bibitem[{Barrieu and El~Karoui(2013)}]{BarrieuElKaroui2013AoP}
Barrieu P., El~Karoui N., 2013. Monotone stability of quadratic
  semimartingales with applications to unbounded general quadratic {BSDE}s.
  Ann. Probab. 41~(3B), 1831--1863.

\bibitem[{Briand et~al.(2003)Briand, Delyon, Hu, Pardoux, and
  Stoica}]{BriandDelyonHu2003SPA}
Briand P., Delyon B., Hu Y., Pardoux E., Stoica L., 2003. {$L^p$}
  solutions of backward stochastic differential equations. Stochastic Process.
  Appl. 108~(1), 109--129.

\bibitem[{Briand and Elie(2013)}]{BriandElie2013SPA}
Briand P., Elie R., 2013. A simple constructive approach to quadratic {BSDE}s
  with or without delay. Stochastic Process. Appl. 123, 2921--2939.

\bibitem[{Briand and Hu(2006)}]{BriandHu2006PTRF}
Briand P., Hu Y., 2006. {BSDE} with quadratic growth and unbounded terminal
  value. Probab. Theory Related Fields 136~(4), 604--618.

\bibitem[{Briand and Hu(2008)}]{BriandHu2008PTRF}
Briand P., Hu Y., 2008. {Q}uadratic {BSDE}s with convex generators and
  unbounded terminal conditions. Probab. Theory Related Fields 141~(3),
  543--567.

\bibitem[{Buckdahn et~al.(2018)Buckdahn, Hu, and Tang}]{BuckdahnHuTang2018ECP}
Buckdahn R., Hu Y., Tang S., 2018. Uniqueness of solution to scalar {BSDE}s
  with ${L}\exp\left(\mu\sqrt{2\log(1+L)}\right)$-integrable terminal values.
  Electron. Commun. Probab. 23, Paper No. 59, 8pp.

\bibitem[{Delbaen et~al.(2011)Delbaen, Hu, and
  Richou}]{DelbaenHuRichou2011AIHPPS}
Delbaen F., Hu Y., Richou A., 2011. On the uniqueness of solutions to
  quadratic {BSDE}s with convex generators and unbounded terminal conditions.
  Ann. Inst. Henri Poincar\'{e} Probab. Stat. 47, 559--574.

\bibitem[{El~Karoui et~al.(1997)El~Karoui, Peng, and
  Quenez}]{ElKarouiPengQuenez1997MF}
El~Karoui N., Peng S., Quenez M.~C., 1997. Backward stochastic differential
  equations in finance. Math. Finance 7~(1), 1--71.

\bibitem[{Fan(2015)}]{Fan2015JMAA}
Fan S., 2015. ${L}^p$ solutions of multidimensional {BSDE}s with weak
  monotonicity and general growth generators. J. Math. Anal. Appl. 432,
  156--178.

\bibitem[{Fan(2016)}]{Fan2016SPA}
Fan S., 2016. Bounded solutions, ${L}^{p}\ (p>1)$ solutions and ${L}^1$
  solutions for one-dimensional {BSDE}s under general assumptions. Stochastic
  Process. Appl. 126, 1511--1552.

\bibitem[{Fan and Hu(2019)}]{FanHu2019ECP}
Fan S., Hu Y., 2019. Existence and uniqueness of solution to scalar {BSDE}s
  with ${L}\exp\left(\mu\sqrt{2\log(1+L)}\right)$-integrable terminal values:
  the critical case. Electron. Commun. Probab. 24, Paper No. 49, 10pp.

\bibitem[{Fan and Hu(2021)}]{FanHu2021SPA}
Fan S., Hu Y., 2021. Well-posedness of scalar BSDEs with sub-quadratic generators and related PDEs. Stochastic Process. Appl., 131, 21--50.

 \bibitem[{Fan et~al.(2020)Fan, Hu, and Tang}]{FanHuTang2020CRM}
Fan S., Hu Y., Tang S., 2020. On the uniqueness of solutions to quadratic
  {BSDE}s with non-convex generators and unbounded terminal conditions. C. R. Math. Acad. Sci. Paris  358~(2),  227--235.

\bibitem[{Fan and Jiang(2012)}]{FanJiang2012JAMC}
Fan S., Jiang L., 2012. ${L}^p$ $(p>1)$ solutions for one-dimensional {BSDE}s
  with linear-growth generators. Journal of Applied Mathematics and Computing
  38~(1--2), 295--304.

\bibitem[{Fan and Jiang(2013)}]{FanJiang2013AMSE}
Fan S., Jiang L., 2013. Multidimensional {BSDE}s with weak monotonicity and
  general growth generators. Acta Mathematica Sinica, English Series 29~(10),
  1885--1907.

\bibitem[{Fan et~al.(2010)Fan, Jiang, and Davison}]{FanJiangDavison2010CRA}
Fan S., Jiang L., Davison M., 2010. Uniqueness of solutions for
  multidimensional {BSDE}s with uniformly continuous generators. C. R. Math.
  Acad. Sci. Paris 348~(11--12), 683--687.

\bibitem[{Hu and Tang(2015)}]{HuTang2015SPA}
Hu Y., Tang S., 2015. Multi-dimensional backward stochastic differential
  equations of diagonally quadratic generators. Stochastic Process. Appl.
  126~(4), 1066--1087.

\bibitem[{Hu and Tang(2018)}]{HuTang2018ECP}
Hu Y., Tang S., 2018. Existence of solution to scalar {BSDE}s with
  ${L}\exp\sqrt{{2\over \lambda}\log(1+L)}$-integrable terminal values.
  Electron. Commun. Probab. 23, Paper No. 27, 11pp.

\bibitem[{Kobylanski(2000)}]{Kobylanski2000AP}
Kobylanski M., 2000. Backward stochastic differential equations and partial
  differential equations with quadratic growth. Ann. Probab. 28~(2), 558--602.

\bibitem[{Lepeltier and San~Martin(1997)}]{LepeltierSanMartin1997SPL}
Lepeltier J.-P., San~Martin J., 1997. Backward stochastic differential
  equations with continuous coefficient. Statist. Probab. Lett. 32~(4),
  425--430.

\bibitem[{Luo and Fan(2018)}]{LuoFan2018SD}
Luo H., Fan S., 2018. Bounded solutions for general time interval {BSDE}s
  with quadratic growth coefficients and stochastic conditions. Stoch. Dynam.
  18~(5), Paper No. 1850034, 24pp.

\bibitem[{Pardoux and Peng(1990)}]{PardouxPeng1990SCL}
Pardoux E., Peng S., 1990. Adapted solution of a backward stochastic
  differential equation. Syst. Control Lett. 14~(1), 55--61.

\bibitem[{Yang(2017)}]{Yanghanlin2017Arxiv}
Yang H., 2017. ${L}^p$ solutions of quadratic {BSDE}s.
  \href{https://arxiv.org/pdf/1506.08146v2.pdf}{arXiv:1506.08146v2 [math.PR] 29
  Sep 2017}.

\end{thebibliography}

\end{document}